\documentclass[11pt, a4paper]{article}
\usepackage[cp1251]{inputenc}
\usepackage{amsmath} \usepackage{euscript} \usepackage{marvosym}
\usepackage[dvipsnames]{xcolor}
\usepackage{mymatrx6}
\usepackage{graphicx}
\usepackage{longtable}
\usepackage{mathrsfs}
\usepackage{amssymb, amscd}
\usepackage{adigraph}
\def\mymatrix{\MyMatrixwithdelims..}
\Linestrue\Autonumtrue
\usepackage{mathrsfs}
\def\bluecommand{\textcolor[cmyk]{1,0,0,0}}
\def\redcommand{\textcolor[cmyk]{0,1,1,0}}
\def\greencommand{\textcolor[cmyk]{1,0,1,0}}
\def\yellowcommand{\textcolor[cmyk]{0,0.25,1,0}}
\def\graycommand{\textcolor[cmyk]{0,0,0,0.2}}

\newtheorem{theorem}{Theorem.}
\newtheorem{lemma}[theorem]{Lemma.}

\begin{document}

\newcounter{bnomer} \newcounter{snomer}
\newcounter{bsnomer}
\setcounter{bnomer}{0}
\renewcommand{\thesnomer}{\thebnomer.\arabic{snomer}}
\renewcommand{\thebsnomer}{\thebnomer.\arabic{bsnomer}}
\renewcommand{\refname}{\begin{center}\large{\textbf{References}}\end{center}}

\setcounter{MaxMatrixCols}{14}

\newcommand\restr[2]{{
  \left.\kern-\nulldelimiterspace 
  #1 
  \right|_{#2} 
}}

\newcommand{\sect}[1]{%
\setcounter{snomer}{0}\setcounter{bsnomer}{0}
\refstepcounter{bnomer}
\par\bigskip\begin{center}\large{\textbf{\arabic{bnomer}. {#1}}}\end{center}}
\newcommand{\sst}[1]{%
\refstepcounter{bsnomer}
\par\bigskip\textbf{\arabic{bnomer}.\arabic{bsnomer}. {#1}}\par}
\newcommand{\defi}[1]{%
\refstepcounter{snomer}
\par\medskip\textbf{Definition \arabic{bnomer}.\arabic{snomer}. }{#1}\par\medskip}
\newcommand{\theo}[2]{%
\refstepcounter{snomer}
\par\textbf{Theorem \arabic{bnomer}.\arabic{snomer}. }{#2} {\emph{#1}}\hspace{\fill}$\square$\par}
\newcommand{\mtheop}[2]{%
\refstepcounter{snomer}
\par\textbf{Theorem \arabic{bnomer}.\arabic{snomer}. }{\emph{#1}}
\par\textsc{Proof}. {#2}\hspace{\fill}$\square$\par}
\newcommand{\mcorop}[2]{%
\refstepcounter{snomer}
\par\textbf{Corollary \arabic{bnomer}.\arabic{snomer}. }{\emph{#1}}
\par\textsc{Proof}. {#2}\hspace{\fill}$\square$\par}
\newcommand{\mtheo}[1]{%
\refstepcounter{snomer}
\par\medskip\textbf{Theorem \arabic{bnomer}.\arabic{snomer}. }{\emph{#1}}\par\medskip}
\newcommand{\theobn}[1]{%
\par\medskip\textbf{Theorem. }{\emph{#1}}\par\medskip}
\newcommand{\theoc}[2]{%
\refstepcounter{snomer}
\par\medskip\textbf{Theorem \arabic{bnomer}.\arabic{snomer}. }{#1} {\emph{#2}}\par\medskip}
\newcommand{\mlemm}[1]{%
\refstepcounter{snomer}
\par\medskip\textbf{Lemma \arabic{bnomer}.\arabic{snomer}. }{\emph{#1}}\par\medskip}
\newcommand{\mprop}[1]{%
\refstepcounter{snomer}
\par\medskip\textbf{Proposition \arabic{bnomer}.\arabic{snomer}. }{\emph{#1}}\par\medskip}
\newcommand{\theobp}[2]{%
\refstepcounter{snomer}
\par\textbf{Theorem \arabic{bnomer}.\arabic{snomer}. }{#2} {\emph{#1}}\par}
\newcommand{\theop}[2]{%
\refstepcounter{snomer}
\par\textbf{Theorem \arabic{bnomer}.\arabic{snomer}. }{\emph{#1}}
\par\textsc{Proof}. {#2}\hspace{\fill}$\square$\par}
\newcommand{\theosp}[2]{%
\refstepcounter{snomer}
\par\textbf{Theorem \arabic{bnomer}.\arabic{snomer}. }{\emph{#1}}
\par\textsc{Sketch of the proof}. {#2}\hspace{\fill}$\square$\par}
\newcommand{\exam}[1]{%
\refstepcounter{snomer}
\par\medskip\textbf{Example \arabic{bnomer}.\arabic{snomer}. }{#1}\par\medskip}
\newcommand{\deno}[1]{%
\refstepcounter{snomer}
\par\textbf{Notation \arabic{bnomer}.\arabic{snomer}. }{#1}\par}
\newcommand{\lemm}[1]{%
\refstepcounter{snomer}
\par\textbf{Lemma \arabic{bnomer}.\arabic{snomer}. }{\emph{#1}}\hspace{\fill}$\square$\par}
\newcommand{\lemmp}[2]{%
\refstepcounter{snomer}
\par\medskip\textbf{Lemma \arabic{bnomer}.\arabic{snomer}. }{\emph{#1}}
\par\textsc{Proof}. {#2}\hspace{\fill}$\square$\par\medskip}
\newcommand{\coro}[1]{%
\refstepcounter{snomer}
\par\textbf{Corollary \arabic{bnomer}.\arabic{snomer}. }{\emph{#1}}\hspace{\fill}$\square$\par}
\newcommand{\mcoro}[1]{%
\refstepcounter{snomer}
\par\textbf{Corollary \arabic{bnomer}.\arabic{snomer}. }{\emph{#1}}\par\medskip}
\newcommand{\corop}[2]{%
\refstepcounter{snomer}
\par\textbf{Corollary \arabic{bnomer}.\arabic{snomer}. }{\emph{#1}}
\par\textsc{Proof}. {#2}\hspace{\fill}$\square$\par}
\newcommand{\nota}[1]{%
\refstepcounter{snomer}
\par\medskip\textbf{Remark \arabic{bnomer}.\arabic{snomer}. }{#1}\par\medskip}
\newcommand{\propp}[2]{%
\refstepcounter{snomer}
\par\medskip\textbf{Proposition \arabic{bnomer}.\arabic{snomer}. }{\emph{#1}}
\par\textsc{Proof}. {#2}\hspace{\fill}$\square$\par\medskip}
\newcommand{\hypo}[1]{%
\refstepcounter{snomer}
\par\medskip\textbf{Conjecture \arabic{bnomer}.\arabic{snomer}. }{\emph{#1}}\par\medskip}
\newcommand{\prop}[1]{%
\refstepcounter{snomer}
\par\textbf{Proposition \arabic{bnomer}.\arabic{snomer}. }{\emph{#1}}\hspace{\fill}$\square$\par}

\newcommand{\proof}[2]{%
\par\medskip\textsc{Proof{#1}}. \hspace{-0.2cm}{#2}\hspace{\fill}$\square$\par\medskip}

\makeatletter
\def\iddots{\mathinner{\mkern1mu\raise\p@
\vbox{\kern7\p@\hbox{.}}\mkern2mu
\raise4\p@\hbox{.}\mkern2mu\raise7\p@\hbox{.}\mkern1mu}}
\makeatother

\newcommand{\okr}[2]{%
\refstepcounter{snomer}
\par\medskip\textbf{{#1} \arabic{bnomer}.\arabic{snomer}. }{\emph{#2}}\par\medskip}

\newcommand{\Ind}[3]{%
\mathrm{Ind}_{#1}^{#2}{#3}}
\newcommand{\Res}[3]{%
\mathrm{Res}_{#1}^{#2}{#3}}
\newcommand{\epsi}{\varepsilon}
\newcommand{\tri}{\triangleleft}
\newcommand{\Supp}[1]{%
\mathrm{Supp}(#1)}
\newcommand{\NSupp}[1]{%
\mathrm{NSupp}(#1)}
\newcommand{\SSu}[1]{%
\mathrm{SingSupp}(#1)}

\newcommand{\lee}{\leqslant}
\newcommand{\gee}{\geqslant}
\newcommand{\reg}{\mathrm{reg}}
\newcommand{\Dyn}{\mathrm{Dyn}}
\newcommand{\Ann}{\mathrm{Ann}\,}
\newcommand{\Cent}[1]{\mathbin\mathrm{Cent}({#1})}
\newcommand{\PCent}[1]{\mathbin\mathrm{PCent}({#1})}
\newcommand{\Irr}[1]{\mathbin\mathrm{Irr}({#1})}
\newcommand{\Exp}[1]{\mathbin\mathrm{Exp}({#1})}
\newcommand{\empr}[2]{[-{#1},{#1}]\times[-{#2},{#2}]}
\newcommand{\sreg}{\mathrm{sreg}}
\newcommand{\ilm}{\varinjlim}
\newcommand{\wdth}{\mathrm{wd}}
\newcommand{\plm}{\varprojlim}
\newcommand{\codim}{\mathrm{codim}\,}
\newcommand{\GKdim}{\mathrm{GKdim}\,}
\newcommand{\chara}{\mathrm{char}\,}
\newcommand{\rk}{\mathrm{rk}\,}
\newcommand{\chr}{\mathrm{ch}\,}
\newcommand{\Ker}{\mathrm{Ker}\,}
\newcommand{\id}{\mathrm{id}}
\newcommand{\Ad}{\mathrm{Ad}}
\newcommand{\Gh}{\mathrm{Gh}}
\newcommand{\col}{\mathrm{col}}
\newcommand{\row}{\mathrm{row}}
\newcommand{\high}{\mathrm{high}}
\newcommand{\low}{\mathrm{low}}
\newcommand{\pho}{\hphantom{\quad}\vphantom{\mid}}
\newcommand{\fho}[1]{\vphantom{\mid}\setbox0\hbox{00}\hbox to \wd0{\hss\ensuremath{#1}\hss}}
\newcommand{\wt}{\widetilde}
\newcommand{\wh}{\widehat}
\newcommand{\ad}[1]{\mathrm{ad}_{#1}}
\newcommand{\tr}{\mathrm{tr}\,}
\newcommand{\GL}{\mathrm{GL}}
\newcommand{\SL}{\mathrm{SL}}
\newcommand{\SO}{\mathrm{SO}}
\newcommand{\Or}{\mathrm{O}}
\newcommand{\Sp}{\mathrm{Sp}}
\newcommand{\SuppD}{\mathbb{S}\mathrm{upp}}
\newcommand{\Sa}{\mathrm{S}}
\newcommand{\Sing}{\mathrm{Sing}}
\newcommand{\Ua}{\mathrm{U}}
\newcommand{\Andre}{\mathrm{Andre}}
\newcommand{\Aord}{\mathrm{Aord}}
\newcommand{\Mat}{\mathrm{Mat}}
\newcommand{\Stab}{\mathrm{Stab}}
\newcommand{\htt}{\mathfrak{h}}
\newcommand{\spt}{\mathfrak{sp}}
\newcommand{\slt}{\mathfrak{sl}}
\newcommand{\sot}{\mathfrak{so}}

\newcommand{\vfi}{\varphi}
\newcommand{\aad}{\mathrm{ad}}
\newcommand{\vpi}{\varpi}
\newcommand{\teta}{\vartheta}
\newcommand{\Bfi}{\Phi}
\newcommand{\Fp}{\mathbb{F}}
\newcommand{\Rp}{\mathbb{R}}
\newcommand{\Zp}{\mathbb{Z}}
\newcommand{\Cp}{\mathbb{C}}
\newcommand{\Ap}{\mathbb{A}}
\newcommand{\Pp}{\mathbb{P}}
\newcommand{\Kp}{\mathbb{K}}
\newcommand{\Np}{\mathbb{N}}
\newcommand{\ut}{\mathfrak{u}}
\newcommand{\at}{\mathfrak{a}}
\newcommand{\glt}{\mathfrak{gl}}
\newcommand{\hei}{\mathfrak{hei}}
\newcommand{\nt}{\mathfrak{n}}
\newcommand{\kt}{\mathfrak{k}}
\newcommand{\mt}{\mathfrak{m}}
\newcommand{\rt}{\mathfrak{r}}
\newcommand{\rad}{\mathfrak{rad}}
\newcommand{\bt}{\mathfrak{b}}
\newcommand{\unt}{\underline{\mathfrak{n}}}
\newcommand{\gt}{\mathfrak{g}}
\newcommand{\vt}{\mathfrak{v}}
\newcommand{\pt}{\mathfrak{p}}
\newcommand{\Xt}{\mathfrak{X}}
\newcommand{\Po}{\mathcal{P}}
\newcommand{\PV}{\mathcal{PV}}
\newcommand{\Uo}{\EuScript{U}}
\newcommand{\Fo}{\EuScript{F}}
\newcommand{\Do}{\EuScript{D}}
\newcommand{\Eo}{\EuScript{E}}
\newcommand{\Jo}{\EuScript{J}}
\newcommand{\Iu}{\mathcal{I}}
\newcommand{\Mo}{\mathcal{M}}
\newcommand{\Nu}{\mathcal{N}}
\newcommand{\Ro}{\mathcal{R}}
\newcommand{\Co}{\mathcal{C}}
\newcommand{\Ko}{\mathcal{K}}
\newcommand{\So}{\mathcal{S}}
\newcommand{\Lo}{\mathcal{L}}
\newcommand{\Ou}{\mathcal{O}}
\newcommand{\Uu}{\mathcal{U}}
\newcommand{\Tu}{\mathcal{T}}
\newcommand{\Au}{\mathcal{A}}
\newcommand{\Vu}{\mathcal{V}}
\newcommand{\Du}{\mathcal{D}}
\newcommand{\Bu}{\mathcal{B}}
\newcommand{\Sy}{\mathcal{Z}}
\newcommand{\Sb}{\mathcal{F}}
\newcommand{\Gr}{\mathcal{G}}
\newcommand{\Xu}{\mathcal{X}}
\newcommand{\Op}{\mathbb{O}}
\newcommand{\chv}{\mathrm{chv}}
\newcommand{\rtc}[1]{C_{#1}^{\mathrm{red}}}

\newcommand{\JSpec}[1]{\mathrm{JSpec}\,{#1}}
\newcommand{\MSpec}[1]{\mathrm{MSpec}\,{#1}}
\newcommand{\PSpec}[1]{\mathrm{PSpec}\,{#1}}
\newcommand{\APbr}[1]{\mathrm{span}\{#1\}}
\newcommand{\APbre}[1]{\langle #1\rangle}
\newcommand{\APro}[1]{\setcounter{AP}{#1}\Roman{AP}}\newcommand{\apro}[1]{{\rm\setcounter{AP}{#1}\roman{AP}}}
\newcommand{\ot}{\xleftarrow[]{}}
\newcounter{AP}


\author{Mikhail Ignatev\and Alexey Petukhov}
\date{}
\title{Coadjoint orbits of low dimension for nilradicals of Borel subalgebras in classical types}\maketitle
\begin{abstract} Let $\gt$ be a classical simple Lie algebra over an algebraically closed field $\Fp$ of characteristic zero or large enough, and let $\nt$ be a maximal nilpotent subalgebra of $\gt$. 
The main tool in representation theory of $\nt$ is the orbit method, which classifies primitive ideals in the universal enveloping algebra $\Ua(\nt)$ and unitary representations of the unipotent group $N=\exp(\nt)$ in terms of coadjoint orbits on the dual space $\nt^*$. In the paper, we describe explicitly coadjoint orbits of low dimension for $\nt$ as above. The answer is given in terms of subsets of positive roots.\break As a corollary, we provide a way to calculate the number of irreducible complex representations of dimensions $q$, $q^2$ and $q^3$ for a maximal unipotent subgroup $N(q)$ in a classical Chevalley group $G(q)$ over a finite field $\Fp_q$ with $q$ elements. It turned out that this number is a polynomial in $q-1$ with nonnegative integer coefficients, which agrees with Isaac's conjecture.

\medskip\noindent{\bf Keywords:} coadjoint orbit, the orbit method, Lie--Dynkin nil-algebra, root system, finite unipotent group, the number of irreducible characters.\\
{\bf AMS subject classification:} 17B08, 17B10, 17B30, 17B65.\end{abstract}

\tableofcontents\newpage

\sect{Introduction}\label{sect:intro}\addcontentsline{toc}{subsection}{\ref{sect:intro}. Introduction}

Representation theory of unipotent Lie or algebraic groups naturally lead to the study of some geometric objects attached to these groups, i.e., coadjoint orbits. 
The main technical tool in this branch of mathematics, the orbit method due to A.A. Kirillov, see~\cite{Kirillov04}, is based on the interplay between symplectic geometry and advanced linear algebra. 
Originally created by A.A. Kirillov in 1962 for unitary representations of nilpotent real and complex Lie groups, it was adopted by D. Kazhdan in 1977 for unipotent algebraic groups over fields of positive characteristic. Furthermore, in 1970's, J.~Dixmier and his school constructed a version of the orbit method for Lie algebras. In this work we study such coadjoint orbits in certain special cases and even count the number of orbits of certain classes under the assumption that ground field is finite.

Let $\nt$ be the Lie algebra of a unipotent group $N$, and $\nt^*$ be the dual space of $\nt$. Clearly, $N$ acts on $\nt$ via the adjoint action, and the dual action of $N$ on $\nt^*$ is called coadjoint. We denote by $\Ua(\nt)$ the universal enveloping algebra of $\nt$. Recall that, by definition, a primitive ideal in $\Ua(\nt)$ is the annihilator of a simple $\nt$-module.
For $\Fp=\Cp$, the orbit method says that there are bijections between the following sets:
\begin{itemize}
\item the continuous unitary simple modules of $N$; 
\item the coadjoint orbits of $N$ in $\nt^*$;
\item the primitive ideals of $\Ua(\nt)$.
\end{itemize}
Furthermore, the bijection between the last two sets is in fact a homeomorphism, where the space of orbits is endowed with the usual Zariski quotient topology, while the space of primitive ideals in $\Ua(\nt)$ is endowed with the Jacobson topology.
For a finite field $\Fp$, the orbit method establishes a bijection between the set of coadjoint orbits and the set of complex finite-dimensional irreducible characters of the group $N$. 


There are several cases for which there is a known description of coadjoint orbits, see, e.g.,~\cite{GoodwinMoschRoehrle16}, \cite{GoodwinLeMagaardPaolini16}, \cite{GoodwinLeMagaard}, \cite{Panov08}, \cite{Sur23},  \cite{Sur25}. 
Unfortunately, a complete description of coadjoint orbits in most cases is a wild problem. Hence, a natural question is how to describe explicitly interesting classes of coadjoint orbits on $\nt^*$. For example, coadjoint orbits of maximal possible dimension for the Lie algebra of upper-triangular matrices were classified in the first Kirillov's work on the orbit method \cite{Kirillov62}. It is interesting that for the root systems $B_n$, $D_n$ and $E_8$ a complete description of orbits of maximal dimension is still unknown. At the same time, B. Kostant established a description of a very nice class of orbits on $\nt^*$, whose dimension is always maximal possible, with $\nt$ being a nilradical of a Borel subalgebra of a simple finite-dimensional Lie algebra~$\gt$ \cite{Kostant12}, \cite{Kostant13}. Orbits of submaximal dimension for the upper-triangular matrices were described in the paper~\cite{IgnatevPanov09}.

On the other hand, a lot of important patterns about orbits in general can be observed via studying orbits of \emph{low} dimensions. 
But the complete description of such orbits was not available before. The goal of this paper is to fill this gap. 
In this paper we provide a complete description of all coadjoint orbits of dimensions 2, 4, 6 for the maximal nilpotent subalgebras of all classical series of simple Lie algebras. 
Note that the dimension of a coadjoint orbit is always even. 

Let $\nt$ be such a Lie algebra. 
As a first step of this classification we reduce the full description to the description of so-called extensive orbits, see Definition~\ref{D:ext} or the text below. 
The point here is that every orbit can be canonically decomposed as a pointwise Minkowski sum of several other orbits standing for nonempty disjoint subpieces of the underlying Dynkin diagram and a character, see Theorem~\ref{T:extdec}. 
An orbit is extensive if exactly one such a piece is involved. 
This construction implies that every orbit can be canonically decomposed as a pointwise direct sum of several extensive orbits and a character, see Theorem~\ref{T:extdec}. 
One of the advantages of this approach is that extensive orbits of every given dimension shows up only in finitely many Lie algebras, see Proposition~\ref{PlWd} for the detail.

Our first main result here is a classification of extensive orbits on~$\nt^*$ of dimension less or equal to 6 for all classical finite-dimensional simple Lie algebras $\gt=\mathrm{Lie}\,G$. This automatically provides the classification of all orbits of dimension less or equal to 6 for all classical finite-dimensional simple Lie algebras $\gt$.

The algebra $\nt$ is spanned by the root vectors $e_{\alpha}$ for $\alpha\in\Phi^+$, where $\Phi$ is the root system of $\gt$ and $\Phi^+$ is the set of positive roots with respect to a fixed Borel subalgebra $\bt=\mathrm{Lie}\,B$. So, it is not surprising that the extensive orbits of low dimensions are classified in terms of subsets of $\Phi^+$. 

Namely, to each linear form $f\in\nt^*$ one can attach its \emph{support} $$\Supp{f}=\{\alpha\in\Phi^+\mid f(e_{\alpha})\neq0\}.$$ 
To a subset $S\subset\Phi^+$ we attach
$$V(S):=\mathrm{Supp}^{-1}(S):=\{f\in\nt^*\mid \Supp{f}=S\}.$$
It is clear that $V(S)$ is a subvariety of $\nt^*$ and $V(S)\cong(\Fp^\times)^{\times |S|}$. 

We provide a finite list of subsets $S_1$, $S_2$, $\ldots\subset \Phi^+$ such that each orbit of dimension $\leq6$ contains exactly one linear form $f\in V(S_1)\sqcup V(S_2)\sqcup\ldots$ and for each $f$ from this union we have $\dim N.f\le 6$. 
This means that we identify orbits with certain collection of labeled subsets.  
The precise statements are given in Theorems~\ref{Tdim2},~\ref{Tdim4},~\ref{Tdim6}.

As a corollary, we compute the number of irreducible complex representations of low dimension for a maximal unipotent subgroup $N(q)$ of a classical Chevalley group $G(q)$ over a finite field $\Fp_q$ with $q$ elements, where the characteristic of~$\Fp_q$ is large enough. According to the orbit method, the dimension of an irreducible character equals $q^e$, where $2e$ is the dimension of the corresponding orbit. Thus, to finish the computation it remains to count the number of subsets $S$ and make a polynomial out of this data by a simple rule, see Theorem~\ref{theo:Isaacs} for the detail. We consider this theorem as our second main result.

It turned out that the number of irreducible characters of degree $q$, $q^2$ or $q^3$ is a polynomial\break in $q-1$ with nonnegative integer coefficients. In 1960, G. Higman conjectured \cite{Higman60} that the number $O_n(q)$ of the irreducible characters of the unitriangular group $U_n(q)$ (i.e., the group of upper-triangular $n\times n$ matrices over $\Fp_q$ with units on the diagonal) is a polynomial in $q$ for all possible $e$. G. Lehrer conjectured \cite{Lehrer74} that the number $O_{n,e}(q)$ of irreducible characters of $U_n(q)$ of degree $q^e$ is a polynomial in $q$ for all possible $e$. More recently, I.M. Isaacs \cite{Isaacs07} put forth the even stronger conjecture that $O_{n,e}(q)$ is a polynomial in $q-1$ with nonnegative integer coefficients (see Section~\ref{sect:number_chars} for further links and detailed discussion of this conjecture). Hence, our results imply that the Isaac's conjecture is true for $e=1,2,3$, as well as its generalizations for other classical Chevalley groups. (For $U_n(q)$, this was proved in \cite{Loukaki11}; see also \cite{GoodwinMoschRoehrle16} and Section~\ref{sect:number_chars} for the detail.)

To compute the respective sets $S_1$, $S_2$, $\ldots$ we use a method which was originally developed for general finite groups and is used to classify all representations of $\nt$ for root systems of small rank, see~\cite{Evseev10} for the detail. 
The concept of extensive orbits allows us to reduce the problem to the root systems of rank 8 or less. 
We reworked the method of~\cite{GoodwinLeMagaardPaolini16} slightly for such groups: it turns out that the method can be encoded via placement of letters `S', `A', `I', `L' in the boxes of some tables corresponding to the classical root systems, see Section~\ref{S:dim2} for the detail. 
We wish to mention that a major part of the respective calculation was done via a program written on GAP~\cite{GAP}. 
Also it is easy to enlarge these results to the setting of some infinite-dimensional locally nilpotent Lie algebras, see Section~\ref{SLDynL}. 
We consider Proposition~\ref{P:crfd} from Section~\ref{SLDynL} as our third main result.

The structure of the paper is as follows. 
In Section 2 we recall the notation for finite-dimensional classical Lie algebras (roots, weights, etc.) used throughout the paper. 
In Section 3 we define the notion of extensive orbits. 
We also show how a classification of coadjoint orbits can be reduced to the classification of extensive orbits, see Theorem~\ref{T:extdec}.  
In Section 4 we discuss the case of elementary orbits as a good and very feasible example in this context (an orbit is elementary if its support consists of a single element).  
We start Section~5 from Theorem~\ref{Tnso} which connects the dimension of an orbit with its support.
Also Section 5 contains the tools from~\cite{GoodwinLeMagaardPaolini16} on which our main proofs and algorithms are based on: notions of C-patterns, C-quatterns and AL-moves. 
In Section 6, 7, 8 we classify the extensive coadjoint orbits in simple classical Lie algebras for dimensions 2, 4, 6 respectively. 
The answer for dimension 6 is rather long and it is printed out in 13 pages of Section 9. 
In Section 10 we discuss a connection of our results with the Isaac's conjecture for conjugacy classes of the respective $p$-groups. 
In Section 11 (Appendix) we generalize our results to the setting of Lie--Dynkin nil-algebras. 
In Section 12 (Appendix) we recall the enumeration of roots in CHEVIE package. This enumeration is a part of the description of orbits of dimensions 2, 4, 6 given in Section 6, 7, 8 respectively. 


\medskip\textsc{Acknowledgements}. We thank Aleksandr Panov, Aleksandr Shevchenko, Matvey Surkov and Mikhail Venchakov for useful discussions. We also thank Mikhail Panov for his useful suggestions about drawing pictures with \TeX.

The first author was supported by RSF (project No. 25--21--00219). 
The second author was supported by RSF (project No. 22--41--02028).

\sect{Finite-dimensional case: notation}\label{Sfdcn}\addcontentsline{toc}{subsection}{\ref{Sfdcn}. Finite-dimensional case: notation}
In this section, we first briefly recall definitions of classical finite-dimensional simple Lie algebras and fix notation for the nilradicals of their Borel subalgebras.
This will be used to state the main results in latter sections. 

Pick $n\in\Zp_{>0}$. Let $\gt$ denote one of the Lie algebras $\slt_n(\Fp)$, $\sot_{2n}(\Fp)$, $\sot_{2n+1}(\Fp)$ or $\spt_{2n}(\Fp)$. We will assume till the end of the paper that the characteristic of $\Fp$ is zero or greater than $\dim\gt$. The algebra $\sot_{2n}(\Fp)$ (respectively, $\sot_{2n+1}(\Fp)$ and $\spt_{2n}(\Fp)$) is realized as the sub\-al\-gebra of $\slt_{2n}(\Fp)$ (respectively, $\slt_{2n+1}(\Fp)$ and $\slt_{2n}(\Fp)$) consisting of all~$x$ such that $$\beta(u,xv)+\beta(xu,v)=0$$ for all $u,v$ in $\Fp^{2n}$ (respectively, in $\Fp^{2n+1}$ and $\Fp^{2n}$), where
\begin{equation*}
\beta(u,v)=\begin{cases}
\sum\nolimits_{i=1}^n(u_iv_{-i}+u_{-i}v_i)&\text{for }\sot_{2n}(\Fp),\\
u_0v_0+\sum\nolimits_{i=1}^n(u_iv_{-i}+u_{-i}v_i)&\text{for }\sot_{2n+1}(\Fp),\\
\sum\nolimits_{i=1}^n(u_iv_{-i}-u_{-i}v_i)&\text{for }\spt_{2n}(\Fp).
\end{cases}
\end{equation*}
Here for $\sot_{2n}(\Fp)$ (respectively, for $\sot_{2n+1}$ and $\spt_{2n}(\Fp))$ we denote by $e_1,\ldots,e_n,e_{-n},\ldots,e_{-1}$ (res\-pec\-tively, by $e_1,\ldots,e_n,e_0,e_{-n},\ldots,e_{-1}$ and $e_1,\ldots,e_n,e_{-n},\ldots,e_{-1}$) the standard basis of $\Fp^{2n}$ (res\-pec\-tively, of $\Fp^{2n+1}$ and $\Fp^{2n}$), and by $x_i$ the coordinate of a vector $x$ corresponding to $e_i$.

The set of all diagonal matrices from $\gt$ is a Cartan subalgebra of $\gt$; we denote it by $\htt$. Let $\Phi$ be the root system of $\gt$ with respect to $\htt$. Note that $\Phi$ is of type $A_{n-1}$ (respectively, $D_n$, $B_n$ and $C_n$) for $\slt_n(\Fp)$ (respectively, for $\sot_{2n}(\Fp)$, $\sot_{2n+1}(\Fp)$ and $\spt_{2n}(\Fp)$). The set of all upper-triangular matrices from $\gt$ is a Borel subalgebra of $\gt$ containing $\htt$; we denote it by $\bt$. Let $\Phi^+$ be the set of positive roots with respect to $\bt$. 
Denote by $\Pi$ the set of simple roots of $\Phi^+$ and by $Dyn$ the Dynkin diagram defined by $\Pi$. 
As usual, we identify $\Phi^+$ with the following subset of $\Rp^n$:
\begin{equation}\predisplaypenalty=0
\begin{split}
A_{n-1}^+&=\{\epsi_i-\epsi_j,~1\leq i<j\leq n\},\\
B_n^+&=\{\epsi_i-\epsi_j,~1\leq i<j\leq n\}\cup\{\epsi_i+\epsi_j,~1\leq i<j\leq n\}\cup\{\epsi_i,~1\leq i\leq n\},\\\label{formula:root_basis_nilradical_fin_dim}
C_n^+&=\{\epsi_i-\epsi_j,~1\leq i<j\leq n\}\cup\{\epsi_i+\epsi_j,~1\leq i<j\leq n\}\cup\{2\epsi_i,~1\leq i\leq n\},\\
D_n^+&=\{\epsi_i-\epsi_j,~1\leq i<j\leq n\}\cup\{\epsi_i+\epsi_j,~1\leq i<j\leq n\}.\\
\end{split}
\end{equation}
Here $\{\epsi_i\}_{i=1}^n$ is the standard basis of $\Rp^n$.

Denote by $\nt$ the algebra of all strictly upper-triangular matrices from $\gt$. Then $\nt$ has a basis consisting of root vectors $e_{\alpha}$, $\alpha\in\Phi^+$, where
\begin{equation*}\predisplaypenalty=0
\begin{split}
e_{\epsi_i}&=\sqrt{2}(e_{0,i}-e_{-i,0}),~e_{2\epsi_i}=e_{i,-i},\\
e_{\epsi_i-\epsi_j}&=\begin{cases}
e_{i,j}&\text{for }A_{n-1},\\
e_{i,j}-e_{-j,-i}&\text{for }B_n,~C_n\text{ and }D_n,
\end{cases}\\
e_{\epsi_i+\epsi_j}&=\begin{cases}
e_{i,-j}-e_{j,-i}&\text{for }B_n\text{ and }D_n,\\
e_{i,-j}+e_{j,-i}&\text{for }C_n,
\end{cases}
\end{split}
\end{equation*}
and $e_{i,j}$ are the usual elementary matrices. For $\sot_{2n}(\Fp)$ (respectively, for $\sot_{2n+1}(\Fp)$ and $\spt_{2n}(\Fp)$) we index the rows (from left to right) and the columns (from top to bottom) of matrices by the numbers $1,\ldots,n,-n,\ldots,-1$ (respectively, by the numbers $1,\ldots,n,0,-n,\ldots,-1$ and $1,\ldots,n,-n,\ldots,-1$). Note that $$\gt=\htt\oplus\nt\oplus\nt_-,$$ where $\nt_-=\langle e_{-\alpha},~\alpha\in\Phi^+\rangle_{\Fp}$, and, by definition, $e_{-\alpha}=e_{\alpha}^T$. (The superscript~$T$ always indicates matrix transposition.) The set $\{e_{\alpha},~\alpha\in\Phi\}$ can be uniquely extended to a Chevalley basis of $\gt$. Clearly, $\nt$ is the nilradical of the Borel subalgebra $\bt$ and $e_\alpha, \alpha\in\Phi^+,$ is a basis of $\nt^*$.

We denote $N=\exp(\nt)$, so that $\nt$ is the Lie algebra of the algebraic group $N$. The group $N$ acts on $\nt$ via the adjoint action, and the dual action of $N$ on $\nt^*$ is called \emph{coadjoint}. 
We denote by $e_\alpha^*$, $\alpha\in\Phi^+,$ the standard dual basis to $\nt^*$. 
Our first goal is to describe coadjoint orbits on $\nt^*$ of dimension less or equal to 6.
To do this, we need to introduce some additional notation.
\defi{Roots $\alpha, \gamma\in\Phi^+$ are called
$\beta$-\emph{singular} if $\alpha+\gamma=\beta$. The set of all $\beta$-singular roots is denoted by $\Sing(\beta)$.}
It is easy to see that
\begin{equation}\label{Esroot}
\begin{split}
\Sing(\epsi_i-\epsi_j)=&\bigcup_{l=i+1}^{j-1}\{\epsi_i-\epsi_l,\epsi_l-\epsi_j\},~1\leq i < j\leq n,\\
\Sing(\epsi_i)=&\bigcup_{l=i+1}^{n}\{\epsi_i-\epsi_l,\epsi_l\},~1\leq i\leq n,\\
\Sing(\epsi_i+\epsi_j)=&\bigcup_{l=i+1}^{j-1}\{\epsi_i-\epsi_l,\epsi_l+\epsi_j\}
\cup\bigcup_{l=j+1}^n\{\epsi_i-\epsi_l,\epsi_j+\epsi_l\}\cup\\
&\bigcup_{l=j+1}^n\{\epsi_i+\epsi_l,\epsi_j-\epsi_l\}\cup S_{ij}
,~1\leq i\le j\leq
n,\mbox{ where}\\
S_{ij}=&\begin{cases}\{\epsi_i, \epsi_j\},&\text{if }\Phi=B_n,\\
\{\epsi_i-\epsi_j, 2\epsi_j\},&\text{if }i\ne j,~\Phi=C_n,\\
\varnothing,&\text{if }\Phi=D_n
\end{cases}.
\end{split}
\end{equation}

Recall that there is a natural partial order on $\Phi^+$ defined as follows: 
for $\alpha, \beta\in\Phi$ we set $\alpha\le\beta$ (or $\beta\ge\alpha$) iff $\beta-\alpha$ is a sum of vectors from $\Phi^+$ or if $\alpha=\beta$. 
For example, $\epsi_i-\epsi_j\leq\epsi_a-\epsi_b$ if and only if $a\leq i$ and $j\leq b$; one can easily obtain similar conditions for other types of roots and root systems.

{\bf Remark.} It is convenient to draw schematically root systems under consideration as follows. 
We will draw $\Phi=A_{n-1}$ as the lower-triangular chessboard $n\times n$, where a root $\alpha=\epsi_i-\epsi_j$, $i<j$, corresponds to the box $(j,i)$. The root system $\Phi=C_n$ is drawn as a half of lower-triangular chessboard $2n\times2n$. 
The choice ``lower-triangular" for the upper-triangular matrices and the chessboard is used because we are mainly interested in the coadjoint representation (which is somehow symmetric to the adjoint representation). 

Precisely, we will index the rows and the columns of the chessboard by the indices 
$$1,\ldots,n,-n,\ldots,-1.$$
A root $\alpha=\epsi_i-\epsi_j$, $i<j$, will be drawn as the box $(j,i)$; a root $\alpha=\epsi_i+\epsi_j$ corresponds to the box $(j,-i)$. (Similar convention is applied for $\Phi=D_n$.) A root $\alpha=2\epsi_i$ corresponds to the box $(-i,i)$. Finally, $\Phi=B_n$ will be drawn as a half of the lower-triangular chessboard $(2n+1)\times(2n+1)$. We will enumerate rows and columns by the indices 
$$1,\ldots,n,0,-n,\ldots,-1.$$
Roots of the form $\alpha=\epsi_i\pm\epsi_j$ satisfy the same conditions as above, while a root $\alpha=\epsi_i$ corresponds to the box $(0,i)$.

The corresponding pictures for $A_4, B_4, C_4, D_4$ are given below.
$$\begin{array}{c}A_4:
~{\Autonumfalse\mymatrix{
\lNote{1}\Note{1}\Bot{2pt}& \Note{2}\pho&\Note{3}\pho&\Note{4}\pho&\\
\lNote{2}\gray\pho\Rt{2pt}&\Note{2}\Bot{2pt}\pho&\pho&&\\
\lNote{3}\gray\pho&\gray \Rt{2pt}&\Bot{2pt}&&\\
\lNote{4}\gray&\gray\star&\gray \Rt{2pt}&\Bot{2pt}&\\
\lNote{5}\gray&\gray&\gray &\gray\Rt{2pt}&\pho\\
}}\\
\star\to\epsi_2-\epsi_4
\end{array},
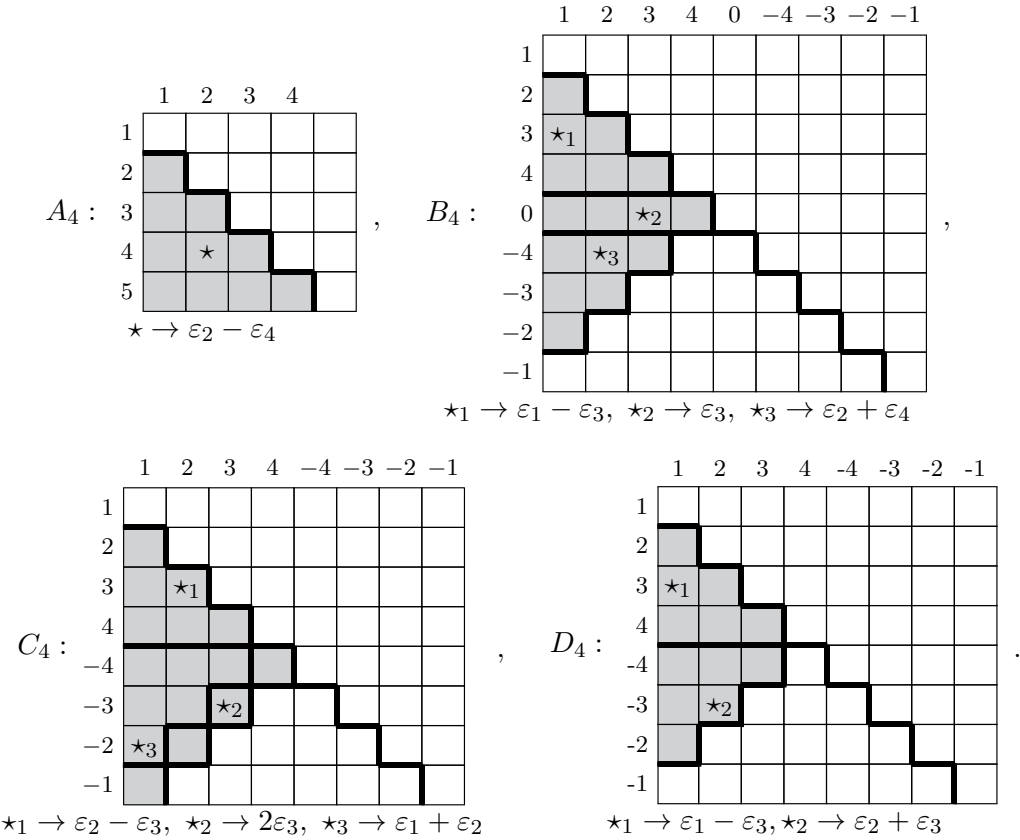
\hspace{10pt}
\begin{array}{c}B_4:
~{\Autonumfalse\mymatrix{
\lNote{1}\Note{1}\Bot{2pt}\pho& \Note{2}\pho& \Note{3}\pho& \Note{4}\pho& \Note{0}\pho& \Note{$-4$}\pho& \Note{$-3$}\pho& \Note{$-2$}\pho& \Note{$-1$}\pho\\
\lNote{2}\gray\Rt{2pt}\pho& \Bot{2pt}\pho& \pho& \pho& \pho& \pho& \pho& \pho& \pho\\
\lNote{3}\gray\star_1& \Rt{2pt}\gray\pho& \Bot{2pt}\pho& \pho& \pho& \pho& \pho& \pho& \pho\\
\lNote{4}\gray\Bot{2pt}\pho& \gray\Bot{2pt}\pho& \gray\Rt{2pt}\Bot{2pt}\pho& \Bot{2pt}\pho& \pho& \pho& \pho& \pho& \pho\\
\lNote{0}\gray\Bot{2pt}\pho& \gray\Bot{2pt}\pho& \gray\Bot{2pt}\star_2& \gray\Rt{2pt}\Bot{2pt}\pho& \Bot{2pt}\pho& \pho& \pho& \pho& \pho\\
\lNote{$-4$}\gray\pho& \gray\star_3& \Bot{2pt}\Rt{2pt}\gray\pho& \pho& \Rt{2pt}\pho& \Bot{2pt}\pho& \pho& \pho& \pho\\
\lNote{$-3$}\gray\pho& \gray\Bot{2pt}\Rt{2pt}\pho& \pho& \pho& \pho& \Rt{2pt}\pho& \Bot{2pt}\pho& \pho& \pho\\
\lNote{$-2$}\gray\Bot{2pt}\Rt{2pt}\pho& \pho& \pho& \pho& \pho& \pho& \Rt{2pt}\pho& \Bot{2pt}\pho& \pho\\
\lNote{$-1$}\pho& \pho& \pho& \pho& \pho& \pho& \pho& \Rt{2pt}\pho& \pho\\}}\\
\star_1\to\epsi_1-\epsi_3,~\star_2\to\epsi_3,~\star_3\to\epsi_2+\epsi_4
\end{array},\hspace{10pt}$$
$$\begin{array}{c}C_4:
{\Autonumfalse\mymatrix{
\lNote{1}\Note{1}\Bot{2pt}\pho& \Note{2}\pho& \Note{3}\pho& \Note{4}\pho& \Note{$-4$}\pho& \Note{$-3$}\pho& \Note{$-2$}\pho& \Note{$-1$}\pho\\
\lNote{2}\gray\Rt{2pt}\pho& \Bot{2pt}\pho& \pho& \pho& \pho& \pho& \pho& \pho\\
\lNote{3}\gray\pho& \gray\Rt{2pt}\star_1& \Bot{2pt}\pho& \pho& \pho& \pho& \pho& \pho\\
\lNote{4}\gray\Bot{2pt}\pho& \gray\Bot{2pt}& \gray\Bot{2pt}\Rt{2pt}\pho& \Bot{2pt}\pho& \pho& \pho& \pho& \pho\\
\lNote{$-4$}\gray\pho& \gray\pho& \gray\Bot{2pt}\Rt{2pt}& \gray\Rt{2pt}& \Bot{2pt}\pho& \pho& \pho& \pho\\
\lNote{$-3$}\gray\pho& \gray\Bot{2pt}\Rt{2pt}\pho& \gray\star_2& \Lft{2pt}\Top{2pt}\pho& \Rt{2pt}\pho& \Bot{2pt}\pho& \pho& \pho\\
\lNote{$-2$}\gray\Bot{2pt}\Rt{2pt}\star_3& \gray\pho& \Lft{2pt}\Top{2pt}\pho& \pho& \pho& \Rt{2pt}\pho& \Bot{2pt}\pho& \pho\\
\lNote{$-1$}\gray\pho& \Lft{2pt}\Top{2pt}\pho& \pho& \pho& \pho& \pho& \Rt{2pt}\pho& \pho\\}}\\
\star_1\to\epsi_2-\epsi_3,~\star_2\to2\epsi_3,~\star_3\to\epsi_1+\epsi_2
\end{array},\hspace{10pt}
\begin{array}{c}D_4:
~{\Autonumfalse
\mymatrix{
\lNote{1} \Note{1}\Bot{2pt}& \Note{2}\pho& \Note{3}\pho& \Note{4}\pho& \Note{-4}\pho& \Note{-3}\pho& \Note{-2}\pho& \Note{-1}\pho\\
\lNote{2} \gray\Note{1}\Rt{2pt}& \Note{2}\Bot{2pt}\pho& \Note{3}\pho&&&&&\\
\lNote{3}\gray\star_1&\gray \Rt{2pt}&\Bot{2pt}\pho&&&&&\\
\lNote{4}\gray\Bot{2pt}&\gray\Bot{2pt}&\gray\Bot{2pt}\Rt{2pt}&\Bot{2pt}&&&&\\
\lNote{-4}\gray&\gray&\gray\Rt{2pt}\Bot{2pt}&\Rt{2pt}&\Bot{2pt}&&&\\
\lNote{-3}\gray&\gray\star_2 \Rt{2pt}\Bot{2pt}&&&\Rt{2pt}&\Bot{2pt}&&\\
\lNote{-2}\gray\Bot{2pt}\Rt{2pt}& \pho& \pho&&&\Rt{2pt}&\Bot{2pt}&\\
\lNote{-1}& \pho& \pho&&&&\Rt{2pt}&\\
}}\\
\star_1\to\epsi_1-\epsi_3, \star_2\to\epsi_2+\epsi_3\end{array}.$$
In these pictures the cells corresponding to the roots of $\Phi^+$ are marked by gray color and for several such cells we provide the explicit presentation of the corresponding roots. 
It is easy to see that most of the cells are not gray and in some cases we do not draw the rows and columns which contain no gray cells, see example below. 
$$\begin{array}{c}A_4:
~{\Autonumfalse\mymatrix{
\lNote{2}\Note{1}\gray\pho\Top{2pt}\Rt{2pt}&\Note{2}\Note{2}\Bot{2pt}\pho&\Note{3}\pho&\Note{4}\pho\\
\lNote{3}\gray\pho&\gray \Rt{2pt}&\Bot{2pt}&\\
\lNote{4}\gray&\gray\star&\gray \Rt{2pt}\pho&\Bot{2pt}\pho\\
\lNote{5}\gray&\gray&\gray &\gray\Rt{2pt}\\
}}\\
\star\to\epsi_2-\epsi_4
\end{array},
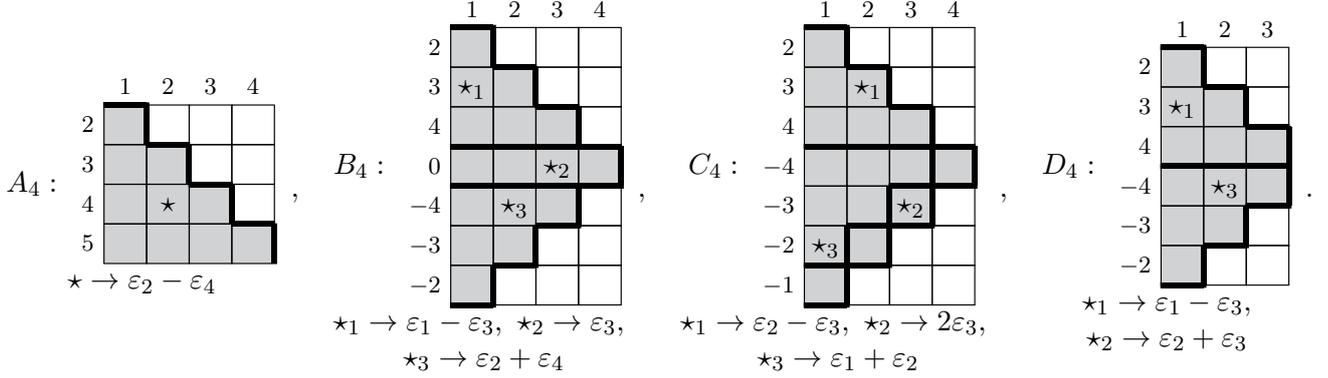
\hspace{6pt}
\begin{array}{c}B_4:
~{\Autonumfalse\mymatrix{
\lNote{2}\Note{1}\gray\Rt{2pt}\Top{2pt}\pho&\Note{2} \Bot{2pt}\pho&\Note{3} \pho&\Note{4} \pho\\
\lNote{3}\gray\star_1& \Rt{2pt}\gray\pho& \Bot{2pt}\pho& \pho\\
\lNote{4}\gray\Bot{2pt}\pho& \gray\Bot{2pt}\pho& \gray\Rt{2pt}\Bot{2pt}\pho& \Bot{2pt}\pho\\
\lNote{0}\gray\Bot{2pt}\pho& \gray\Bot{2pt}\pho& \gray\Bot{2pt}\star_2& \gray\Rt{2pt}\Bot{2pt}\pho\\
\lNote{$-4$}\gray\pho& \gray\star_3& \Bot{2pt}\Rt{2pt}\gray\pho& \pho\\
\lNote{$-3$}\gray\pho& \gray\Bot{2pt}\Rt{2pt}\pho& \pho&\\
\lNote{$-2$}\gray\Bot{2pt}\Rt{2pt}\pho& \pho& \pho&\\
}}\\
\star_1\to\epsi_1-\epsi_3,~\star_2\to\epsi_3,\\~\star_3\to\epsi_2+\epsi_4
\end{array},\hspace{6pt}
\begin{array}{c}C_4:
~{\Autonumfalse\mymatrix{
\lNote{2}\Note{1}\gray\Rt{2pt}\Top{2pt}\pho&\Note{2} \Bot{2pt}\pho&\Note{3} \pho&\Note{4} \pho\\
\lNote{3}\gray& \Rt{2pt}\gray\star_1& \Bot{2pt}\pho& \pho\\
\lNote{4}\gray\Bot{2pt}\pho& \gray\Bot{2pt}\pho& \gray\Rt{2pt}\Bot{2pt}\pho& \Bot{2pt}\pho\\
\lNote{$-4$}\gray\pho& \gray\pho& \Bot{2pt}\Rt{2pt}\gray\pho& \gray\Bot{2pt}\Rt{2pt}\pho\\
\lNote{$-3$}\gray\pho& \gray\Bot{2pt}\Rt{2pt}\pho& \gray\Bot{2pt}\Rt{2pt}\star_2&\\
\lNote{$-2$}\gray\Bot{2pt}\Rt{2pt}\star_3& \gray\Bot{2pt}\Rt{2pt}\pho& \pho&\\
\lNote{$-1$}\gray\Bot{2pt}\Rt{2pt}\pho& \pho& \pho&\\
}}\\
\star_1\to\epsi_2-\epsi_3,~\star_2\to2\epsi_3,\\~\star_3\to\epsi_1+\epsi_2
\end{array},\hspace{6pt}
\begin{array}{c}D_4:
~{\Autonumfalse\mymatrix{
\lNote{2}\Note{1}\Top{2pt}\gray\Rt{2pt}\pho&\Note{2} \Bot{2pt}\pho&\Note{3} \pho\\
\lNote{3}\gray\star_1& \Rt{2pt}\gray\pho& \Bot{2pt}\pho\\
\lNote{4}\gray\Bot{2pt}\pho& \gray\Bot{2pt}\pho& \gray\Rt{2pt}\Bot{2pt}\pho\\
\lNote{$-4$}\gray\pho& \gray\star_3& \Bot{2pt}\Rt{2pt}\gray\pho\\
\lNote{$-3$}\gray\pho& \gray\Bot{2pt}\Rt{2pt}\pho& \pho\\
\lNote{$-2$}\gray\Bot{2pt}\Rt{2pt}\pho& \pho& \pho\\
}}\\
\star_1\to\epsi_1-\epsi_3,\\\star_2\to\epsi_2+\epsi_3
\end{array}.$$

\sect{Extensive orbits}\label{sect:ext_orbits}\addcontentsline{toc}{subsection}{\ref{sect:ext_orbits}. Extensive orbits}
We use notation of Section~\ref{Sfdcn}. 
Recall that $\Supp{f}=\{\alpha\in\Phi^+\mid f(e_{\alpha})\neq0\}$. Set $$\NSupp{f}=\bigcup_{f'\in N.f}\Supp{f'}.$$

Consider $f\in \nt^*$. 
The following lemma provides an alternative description of $\NSupp{f}$.
\lemmp{\label{L:suppg}There is a non-empty open in Zariski topology subset $U\subset N.f$ such that $\Supp{f'}=\NSupp{f}$ for all $f'\in U$.}{If $f'(e_\alpha)\ne0$ for some $f'\in N.f$ and $\alpha\in\Phi^+$ then $e_\alpha\ne0$ on a certain open subset $U_\alpha$ of $N.f$. Thanks to the fact that $U$ is unipotent, $N.f$ is irreducible and hence $$U:=\bigcap_{\alpha\in\NSupp{f}}U_\alpha$$
is a non-empty open subset of $N.f$ which clearly has the desired properties.}

We need the following well-known observation.
\propp{\label{P:bou-p19} Let $\alpha, \beta\in\Phi^+$ be roots such that $\alpha\ge\beta$. Then there exists a sequence of roots $\gamma_1, \gamma_2,\ldots, \gamma_u\in\Phi^+$ such that $\beta=\alpha+\gamma_1+\ldots+\gamma_u$ and every partial sum $\alpha+\gamma_1+\ldots+\gamma_i$ is also a root.}
{It is implied by \cite[Chapter VI, Proposition 19]{Bou}. 
Indeed, by definition $\beta=\gamma_1+\gamma_2+\ldots+\gamma_u+\alpha$ for some $\gamma_1, \gamma_2, \ldots,\gamma_u\in\Phi^+$ and the above Proposition~19 implies that for some permutation of the multiset $\{\gamma_1, \gamma_2, \ldots, \gamma_u, \alpha\}$ every partial sum of this multiset is a positive root. 
Pick the last such sum not containing $\alpha$ and denote it by $\gamma_0$. 
Then replace the summands corresponding to $\gamma_0$ by $\gamma_0$ shortening the sequence. Then the new sequence will satisfy the needed property.}
\propp{\textup{a)} For $\alpha\le\beta$ the condition $\beta\in \NSupp{f}$~implies~$\alpha\in\NSupp{f}$.\\
\textup{b)} The Lie algebra $$\nt_{\Phi^+\setminus\NSupp{f}}:=\bigoplus_{\alpha\notin\NSupp{f}}\Fp e_\alpha$$ is an ideal of $\nt$.
\\
\textup{c)} One has $\NSupp{f}=\{\alpha\in\Phi^+\mid \exists \beta\in\Supp{f}\colon \alpha\le\beta\}$.
}{a) Assume to the contrary that there exist $\alpha, \beta\in\Phi^+$ such that $\beta\in\NSupp{f}$, $\beta\ge\alpha$, $\alpha\notin\NSupp{f}$. 
We may assume that $e_\beta(f)\ne0$, see Lemma~\ref{L:suppg}. 
Thanks to Proposition~\ref{P:bou-p19} we can assume that $\gamma=\beta-\alpha\in\Phi^+$. 
We have $(N.f)(e_\alpha)=0$. Thus, $f(N.e_\alpha)=0$. 
This implies that $f([\nt, e_\alpha])=0$. 
On the other hand, $[e_\alpha, e_\gamma]$ is a nonzero multiple of $e_\beta$, and thus $e_\beta(f)=0$. This is a contradiction. 

b) is implied by a). 

c) Thanks to a) we have $\NSupp{f}\supset\{\alpha\in\Phi^+\mid \exists \beta\in\Supp{f}: \alpha\le\beta\}$. 
We also have 
$$f(\nt_{\Phi^+\setminus\NSupp{f}})=0\implies N.f(\nt_{\Phi^+\setminus\NSupp{f}})=0\implies$$
$$\implies\NSupp{f}\subset\{\alpha\in\Phi^+\mid \exists \beta\in\Supp{f}\colon \alpha\le\beta\},$$ as needed.
}

\defi{\label{D:ext} We will say that two roots $\alpha, \beta\in\Pi$ are {\it adjacent} if $(\alpha, \beta)\ne0$. 
We will say that two roots $\alpha, \beta\in\Pi$ are {\it $f$-adjacent} if they are adjacent (in $\Pi$) and $\alpha+\beta\in\NSupp{f}$. We say that $f\in\nt^*$ is {\em extensive} if $\alpha, \beta\in\Pi$ are adjacent if and only if they are $f$-adjacent.}
The notion of $f$-adjacent roots defines the Dynkin subdiagram $\Dyn(f)$ of the Dynkin diagram of $\Pi$ (so, $\Dyn(f)$ is a subgraph of $\Dyn$ with the same vertices). 
Further, put $\Phi(f)$ to be the root system generated by $\Dyn(f)$; let $\nt(f)$ be the respective Lie algebra, and $N(f)$ be the corresponding unipotent group. 
\lemmp{We have $\NSupp{f}\subset\Phi(f)$.}{
Pick $\beta\in\NSupp{f}$. The case when $\beta$ is simple is clear, so let $\beta\notin\Pi$. It is a sum of several simple roots $\alpha_1, \alpha_2, \ldots$, where $\alpha_1\notin\alpha_2$. 
We claim that these roots are $f$-adjacent if and only if they are adjacent. This will imply the desired result. 

To prove the claim note that the fact that two roots $\alpha_1, \alpha_2$ are adjacent is equivalent to the condition $\alpha_1+\alpha_2\in\Phi^+$. By definition we have $\beta\ge \alpha_1+\alpha_2$. This and the definition of $\NSupp{f}$ prove the claim.  
}
Fix the following notation.

$\bullet$ $\Pi_1, \ldots, \Pi_l$ the irreducible components of $\Dyn(f)$.

$\bullet$ $\Phi_1,\ldots, \Phi_l$ the respective root subsystems of $\Phi$.

$\bullet$ $\nt_1,\ldots, \nt_l$ the respective subalgebras of $\nt$.

$\bullet$ $N_1,\ldots, N_l$ the respective subgroups of $N$.\\
Finally, we set $f_i:=f|_{\nt_i}$. Without loss of generality we can assume that $|\Pi_i|>1$ for the first $k$ entries and $|\Pi_i|=1$ elsewhere. 

For every $i$ we have a map $$\nt\to\nt_i\colon \begin{cases}e_\alpha\mapsto e_\alpha& \alpha\in\Phi_i^+,\\e_\alpha\mapsto 0&\alpha\notin\Phi_i^+.\end{cases}$$
It is easy to check that this map is a homomorphism of Lie algebras. This allows us to identify $\nt_i^*$ with an $\nt$-submodule of $\nt^*$, so one can consider $f_i$ as an element of $\nt^*$. This also shows that $$N.f\subset\bigoplus_{i=1}^l\nt_i^*\subset\nt^*.$$ 
\theop{\label{T:extdec}We have\\
\textup{a)} $f_i$ is a character of $\nt$ for $i>k$\textup, i.e., $N.f_{i}=\{f_{i}\}$\textup;\\
\textup{b)} $N.f_i=N_i.f_i$ for all $i$\textup;\\
\textup{c)} $f_i$ is extensive in $\nt_i^*$ for all $i$\textup;\\
\textup{d)} $N.f=N.f_1+ N.f_2+ \ldots+ N.f_k+(f_{k+1}+\ldots +f_{l})$ as $N$-varieties.}{
By definition we have $\Supp{f_i}\subset\Pi_i$ and $|\Pi_i|=1$ for $i>k$. 
Thus, $f_{i}$ is $N$-invariant for $i>k$. 
This implies a). 

Part b) is implied by the fact that the composition of maps $\nt_i\to\nt, \nt\to\nt_i$ is an isomorphism. 

Part c) follows from the definition of $f$-adjacent roots.

To prove part d) note that the direct sum $\nt\to\nt_1\oplus\nt_2\oplus\ldots\oplus\nt_l$ of maps $\nt\to\nt_i$ is surjective and moreover that the composition of maps $$\nt_1\oplus\nt_2\oplus\ldots\oplus\nt_l\to\nt, \nt\to \nt_1\oplus\nt_2\oplus\ldots\oplus\nt_l$$
is an isomorphism. This and a) implies d).
}
Theorem~\ref{T:extdec} shows that every coadjoint orbit of $N$ can be canonically decomposed as a direct sum of several (may be one or zero) extensive orbits for smaller subalgebras and a character. 
We would like to provide a rude estimate on the dimension of $N.f$ based on this decomposition. To do this we need the following definition. 
\defi{Set $\wdth(A_n):=2[\frac n2], \wdth(B_n):=2[\frac n2], \wdth(C_n):=2[\frac n2], \wdth(D_n):=2[\frac{n-1}2],\break \wdth(E_n):=2[\frac{n-1}2] (n=6, 7, 8), \wdth(F_4):=4, \wdth(G_2):=2$. For a nonconnected Dynkin diagram $Dyn$ we set $\wdth(Dyn)$ to be the sum of $\wdth(\cdot)$ for its parts.  
Finally we set $\wdth(f):=\wdth(\Dyn(f)).$}

\mprop{\label{PlWd}Let $f\in\nt^*$. Then $\dim N_i.f_i\ge\wdth(\Pi_i)$\textup, $1\le i\le l$\textup, hence $\dim N.f\ge\wdth(f)$.}
The proof of Proposition~\ref{PlWd} will be deduced from Lemma~\ref{Lbrank}.

Theorem~\ref{T:extdec} shows that it is enough to deal with extensive orbits $N_i.f_i$ for connected Dynkin diagrams. In particular, to describe the orbits of dimension 6 or less it is enough to describe the extensive orbits of dimensions 6 or less for such diagrams. And thanks to Proposition~\ref{PlWd} such orbits exist in classical cases only for Lie algebras of rank 8 or less. 

We will use an approach of~\cite{GoodwinLeMagaardPaolini16} to describe the orbits of dimension 6 or less for the classical cases. 
We will start from extensive orbits of dimension 2 in Section~\ref{S:dim2} then proceed to extensive orbits of dimension 4 in Section~\ref{S:dim4} and finally we will deal with extensive orbits of dimension 6 in Section~\ref{S:dim6}.

\sect{Elementary orbits}\label{sect:elementary}\addcontentsline{toc}{subsection}{\ref{sect:elementary}. Elementary orbits}
It is convenient to start the description of coadjoint orbits from the orbits whose support consist of a single root $\alpha\in\Phi^+$ (we call such orbits {\it elementary}). This case is very well-known, see, e.g.,~\cite{Mukherjee05}, but we wish to discuss it in detail to illustrate how our methods and ideas work in this (relatively simple) case. Recall that, given a linear form $f\in\nt^*$, we denote by $N.f$ its coadjoint orbit. Being an orbit of a connected unipotent group on an affine space, $N.f$ is an affine variety (in fact, it is isomorphic to an affine space); we denote by $\dim N.f$ its dimension. 
\mlemm{\label{Lrkalpha} Let $\alpha\in\Phi^+$ be a root. Then $\dim N.e_\alpha^*=|\Sing(\alpha)|$.}
To prove this lemma, we will reduce the calculation of $\dim N.e_\alpha^*$ to the computation of the rank of a certain skew-symmetric bilinear form.
\defi{Let $f\in\nt^*$ be a linear form. We set $B_f$ to be the following bilinear skew-symmetric form: $\nt\times \nt \to \mathbb C\colon x, y\mapsto f([x, y])$. 
Denote by $\ker B_f$ the kernel of $B_f$. 
Denote by $\rk B_f$ the rank of $B_f$. 
}

The following lemma is well known.
\lemm{\label{Lbrank}We have $\dim N.f=\rk B_f$.}
\proof{~of~Lemma~\ref{Lrkalpha}}{ It is enough to check that $\rk B_{e_\alpha^*}=|\Sing(\alpha)|$. 
It is easy to verify that $e_{\beta}\in\ker B_f$ if and only if $\beta\notin \Sing(\alpha)$. Further we have that $B_{e_\alpha^*}(e_{\alpha_1}, e_{\alpha_2})\ne0$ if and only if $\alpha_1+\alpha_2=\alpha$. 
Therefore, the matrix of $B_{e_\alpha^*}$ in the basis $\{e_\alpha,~\alpha\in\Phi^+\}$, splits into $\dfrac{|\Sing(\alpha)|}2$ two-by-two blocks $\begin{pmatrix}0&c\\-c&0\end{pmatrix}$ with $c\ne0$. 
This implies the desired result.}
\begin{proof}{~of Proposition~\ref{PlWd}}{
Without loss of generality we replace $f$ by a generic element of $N.f$. 
Also without loss of generality we assume that the underlying Dynkin diagram $Dyn$ is connected and $Dyn=\Dyn(f)$. 
Then thanks to Lemma~\ref{L:suppg} we have $B_f(e_\alpha, e_\beta)\ne0$ for $\alpha, \beta\in\Pi$ whenever $\alpha, \beta$ are adjacent. 
We consider the minor $M$ of $B_f$ defined by the set of vectors $e_\alpha, \alpha\in\Pi$. 
We claim that the rank of this minor is exactly $\wdth(Dyn)$. 
The required result can be easily deduced from this statement. 
Without loss of generality we assume that $Dyn$ is connected. To prove this claim we assume first that $Dyn$ consists of one point (i.e., it is $A_1$). Then the claim is trivial. 

In the other case let $\alpha$ be the corner edge of the respective Dynkin diagram $Dyn$, i.e., such an edge that there exists unique $\beta\in\Pi\setminus\{\alpha\}$ such that the scalar product of $\beta$ and $\alpha$ is nonzero. 
Then the $e_\alpha$-row of $M$ has the only non-zero element at the intersection with the $e_\beta$-column. 
This implies that if we drop $e_\alpha$-row and $e_\beta$-column then the rank will decrease by 1. The same is true for $e_\alpha$-column and $e_\beta$-row. 
This procedure removes two adjacent vertices $\alpha, \beta$ (one of which has to be the corner edge) from our Dynkin diagram and reduces the respective rank by 2. Repeating this procedure up to the moment when we get zero-matrix (which has rank 0), we prove the claim.
}\end{proof}

It is quite straightforward to deduce the list of orbits of dimensions $2, 4, 6$ of the form $N.e_{\alpha}^*$ from Lemma~\ref{Lrkalpha}. 
We provide the lists of such $\alpha$'s in Tables~\ref{T:dim2elm},~\ref{T:dim4elm},~\ref{T:dim6elm} below, c.f. tables of Theorems~\ref{Tdim2},~\ref{Tdim4},~\ref{Tdim6}. 
\begin{longtable}{|c|c|c|}
\hline \multicolumn{3}{|c|}{Table~\label{T:dim2elm}\ref{T:dim2elm}: $\alpha$ with $\dim N.e_{\alpha}^*=2$}\\
\hline Root system&Root~$\alpha$& $\Sing(\alpha)$\\
\hline $A_n$& $\varepsilon_i-\varepsilon_{i+2}$& $\varepsilon_i-\varepsilon_{i+1}, \varepsilon_{i+1}-\varepsilon_{i+2}$\\
\hline $B_n$& $\varepsilon_i-\varepsilon_{i+2}$& $\varepsilon_i-\varepsilon_{i+1}, \varepsilon_{i+1}-\varepsilon_{i+2}$\\
\hline $B_n$& $\varepsilon_{n-1}+\varepsilon_{n}$& $\varepsilon_{n-1}, \varepsilon_{n}$\\

\hline $B_n$& $\varepsilon_{n-1}$& $\varepsilon_{n-1}-\varepsilon_{n}, \varepsilon_{n}$\\
\hline $C_n$& $\varepsilon_i-\varepsilon_{i+2}$& $\varepsilon_i-\varepsilon_{i+1}, \varepsilon_{i+1}-\varepsilon_{i+2}$\\
\hline $C_n$& $\varepsilon_{n-1}+\varepsilon_{n}$& $\varepsilon_{n-1}-\varepsilon_n, 2\varepsilon_{n}$\\
\hline $C_n$& $2\varepsilon_{n-1}$& $\varepsilon_{n-1}-\varepsilon_{n}, \varepsilon_{n-1}+\varepsilon_{n}$\\
\hline $D_n$& $\varepsilon_i-\varepsilon_{i+2}$& $\varepsilon_i-\varepsilon_{i+1}, \varepsilon_{i+1}-\varepsilon_{i+2}$\\
\hline $D_n$& $\varepsilon_{n-2}+\varepsilon_{n}$& $\varepsilon_{n-2}-\varepsilon_{n-1}, \varepsilon_{n-1}+\varepsilon_{n}$\\
\hline \end{longtable}
\begin{longtable}{|c|c|c|}
\hline \multicolumn{3}{|c|}{Table~\label{T:dim4elm}\ref{T:dim4elm}: $\alpha$ with $\dim N.e_{\alpha}^*=4$}\\
\hline Root system&Root~$\alpha$& $\Sing(\alpha)$\\
\hline $A_n$& $\varepsilon_i-\varepsilon_{i+3}$& $\varepsilon_i-\varepsilon_{i+1}, \varepsilon_{i}-\varepsilon_{i+2}, \varepsilon_{i+1}-\varepsilon_{i+3}, \varepsilon_{i+2}-\varepsilon_{i+3}$\\

\hline $B_n$& $\varepsilon_i-\varepsilon_{i+3}$& $\varepsilon_i-\varepsilon_{i+1}, \varepsilon_{i}-\varepsilon_{i+2}, \varepsilon_{i+1}-\varepsilon_{i+3}, \varepsilon_{i+2}-\varepsilon_{i+3}$\\
\hline $B_n$& $\varepsilon_{n-2}+\varepsilon_{n}$& $\varepsilon_{n-2}, \varepsilon_{n}, \varepsilon_{n-2}-\varepsilon_{n-1}, \varepsilon_{n-1}+\varepsilon_{n}$\\
\hline $B_n$& $\varepsilon_{n-2}$& $\varepsilon_{n-2}-\varepsilon_{n}, \varepsilon_{n}, \varepsilon_{n-2}-\varepsilon_{n-1}, \varepsilon_{n-1}$\\

\hline $C_n$& $\varepsilon_i-\varepsilon_{i+3}$& $\varepsilon_i-\varepsilon_{i+1}, \varepsilon_{i}-\varepsilon_{i+2}, \varepsilon_{i+1}-\varepsilon_{i+3}, \varepsilon_{i+2}-\varepsilon_{i+3}$\\
\hline $C_n$& $\varepsilon_{n-2}+\varepsilon_{n}$& $\varepsilon_{n-2}-\varepsilon_n, 2\varepsilon_{n}, \varepsilon_{n-2}-\varepsilon_{n-1}, \varepsilon_{n-1}+\varepsilon_{n-2}$\\
\hline $C_n$& $2\varepsilon_{n-2}$& $\varepsilon_{n-2}\pm\varepsilon_{n}, \varepsilon_{n-2}\pm\varepsilon_{n-1}$\\

\hline $D_n$& $\varepsilon_i-\varepsilon_{i+3}$& $\varepsilon_i-\varepsilon_{i+1}, \varepsilon_{i}-\varepsilon_{i+2}, \varepsilon_{i+1}-\varepsilon_{i+3}, \varepsilon_{i+2}-\varepsilon_{i+3}$\\
\hline $D_n$& $\varepsilon_{n-3}+\varepsilon_{n}$& $\varepsilon_{n-3}-\varepsilon_{n-1}, \varepsilon_{n-1}+\varepsilon_{n}, \varepsilon_{n-3}-\varepsilon_{n-2}, \varepsilon_{n-2}+\varepsilon_{n}$\\
\hline $D_n$& $\varepsilon_{n-2}+\varepsilon_{n-1}$& $\varepsilon_{n-2}+\varepsilon_{n}, \varepsilon_{n-1}-\varepsilon_{n}, \varepsilon_{n-2}-\varepsilon_{n}, \varepsilon_{n-2}+\varepsilon_{n}$\\

\hline \end{longtable}

\begin{longtable}{|c|c|c|}
\hline \multicolumn{3}{|c|}{Table~\label{T:dim6elm}\ref{T:dim6elm}: $\alpha$ with $\dim N.e_{\alpha}^*=6$}\\
\hline Root system&Root~$\alpha$& $\Sing(\alpha)$\\
\hline $A_n$& $\varepsilon_i-\varepsilon_{i+4}$& $\varepsilon_i-\varepsilon_{i+1}, \varepsilon_{i}-\varepsilon_{i+2}, \varepsilon_{i}-\varepsilon_{i+3}, \varepsilon_{i}-\varepsilon_{i+4}, \varepsilon_{i+1}-\varepsilon_{i+4}, \varepsilon_{i+2}-\varepsilon_{i+4}, \varepsilon_{i+3}-\varepsilon_{i+4}$\\

\hline $B_n$& $\varepsilon_i-\varepsilon_{i+4}$& $\varepsilon_i-\varepsilon_{i+1}, \varepsilon_{i}-\varepsilon_{i+2}, \varepsilon_{i}-\varepsilon_{i+3}, \varepsilon_{i}-\varepsilon_{i+4}, \varepsilon_{i+1}-\varepsilon_{i+4}, \varepsilon_{i+2}-\varepsilon_{i+4}, \varepsilon_{i+3}-\varepsilon_{i+4}$\\
\hline $B_n$& $\varepsilon_{n-3}+\varepsilon_{n}$& $\varepsilon_{n-3}, \varepsilon_{n}, \varepsilon_{n-3}-\varepsilon_{n-2}, \varepsilon_{n-3}-\varepsilon_{n-1}, \varepsilon_{n-1}+\varepsilon_{n}, \varepsilon_{n-2}+\varepsilon_{n}$\\
\hline $B_n$& $\varepsilon_{n-2}+\varepsilon_{n-1}$& $\varepsilon_{n-2}, \varepsilon_{n-1}, \varepsilon_{n-2}\pm\varepsilon_{n}, \varepsilon_{n-1}\pm\varepsilon_n$\\
\hline $B_n$& $\varepsilon_{n-3}$& $\varepsilon_{n-3}-\varepsilon_{n-2}, \varepsilon_{n-3}-\varepsilon_{n-1}, \varepsilon_{n-3}-\varepsilon_{n}, \varepsilon_{n-2}, \varepsilon_{n-1}, \varepsilon_{n-2}$\\

\hline $C_n$& $\varepsilon_i-\varepsilon_{i+4}$& $\varepsilon_i-\varepsilon_{i+1}, \varepsilon_{i}-\varepsilon_{i+2}, \varepsilon_{i}-\varepsilon_{i+3}, \varepsilon_{i}-\varepsilon_{i+4}, \varepsilon_{i+1}-\varepsilon_{i+4}, \varepsilon_{i+2}-\varepsilon_{i+4}, \varepsilon_{i+3}-\varepsilon_{i+4}$\\
\hline $C_n$& $\varepsilon_{n-3}+\varepsilon_{n}$& $\varepsilon_{n-3}-\varepsilon_n, 2\varepsilon_{n}, \varepsilon_{n-3}-\varepsilon_{n-2}, \varepsilon_{n-3}-\varepsilon_{n-1}, \varepsilon_{n-2}+\varepsilon_n, \varepsilon_{n-1}+\varepsilon_n$\\
\hline $C_n$& $\varepsilon_{n-2}+\varepsilon_{n-1}$& $\varepsilon_{n-2}-\varepsilon_{n-1}, 2\varepsilon_{n-1}, \varepsilon_{n-3}\pm\varepsilon_n, \varepsilon_{n-1}\pm\varepsilon_n$\\
\hline $C_n$& $2\varepsilon_{n-3}$& $\varepsilon_{n-3}\pm\varepsilon_{n-2}, \varepsilon_{n-3}\pm\varepsilon_{n-1}, \varepsilon_{n-3}\pm\varepsilon_{n}$\\

\hline $D_n$& $\varepsilon_i-\varepsilon_{i+4}$& $\varepsilon_i-\varepsilon_{i+1}, \varepsilon_{i}-\varepsilon_{i+2}, \varepsilon_{i}-\varepsilon_{i+3}, \varepsilon_{i}-\varepsilon_{i+4}, \varepsilon_{i+1}-\varepsilon_{i+4}, \varepsilon_{i+2}-\varepsilon_{i+4}, \varepsilon_{i+3}-\varepsilon_{i+4}$\\
\hline $D_n$& $\varepsilon_{n-4}+\varepsilon_{n}$& $\varepsilon_{n-4}-\varepsilon_{n-3}, \varepsilon_{n-4}-\varepsilon_{n-2}, \varepsilon_{n-4}-\varepsilon_{n-1},  \varepsilon_{n-3}+\varepsilon_{n}, \varepsilon_{n-2}+\varepsilon_{n}, \varepsilon_{n-1}+\varepsilon_{n}$\\
\hline $D_n$& $\varepsilon_{n-3}+\varepsilon_{n-1}$& $\varepsilon_{n-3}-\varepsilon_{n-2}, \varepsilon_{n-3}\pm\varepsilon_{n}, \varepsilon_{n-2}+\varepsilon_{n-1}, \varepsilon_{n-1}\pm\varepsilon_{n}$\\
\hline \end{longtable}
And finally we provide a list of all extensive elementary orbits of dimensions 2, 4, 6.
\begin{longtable}{|c|c|c|c|c|c|}
\hline \multicolumn{6}{|c|}{Table~\label{Tbl:dim246ExtElm}\ref{Tbl:dim246ExtElm}: $\alpha$ with extensive $N.e_{\alpha}^*$ of dimensions $2, 4, 6$}\\
\hline \multicolumn{2}{|c|}{$\dim N.e_{\alpha}^*=2$}&\multicolumn{2}{|c|}{$\dim N.e_{\alpha}^*=4$}&\multicolumn{2}{|c|}{$\dim N.e_{\alpha}^*=6$}\\
\hline Root system&Root~$\alpha$&Root system&Root~$\alpha$&Root system&Root~$\alpha$\\
\hline $A_2$&$\varepsilon_1-\varepsilon_3$&$A_3$&$\varepsilon_1-\varepsilon_4$&$A_4$&$\varepsilon_1-\varepsilon_5$\\
\hline $B_2$&$\varepsilon_1$&$B_3$&$\varepsilon_1$&$B_3$&$\varepsilon_1+\varepsilon_2$\\
\hline $B_2$&$\varepsilon_1+\varepsilon_2$&$B_3$&$\varepsilon_1+\varepsilon_3$&$B_4$&$\varepsilon_1$\\
\hline $C_2$&$\varepsilon_1+\varepsilon_2$&$C_3$&$\varepsilon_1+\varepsilon_3$&$B_4$&$\varepsilon_1+\varepsilon_4$\\
\hline $C_2$&$2\varepsilon_1$&$C_3$&$2\varepsilon_1$&$C_3$&$\varepsilon_1+\varepsilon_2$\\
\hline &&$D_4$&$\varepsilon_1\pm\varepsilon_4$&$C_4$&$2\varepsilon_1$\\
\hline &&&&$D_4$&$\varepsilon_1\pm\varepsilon_4$\\
\hline &&&&$D_4$&$\varepsilon_1+\varepsilon_3$\\
\hline &&&&$D_5$&$\varepsilon_1\pm\varepsilon_5$\\
\hline \end{longtable}

\sect{A restriction on the supports in the classical cases}\label{sect:restrict}\addcontentsline{toc}{subsection}{\ref{sect:restrict}. A restriction on the supports in the classical cases}
The following result provides a uniform restriction on the supports of orbits of dimension $2d$ for the classical Lie algebras.
\mtheo{\label{Tnso} Assume that the defining root system $\Phi$ for $\nt$ is classical. Set $\delta_0$ to be the highest root of $\Phi^+$ and pick $d\in\mathbb Z_{\ge 0}$ such that $2d\le |\Sing(\delta_0)|$. Then the union of the supports of all orbits of dimension $2d$ is equal to the set 
$$\SuppD_d:=
\{\alpha\in\Phi^+\mid \exists \alpha'\in\Phi^+\colon \alpha'\ge\alpha\text{ and }|\Sing(\alpha')|=2d\}.$$
}
Note that for root systems of types $A_n, B_n, D_n$ the condition of Theorem~\ref{Tnso} is equivalent to $|\Sing(\alpha)|\le 2d$.
Theorem~\ref{Tnso} follows from the next proposition. 
\mprop{\label{Psupp} \textup{a)} Let $f\in\nt^*$ and $\alpha\in \NSupp{f}$ be a root. 
If $\Phi$ is of type $A_n, B_n, C_n$ or $D_n$ and $\alpha$ is maximal in $\NSupp{f}$ with respect to $\le$ then $|\Sing(\alpha)|\le \dim N.f$.\\
\textup{b)} Let $f\in\nt^*$ and $\alpha\in \NSupp{f}$ be a root. 
If $\Phi$ is of type $A_n, B_n$ or $D_n$ then $|\Sing(\alpha)|\le \dim N.f$.
}
The statement of Proposition~\ref{Psupp}b) is a bit more elegant than the statement of Proposition~\ref{Psupp}a) and Proposition~\ref{Psupp}b) avoids $C_n$ case. 
The issue here is that the statement of Proposition~\ref{Psupp}b) is wrong for $C_n, n\ge 3$.


{
\begin{proof}{~of Proposition~\ref{Psupp} for types $A$ and $D$}{~
We will show that there are no orbits $N.f$ such that $\dim(N.f)=2d$ and $\NSupp{f}\not\subset\SuppD_d$. 

Assume to the contrary that such an orbit $N.f$ exists and pick a maximal with respect to $\le$ root $\gamma\in\NSupp{f}\setminus\SuppD_d$. 
Without loss of generality we can assume that $\gamma\in\Supp{f}$. 
Recall that $\dim(N.f)=\rk B_f$ and consider the minor of $B_f$ standing for basis vectors $e_\alpha, \alpha\in \Sing(\gamma)$. 

We claim that in $A_n$ case and $D_n$ case this minor splits into $|\Sing(\gamma)|/2$ two-by-two non-zero skew-symmetric blocks defined by pairs $\alpha, \beta\in\Phi^+$ with $\alpha+\beta=\gamma$. The rank of each such a block is equal to 2 and thus the rank of the whole minor is $|\Sing(\gamma)|$. This implies $\dim N.f\ge|\Sing(\gamma)|$. 

To prove the claim we need to show that all other entries of the above minor are 0. To do this we prove the following lemma.
\lemmp{\label{L:sumr3}Assume that for some $\alpha_1, \alpha_2, \beta_1, \beta_2,\gamma\in\Phi^+$ we have $$\alpha_1+\beta_1=\alpha_2+\beta_2=\gamma,~\{\alpha_1, \beta_1\}\ne\{\alpha_2, \beta_2\}.$$
Then $\alpha_1+\alpha_2\notin\Phi^+$.}
{Assume to the contrary that $\alpha_1+\alpha_2=\tau\in\Phi^+$. 
Consider the space $V$ generated by $\gamma, \alpha_1, \alpha_2$ and set $\underline{\Phi}:=V\cap\Phi$. It is clear that $\underline{\Phi}$ is a root system and that $\underline{\Phi}$ contains all roots under consideration: $\alpha_1, \beta_1, \alpha_2, \beta_2, \gamma, \tau$. 

By definition $\underline{\Phi}$ is of rank 3 or less.
Note also that it is simply-laced because $\underline{\Phi}$ is a subset of a root system of type $A_n$ or $D_n$. A straightforward check shows that there are no such roots $\alpha_1, \alpha_2, \beta_1, \beta_2$ in the simply-laced root systems of rank 3 or less.
}
Thus, the required result is proved.
}
\end{proof}
We require a different proof for $B/C$ cases of Proposition~\ref{Psupp} because the analogue of Lemma~\ref{L:sumr3} is strictly wrong for $B_3$ and $C_3$ and thus for $B_n, C_n, n\ge 3$. 
To deal with this issue we need more notation and concepts. 
\sst{C-patterns and C-quatterns}
Let $X$ be a subset of $\Phi^+$. Set $\nt_X:=\bigoplus_{\alpha\in X}\Fp e_\alpha$ with the (Lie) bracket given by the formula $$[e_\alpha, e_\beta]=\begin{cases}[e_\alpha, e_\beta],& {\rm~if~}\alpha+\beta\in X,\\0,&\text{if }\alpha+\beta\notin X\end{cases}$$ for all $\alpha, \beta\in X$
(note that if $\alpha, \beta, \alpha+\beta\in \Phi$ then $[e_\alpha, e_\beta]$ is proportional to $e_{\alpha+\beta})$. Such a bilinear bracket $[\cdot, \cdot ]$ do not always define a Lie algebra but it does define a Lie algebra under the assumption that $X$ is a $C$-pattern (a shorthand for ``combinatorial counterpart of a pattern subgroup'') or $C$-quattern (a shorthand for ``combinatorial counterpart of a quattern subgroup''), see the next definitions. 

\defi{We will say that $X$ is a {\it C-pattern}  if $X$ satisfies the following condition: $$\text{if }\alpha, \beta\in X\text{ and }\alpha+\beta\in\Phi^+\text{ then }\alpha+\beta\in X.$$
This condition is clearly equivalent to the statement that $\nt_X$ is a sub-Lie algebra of $\nt$.}
\defi{ We will say that $X$ is a {\it C-quattern} if $X=X_+\setminus X_-$ for two C-patterns $X_+, X_-$ such that $X_-\subset X_+$ and
$$\text{if }\alpha_-\in X_-,\alpha_+\in X_+~\text{ and }\alpha_-+\alpha_+\in X_+ \text{ then }\alpha_-+\alpha_+\in X_-.$$
This is clearly equivalent to the statement that $\nt_{X_-}$ is an ideal of the Lie algebra $\nt_{X_+}$ and in this case we have $\nt_X\cong \nt_{X_+}/\nt_{X_-}$ where the right hand side is a quotient of Lie algebra by its ideal and thus the left hand side is a Lie algebra.}
From now on we assume that $X$ is always a C-quattern. We denote by $N_X$ the unipotent group with Lie algebra $\nt_X$.  
\defi{We will say that a quattern $X$ is {\it large} if there exists the above presentation $X=X_+\setminus X_-$ with $X_+=\Phi^+$. This means that $\nt_X$ is a quotient of $\nt$ and $N_X$ is a quotient of $N$.}
Note that by Proposition~\ref{P:bou-p19} the set $\NSupp{f}$ is a large C-quattern for all $f\in\nt^*$ .

Every C-pattern $X$ is a C-quattern for $X_+=X, X_-={\varnothing}$ and for every C-quattern $X$ the data $(\nt_X, [\cdot, \cdot])$ provides a structure of Lie algebra on $\nt_X$. 
Quattern and pattern subgroups were considered in~\cite{GoodwinLeMagaardPaolini16} and we just rephrase the respective definitions to make them a bit more combinatorially focused.








\sst{Saturated linear forms and AL-moves} We use the notation and conventions of the previous subsection. 
Set 
$$Z(X):=\{\alpha\in X\mid (\alpha+X)\cap X={\varnothing}\}.$$ 
Then it is easy to verify that $\nt_{Z(X)}$ coincides with the center of $\nt_X$. In particular, this implies that $Z(X)\ne0$ for any C-quattern $X$.

\defi{Pick $f\in \nt_X^*$ and consider $Z\subset Z(X)$. We will say that $f$ is {\it $Z$-saturated} if $f(e_\alpha)\ne0$ for all $\alpha\in Z$. We will say that $f$ is {\it saturated} if $f$ is $Z(X)$-saturated. 
We denote the variety of $Z$-saturated linear forms by $\nt^*_{X; Z}$.}
It is clear that each $\nt_{X; Z}^*$ is $N$-stable and thus is a collection of coadjoint orbits. 

We map $\nt^*_{X; Z}$ to $\nt^*$ by the formula $e_\alpha^*\mapsto e_\alpha^*$ and denote by $\underline{\nt}^*_{X; Z}$ the image of this map. 
We will frequently identify $\nt_{X; Z}^*$ with $\underline{\nt}_{X; Z}^*$. 
\lemmp{Let $X$ be a large C-quattern. 
Pick $f\in \nt_X^*$ and denote by $\hat f$ the preimage of $f$ in $\nt^*$. Then $f$ is saturated if and only if $\NSupp{\hat f}=X$.}{
Follows from Proposition~\ref{P:bou-p19}.
}
\defi{Let $Z\subset  Z(X)$ be a subset and let $Y\subset \nt_{X; Z}^*$ be a subvariety. 
We say that $Y$ is a {\em set-section} for $\nt_{X; Z}^*$ if $Y$ intersects each $N_X$-orbit of $\nt_{X; Z}^*$ in exactly one point (note that $e_\alpha, \alpha\in Z,$ is constant on $N_X$-orbits if we consider $e_\alpha$ as a function defined on $\nt_X^*$).}
It is clear that a set-section is a section of the native quotient map $\nt_{X; Z}^*\to \nt_{X; Z}^*/N_X$  but only in set-theoretic terms. 
Also note that all the set-sections which we will consider explicitly will be unions of the form $V(S_1)\sqcup V(S_2)\sqcup\ldots$ for some subsets $S_1, S_2, \ldots\subset\Phi^+$. 
\lemmp{\label{L:tpS}Let $X$ be a quattern. Consider $\gamma\in Z(X)$, $Z\subset Z(X)\setminus\{\gamma\}$. Let $Y_0$ be a set-section of $\nt_{X\setminus\{\gamma\}; Z}^*$ and $Y_+$ be a set-section for $\nt_{X; Z\cup\gamma}^*$
Then $Y_0\cap Y_{+}=\varnothing$ and $\phi(Y_0)\cup Y_{+}$ is a set-section for $\nt_{X; Z}^*$ where $\phi(Y_0)$ is the image of $Y_0$ under the natural inclusion $\nt_{X\setminus\{\gamma\}}^*\to\nt_X^*$.}{This is implied by the following facts: all the coadjoint orbits of $\nt_X$ split into two classes --- those for which $f(e_\gamma)\ne0$ and those for which $f(e_\gamma)=0$. By definition $Y_+$ is a set-section for the orbits from the first class, while $Y_0$ is a set-section for the second one. 
Clearly, these classes do not intersect and $\phi(Y_0)\cap Y_+={\varnothing}$. 
Thus $\phi(Y_0)\sqcup Y_+$ is a set section for $\nt_{X;Z\cup\{\gamma\}}$.
}
In the next proposition we refer to the notion of polarization, see~\cite[\S 1.12]{Dixmier1} for its definition and properties.
\mprop{\label{P:gdid31}Let $X$ be a quattern and let $\gamma \in Z(X)$. Suppose that there exist $\delta, \beta \in X \setminus \{\gamma\})$ with $\delta+\beta=\gamma$ such that 
\begin{center}$X\cap (\beta-X)={\varnothing}$ and $X\cap (\delta+X)=\{\gamma\}$.\end{center}
\textup{a)} Then $X\setminus\{\alpha, \beta\}$ is a C-quattern and $Z(X)\subset Z(X\setminus\{\alpha, \beta\})$.\\
\textup{b)} Assume that $Y$ is a set-section for $\nt_{X\setminus\{\alpha, \beta\};Z(X)}^*$. Then $\phi(Y)$ is a set-section for $\nt_{X; Z}^*$\textup, where $\phi$ is the linear map $\nt_{X\setminus\{\alpha, \beta\}}^*\to\nt_X^*$ defined as follows\textup: $$\phi(e_\tau^*)=e_\tau^*\in\nt_X^*.$$\\
\textup{c)} Pick $\lambda\in Y$ and assume that $\mathfrak p\subset \nt_{X\setminus\{\alpha, \beta\}}$ is a polarization of $\lambda$. Then $\mathfrak p\oplus\Fp e_{\delta}$ is a polarization of~$\phi(\lambda)$. This also implies that $\dim (N_X.\phi(\lambda))=2+\dim (N_{X\setminus\{\alpha, \beta\}}.\lambda)$.
}
We will say that $\beta$ from the above lemma is an `Arm' and $\delta$ is a `Leg', see~\cite[page~13]{GoodwinLeMagaardPaolini16}. 
Proposition~\ref{P:gdid31}b) reduce the classification of the coadjoint orbits of $\nt_{X; Z}^*$ to the classification of the coadjoint orbits of $\nt_{X\setminus\{\alpha, \beta\};Z(X)}^*$. We call such a reduction an \emph{AL-move}. 

Overall, Proposition~\ref{P:gdid31} is a version of~\cite[Lemma~3.1]{GoodwinLeMagaardPaolini16} for coadjoint orbits. 
And the fastest proof here is a reduction to~\cite[Lemma~3.1]{GoodwinLeMagaardPaolini16} by a version of the Kirillov orbit method.  
Nevertheless, it is not hard to prove it directly using standard to coadjoint orbits approach and we present such a proof here.


\begin{proof}{~of~Proposition~\ref{P:gdid31}}{ 
Part a). The conditions on $\beta$ implies that $\nt_{X\setminus\{\beta\}}$ is a subalgebra of $\nt_X$ of codimension 1 and thus it is normal in $\nt_X$.
Moreover $e_{\delta}$ belongs to the center of $\nt_{X\setminus\{\beta\}}$ and it does not belong to the center of $\nt_X$. 

Pick $f\in\nt_{X; Z(X)}^*$. There are two options for the map $\psi: N_X.f\to\nt_{X\setminus\{\beta\}}^*$, see~\cite[Remark~2.9]{IgnatyevPetukhov21}:

1) $\psi(N.f)$ is a single $N_{X\setminus \{\beta\}}$-orbit and $\psi$ is an isomorphism between $N_X.f$ and $\psi(N.f)$;

2) $\dim N_X.f=\dim \psi(N_X.f)+1, N_X.f=\psi^{-1}(\psi(N_X.f))$ and $\psi(N_X.f)$ splits into a one-dimensional family of $N_{X\setminus \beta}$-orbits with each orbit being of dimension $\dim N.f-2$.

Taking into account that $[e_\delta, e_\beta]$ is proportional to $e_\gamma$ with nonzero scalar and $[e_\gamma, e_\beta]=0$, we have 
\begin{equation}\label{Ealm} (\exp{(te_\delta)}.f)(e_\beta)=f(e_\delta)-t f([e_\delta, e_\beta]),~t\in\Fp.\end{equation}
This implies that $e_\delta$ can take all possible values on $\psi(N_X.f)$ and hence it is not constant. Thus, $\psi(N_X.f)$ can not be a single $N_{X\setminus\{\beta\}}$-orbit and we are in case 2).

In this case we have $\dim \psi(N_X.f)=\dim N_X.f-1$. Note that $\psi(N_X.f)$ is a single $N_X$-orbit and hence it is closed and irreducible and is even isomorphic to an affine space. Thus the divisors defined by the equation $e_\gamma(\cdot)=c$ are of dimension $\dim N_X.f-2$ for all $c\in\Fp$. 
Thanks to the fact that in case 2) all $N_{X\setminus\{\beta\}}$-orbits of $\psi(N_X.f)$ are of dimension $\dim N_X.f-2$, each fiber of $\psi$ is a union of several orbits. Formula (\ref{Ealm}) implies that there is a projection from $\psi(N_X.f)$ to the zero-locus of $e_\delta$ and this implies that such a zero-locus is irreducible. Hence it is a single $N_{X\setminus\{\beta\}}$-orbit. 
We left to mention that $N_{X\setminus\{\beta\}}$-orbits on $\nt_{X\setminus\{\beta\}; Z}^*$-orbits for which $e_\delta=0$ can be identified with $N_{X\setminus \{\delta, \beta\}}$-orbits on $\nt_{X\setminus \{\delta, \beta\}; Z}^*$.

The above discussion implies that the map
\begin{equation}\label{Ealm2}N_X.f\to\text{zero-locus~of~}e_\delta{\rm~on~}\psi(N_X.f)\end{equation}
provides a bijection between $N_X$-orbits on $\nt_{X; Z}^*$ and $N_{X\setminus \{\delta, \beta\}}$-orbits on $\nt_{X\setminus\{\delta, \beta\} ; Z}^*$.

Part b) follows from the explicit form (\ref{Ealm2}) of the bijection of part a).

To prove part c) we note that $\mathfrak p\oplus\Fp e_\delta$ is a polarization of $\phi(\lambda)|_{\nt_{X\setminus\{\beta\}}}$ and thus it will be also a polarization for $\lambda$ by dimension reasons. 
}\end{proof}
\sst{Proof~of Proposition~\ref{Psupp} for types $B$ and $C$}
\begin{proof}{~of Proposition~\ref{Psupp}, type $C$}{
As in types A/D, it is enough to show that there is no $N.f$ such that $\dim N.f=2d$ and $\Supp{f}\not\subset\SuppD_d$. 
Assume to the contrary that such $f\in\nt^*$ exists. 
As in types A/D, we choose a maximal root $\gamma\in \Supp{f}\setminus\mathbb \SuppD_d$. 
If $\gamma=\epsi_i-\epsi_j$ then $\alpha_1, \alpha_2, \beta_1, \beta_2, \gamma$ are short in $\Phi$. 
As in types $A/D$, we reduce this case to the straightforward analogue of Lemma~\ref{L:sumr3} with $\gamma=\epsi_i-\epsi_j$. 
We prove this analogue via a similar argument reducing it to a root systems of rank 3 or less. 
These root systems have irreducible components of types $A$ and $C$ only, because they are contained in root system of type $C_n$. 
Taking into account that $\alpha_1, \alpha_2, \beta_1, \beta_2, \gamma$ are short, it is easy to finish the proof by a straightforward check. 

We left to consider the cases $\gamma=2\epsi_i$, $i\le n$ and $\gamma=\epsi_i+\epsi_j$ for $i< j\le n$. 
Assume $\gamma=2\epsi_i$. 
For the subdiagram $C_{n+1-i}$ of the Dynkin diagram $C_n$, denote by $\mathfrak m$ the respective root subalgebra of $\nt$ as well as by $M$ the respective subgroup of $N$. 
We set $\underline{f}:=f|_{\frak m}$. 
It is clear that $\dim M.\underline{f}\le \dim N.f$ and $\gamma\in\NSupp{\underline{f}}$. 
This allows us to replace $n$ by $n+1-i$, $N$ by $M$, $\gamma$ by $2\varepsilon_1$ and $f$ by $\underline f$ respectively. 

Using AL-moves for the sequence of pairs of roots $(\beta, \delta)$ of the form 
$$(\epsi_1-\epsi_2, \epsi_1+\epsi_2),~(\epsi_1-\epsi_3, \epsi_1+\epsi_3),~\ldots,~(\epsi_1-\epsi_n, \epsi_1+\epsi_n),$$ we show that $\dim M.\underline{f}=2(n-1)+d'$ where $d'$ is the dimension of an orbit of a certain root subalgebra of $\frak m$. This clearly implies that $\dim N.f\ge |\Sing(\gamma)|=2(n-1)$, as needed in this case. 

We proceed to the case $\gamma=\epsi_i+\epsi_j$ for $i< j\le n$. 
As before, consider the subdiagram $C_{n+1-i}$ of the Dynkin diagram $C_n$ and denote by $\frak m$ the respective root subalgebra of $\nt$ as well as by $M$ the respective subgroup of $N$. 
We set $\underline{f}:=f|_{\frak m}$. 
It is clear that $\dim M.\underline{f}\le \dim N.f$ and $\gamma\in\NSupp{\underline{f}}$. 
This allows as to replace $n$ by $n+1-i$, $N$ by $M$, $j$ by $j+1-i$, $\gamma$ by $\varepsilon_1+\varepsilon_{j}, 1<j,$ and $f$ by $\underline f$ respectively. 

Using AL-moves for the sequence of pairs of roots $(\beta, \delta)$ of the form 
$$(\epsi_j-\epsi_{j+1}, \epsi_1+\epsi_{j+1}),~(\epsi_j-\epsi_{j+2}, \epsi_1+\epsi_{j+2}),~\ldots,~(\epsi_j-\epsi_n, \epsi_1+\epsi_n),$$
and then for the sequence
$$(\epsi_j+\epsi_{n}, \epsi_1-\epsi_{n}),~(\epsi_j+\epsi_{n-1}, \epsi_1-\epsi_{n-1}),~\ldots,~(\epsi_j+\epsi_{j+1}, \epsi_1-\epsi_{j+1}),$$
and then for pair $(2\epsi_j, \epsi_1-\epsi_j)$, and then for the sequence
$$(\epsi_{j-1}+\epsi_{j}, \epsi_1-\epsi_{j-1}),~(\epsi_{j-2}+\epsi_{j}, \epsi_1-\epsi_{j-2}),~\ldots,~(\epsi_2+\epsi_{j}, \epsi_1-\epsi_{2}),$$
we show that $\dim M.\underline{f}=2(n-j+n-1)+d'$ where $d'$ is the dimension of an orbit of a certain root subalgebra of $\frak m$. This clearly implies that $\dim N.f\ge |\Sing(\gamma)|=2(2n-j-1)$, as required.}
\end{proof}

\begin{proof}{~of Proposition~\ref{Psupp}, type $B$}{
The proof for type $B$ is a mixture of proofs for types $A/D$ and $C$. 

As in the other types it is enough to show that there is no $N.f$ such that $\dim N.f=2d$ and $\Supp{f}\not\subset\SuppD_d$. 
Assume to the contrary that such $f\in\nt^*$ exists. 
As in the other types we choose a maximal root $\gamma\in \Supp{f}\setminus\SuppD_d$. 
If $\gamma=\epsi_i-\epsi_j$ then $\alpha_1, \alpha_2, \beta_1, \beta_2, \gamma$ are long in $\Phi$. 
As in types $A/D$ we reduce this case to the straightforward analogue of Lemma~\ref{L:sumr3} with $\gamma=\epsi_i-\epsi_j$. 
We prove this analogue via a similar argument reducing it to a root systems of rank 3 or less. 
These root systems have irreducible components of types $A$ and $B$ only, because they are contained in root system of type $B_n$. 
Taking into account that $\alpha_1, \alpha_2, \beta_1, \beta_2, \gamma$ are long it is easy to finish the proof by a straightforward check. 

We left to consider the cases $\gamma=\epsi_i, i\le n$ and $\gamma=\epsi_i+\epsi_j$ for $i< j\le n$. 
Assume $\gamma=\epsi_i$. 
For the subdiagram $B_{n+1-i}$ of the Dynkin diagram $B_n$ denote by $\frak m$ the respective root subalgebra of $\nt$ as well as by $M$ the respective subgroup of $N$. 
We set $\underline{f}:=f|_{\frak m}$. 
It is clear that $\dim M.\underline{f}\le \dim N.f$ and $\gamma\in\NSupp{\underline{f}}$. 
This allows as to replace $n$ by $n+1-i$, $N$ by $M$, $\gamma$ by $\varepsilon_1$ and $f$ by $\underline f$ respectively. 

Using AL-moves for the sequence of pairs of roots $(\beta, \delta)$ of the form 
$$(\epsi_n, \epsi_1-\epsi_n),~(\epsi_{n-1}, \epsi_1-\epsi_{n-1}),~\ldots,~(\epsi_2, \epsi_1-\epsi_2),$$ we show that $\dim M.\underline{f}=2(n-1)+d'$ where $d'$ is the dimension of an orbit of a certain root subalgebra of $\frak m$. This clearly implies that $\dim N.f\ge |\Sing(\gamma)|=2(n-1)$ as needed in this case. 

We proceed to the case $\gamma=\epsi_i+\epsi_j$ for $i< j\le n$. 
As before consider the subdiagram $B_{n+1-i}$ of the Dynkin diagram $B_n$ and denote by $\frak m$ the respective root subgroup of $\nt$ as well as by $M$ the respective subgroup of $N$. 
We set $\underline{f}:=f|_{\frak m}$. 
It is clear that $\dim M.\underline{f}\le \dim N.f$ and $\gamma\in\NSupp{\underline{f}}$. 
This allows as to replace $n$ by $n+1-i$, $N$ by $M$, $j$ by $j+1-i$, $\gamma$ by $\varepsilon_1+\varepsilon_{j}, 1<j,$ and $f$ by $\underline f$ respectively. 

Using AL-moves for the following sequence of pairs of roots
$$(\epsi_j-\epsi_{j+1}, \epsi_1+\epsi_{j+1}),~(\epsi_j-\epsi_{j+2}, \epsi_1+\epsi_{j+2}),~\ldots,~(\epsi_j-\epsi_n, \epsi_1+\epsi_n),$$
and then for pair $(\epsi_j, \epsi_1)$, and then for the following sequence$$(\epsi_j+\epsi_{n}, \epsi_1-\epsi_{n}),~(\epsi_j+\epsi_{n-1}, \epsi_1-\epsi_{n-1}),~\ldots,~(\epsi_j+\epsi_{j+1}, \epsi_1-\epsi_{j+1}),$$
we show that $\dim M.\underline{f}=2(n-j+1+n-j)+d'$ where $d'$ is the dimension of an orbit of $\nt_{X_j}$ with $X_j$ being $\Phi^+$ minus all the above roots used for AL-moves and minus all the roots of the form $\epsi_1-\epsi_i$ with $i>j$. 

We left to show that $d'\ge 2(j-1)$. To do this we use an analogue of the argument from types $A/D$. 
Indeed, we consider the minor of $B_f$ standing for basis vectors $e_\alpha, \alpha\in \Sing(\gamma)$ with $\alpha$ being
\begin{equation}\epsi_1-\epsi_2,~\epsi_1-\epsi_3,~\ldots,~\epsi_1-\epsi_{j-1}\label{E:estB1}\end{equation}
and
\begin{equation}\epsi_{2}+\epsi_j,~\epsi_3+\epsi_j,~\ldots,~\epsi_{j-1}+\epsi_j.\label{E:estB2}\end{equation}

We claim that this minor splits into $j-1$ two-by-two non-zero skew-symmetric blocks defined by pairs $\alpha, \beta\in\Phi^+$ with $\alpha+\beta=\epsi_1+\epsi_j$. The rank of each such a block is equal to 2 and thus the rank of the whole minor is $2(j-1)$. This implies $d'\ge2(j-1)$, as needed. 

To prove the claim we note that a sum of two vectors of~(\ref{E:estB1}) is of the form [$2\epsi_1-$(other $\epsi_*$)]. 
There are no root of $B_n$ of this form. 
A sum of two vectors of~(\ref{E:estB2}) is of the form [$2\epsi_j+$(other $\epsi_*$)]. 
There are no root of $B_n$ of this form. 
A sum of a vectors of~(\ref{E:estB1}) with a vector of~(\ref{E:estB2}) is of the form [$\epsi_1+\epsi_j+$(other $\epsi_*$)-(other $\epsi_*$)]. 
The only root of $B_n$ of this form is $\epsi_1+\epsi_j$ as required.} 
\end{proof}

\sect{Orbits of dimension 2}\label{S:dim2}\addcontentsline{toc}{subsection}{\ref{S:dim2}. Orbits of dimension 2}
In this section we would like to enlist all extensive orbits of dimension 2. 
The technique developed for the case of dimension 2 will be used for dimensions 4 and 6 with some improvements. 

Thanks to Proposition~\ref{PlWd} the only simple algebras for which there is a chance to find extensive orbits of dimension 2 are $A_2, A_3, B_2=C_2, B_3, C_3, D_4, G_2$ (and they all indeed allow such orbits). 
The Dynkin diagram of type $G_2$ can not occur as a subdiagram of the Dynkin diagrams of the classical series and thus the classification of orbits of dimension 2 for the classical series does not involve this exceptional case.  
For each such a case we provide a list of subsets $S_1, S_2,\ldots$ of $\Phi^+$ such that each extensive coadjoint orbit of dimension 2 in $\nt^*$ intersects $V(S_1)\sqcup V(S_2)\sqcup\ldots$ in exactly one point. 

All the orbits of $\nt^*$ for Lie algebras $A_2, A_3, B_2=C_2, B_3, C_3, D_4$ were classified before, see~\cite{GoodwinLeMagaardPaolini16} for more details. Moreover, it is not hard to dig out set-sections for the respective description, but not all irreducible components of these set-sections will be of the form $V(S)$. 
The same holds for the orbits of dimension 4 and even for most of the cases of dimension 6 with the only exception for $D_8$ (up to our knowledge, the classification of orbits for this Lie algebra is not achieved yet). 

The above sets $S_1, S_2, \ldots$ are encoded in this article by a string of length $|\Phi^+|$ consisting of letters `S', `A', `I', `L'. 
The $j$th root in the enumeration convention of CHEVIE package (see also Section~\ref{S:chv}) belongs to $S_i$ if and only if the $j$th letter of $i$th string is `S'. Other letters encode additional data about the corresponding orbits. 
For example, twice the amount of letters `A' equals the dimension of our orbit; it is also equal to twice the amount of letters `L' for each string. 
The letter `I' means that the corresponding position can be `Ignored'.
This notation with strings and sets will be used also for dimensions 4 and 6. 
\mtheo{\label{Tdim2}The following strings enlists \textup(in the above sense\textup) all extensive orbits of dimension \textup2 for Lie algebras of types $A_2, A_3, B_2=C_2, B_3, C_3, D_4$.}
\begin{longtable}{|c|c|c|c|c|c|}
\hline \multicolumn{6}{|c|}{Table~\label{T:dim2all}\ref{T:dim2all}}\\
\hline Root system& String&Root system& String&Root system& String\\
\hline $A_2$& LAS&$C_3$&SLASSSIII&$D_4$&ASLSSSSIIIII\\
\hline $A_3$& ALSSSI&---&SLAISSIII&---&ASLISSSIIIII\\
\hline ---& ALISSI&---& ILASSSIII&---&AILSSSSIIIII\\
\hline $B_2$& LSAS&---&ILAISSIII&---&AILISSSIIIII\\
\hline ---& LIAS&---&ALSSSIIII&&\\
\hline ---& LASI&---&ALISSIIII&&\\
\hline $B_3$& ALSSSIIII&& &&\\
\hline ---& ALISSIIII&& &&\\
\hline
\end{longtable}
Note also that an orbit is elementary if and only if the corresponding string has exactly one letter~`S'. 

To prove Theorem~\ref{Tdim2}, we need to discuss in detail the meaning of the above strings together with the underlying linear algebra.
\sst{The algorithm for dimension 2}
We will now show how one can get the desired strings of Theorem~\ref{Tdim2}. 
As a general concept: we split the entire set of coadjoint orbits into some pieces and simplify each piece as much as we can. 
This algorithm is based on the algorithm given in~\cite{GoodwinLeMagaardPaolini16} and can be considered as a minor modification of it.

Each string $\mathcal S$ will consist of letters `S', `A', `I', `L' and `Q'. For $\alpha\in\Phi^+$ we denote by $\mathcal S|_\alpha$ the letter used on the position attached to $\alpha$. 
Each string $\mathcal S$ corresponds to $\nt_{S; Z}^*$ where \begin{center}$S:=\{\alpha\in\Phi^+\mid \mathcal S|_\alpha=~\mbox{`S'}\mbox{~or~}~\mbox{`Q'}\}$ and $Z:=\{\alpha\in\Phi^+\mid \mathcal S|_\alpha=~\mbox{`S'}\}$.\end{center} It is easy to verify that this notation is consistent: $S$ is a C-quattern and $Z\subset Z(S)$ for all the strings which show up in our algorithm. 

We now describe the algorithm step-by-step.

{\bf Step 1:} consider the set $\SuppD_2$ for each root system in question defined in Theorem~\ref{Tnso} and attach to it a string $MaskD$ consisting of letters `Q' and `I', where letters `Q' (`Question') are used on  the $i$th  position if the $i$th root of $\Phi^+$ belong to $\SuppD_2$ and letter `I' is used otherwise (`Ignored'). 
If $2>|\Sing(\delta_0)|$ in the notation of Theorem~\ref{Tnso} then $MaskD$ set to be the string of `Q' letters only. 

The set $\SuppD_2$ is by definition a large C-quattern and moreover the map $\nt_{\SuppD_2}^*\to\nt^*$ identifies 2-dimensional coadjoint orbits of $\nt^*$ and of $\nt_{\SuppD_2}^*$. 
From this step we consider only the coadjoint orbits of $\nt_{\SuppD_2}^*$. 

{\bf Example for $A_3$:} the string attached to $\SuppD_2$ equals `QQQQQI'. 

{\bf Step 2 (recursive):} during this step we will produce out of $MaskD$ a list of strings which encodes all sets $S=\NSupp{f}\subset\SuppD_2$ for which there is a chance that there exists an extensive orbit of dimension 2 with the given $S=\NSupp{f}$. 
Note that in this case $N.f$ can be identified with the 2-dimensional orbit of $\nt_{\NSupp{f}}^*$. 
Thus every case corresponds to the 2-dimensional coadjoint orbits for a separate Lie algebra.

To each string $\mathcal S$ produced on this step we attach $S, Z$ as above. 
Next, we evaluate $Z(S)$. If $Z=Z(S)$ then we proceed to other strings. If all strings satisfy this condition then we finish Step 2. 

In another case we apply the procedure which we call {\em I/S-move}. This is an analogue of Type~S reduction of~\cite{GoodwinLeMagaardPaolini16}. 
Pick $\alpha\in Z(S)\setminus Z$ and add two new strings: one in which $\mathcal S|_\alpha=$`Q' is replaced by `S' and in another one the same `Q' is replaced by `I'. The formula behind this procedure is as follows: 
$$\underline{\nt}_{S; Z}^*=\underline{\nt}_{S; Z\cup\{\alpha\}}^*\sqcup\underline{\nt}_{S\setminus\{\alpha\}; Z\setminus\{\alpha\}}^*.$$
Recall that $\nt_{S; Z(S)}^*$ coincides with the set of orbits such that $\NSupp{f}=S$ for each large C-quattern~$S$. 
Note that $S$ will be a C-quattern on this step.  
Note that the amount of letters `Q' in an `active' string is reduced by 1 at each step and thus this recursive procedure will terminate at some point.

Finally we get rid of all $\mathcal S$ with non-extensive $S$: i.e., we get rid of those for which the sums of adjacent simple roots are not in $S$. 

We save the result to the list $LrgQt$. 

{\bf Example for $A_3$:} The result for $A_3$ is as follows.
$$\NewAdigraph{MyGraph}{
 QQQQQI, black, 4: 4, 0;
 QQQSQI, black, 4: 3, -2;
 QQQIQI, black, 1: 5, -2;
 QQQSSI, black, 4: 2, -4;
 QQQSII, black, 1: 4, -4;
 non-ext, red, 0.1: 6, -4;
}
{
 QQQQQI, QQQSQI;
 QQQQQI, QQQIQI;
 QQQSQI, QQQSSI;
 QQQSQI, QQQSII;
 QQQIQI, non-ext;
 QQQSII, non-ext;
}
\MyGraph{}$$
Here we start from the string `QQQQQI' and then proceed to a string with one less `Q' on each level. 
Strings `QQQIQI' and `QQQSII' have `I' on the forth or fifth position and this means that the respective set $S$ does not contain at least one of the roots $$\epsi_1-\epsi_3=(\epsi_1-\epsi_2)+(\epsi_2-\epsi_3),~\epsi_2-\epsi_4=(\epsi_2-\epsi_3)+(\epsi_3-\epsi_4).$$ 
Thus all the subsequent strings will correspond to non-extensive orbits (and hence we ignore these cases).

As the final output the list $LrgQt$ consists of only one string in this case: `QQQSSI'. Other strings either stand for non-extensive orbits or can be processed via an I/S-move. 

The tables standing for strings `QQQQQI'$\to$`QQQSQI'$\to$`QQQSSI':
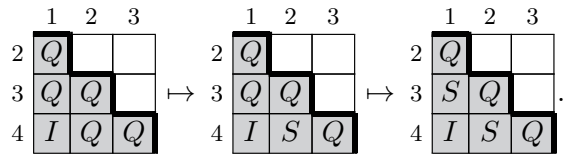
$$
{\Autonumfalse\mymatrix{
\lNote{2}Q\Note{1}\gray\Top{2pt}\Rt{2pt}&\Note{2}\Note{2}\Bot{2pt}\pho&\Note{3}\pho\\
\lNote{3}Q\gray&Q\gray \Rt{2pt}&\Bot{2pt}\\
\lNote{4}I\gray&Q\gray&Q\gray \Rt{2pt}\\
}}\mapsto
{\Autonumfalse\mymatrix{
\lNote{2}Q\Note{1}\gray\Top{2pt}\Rt{2pt}&\Note{2}\Note{2}\Bot{2pt}\pho&\Note{3}\pho\\
\lNote{3}Q\gray&Q\gray \Rt{2pt}&\Bot{2pt}\\
\lNote{4}I\gray&S\gray&Q\gray \Rt{2pt}\\
}}
\mapsto
{\Autonumfalse\mymatrix{
\lNote{2}Q\Note{1}\gray\Top{2pt}\Rt{2pt}&\Note{2}\Note{2}\Bot{2pt}\pho&\Note{3}\pho\\
\lNote{3}S\gray&Q\gray \Rt{2pt}&\Bot{2pt}\\
\lNote{4}I\gray&S\gray&Q\gray \Rt{2pt}\\
}}.$$

{\bf Step 3.1 (recursive):} on this step we work separately with each string $\mathcal S$ from Step 2. 
We attach $S, Z$ to it as before. 
To each string obtained on this step we apply 1) I/S-moves as much as we can 2) AL-moves as much as we can.

On each AL-move we replace a string $\mathcal S$ by a string $\mathcal S'$ in which $\mathcal S|_\beta$ is replaced by `A' and $\mathcal S|_\delta$ is replaced by `L'. 
Also for the new data $\mathcal S'$ we need to describe orbits of dimension decreased by 2 for each letter `A' in $\mathcal S'$ (equivalently, for each letter `L').  

Note that the amount of letters `Q' in an `active' string is reduced by 1 with an I/S-move and by 2 with an AL-move on each step and thus this recursive procedure will terminate at some point.




If a list of strings produced out of string $\mathcal S$ contains no letters `Q' (i.e., all letters `Q' were eliminated on this step) then we call the string $\mathcal S$ from $LrgQt$ {\em abelian} and save it to the list $LrgAQt$. 
In the opposite case we call $\mathcal S$ {\it non-abelian} and save it to $LrgNAQt$. 
Note that in any abelian case we essentially provide a classification of all orbits for this case and not only of orbits of dimension 2.

We save the list of strings produced out of $\mathcal S$ in the cases of $LrgNAQt$ to $ClsPat$ (it stands for a set-section of the whole set of orbits). We copy elements of $ClsPat$ with exactly one letter `A' to the list $ClsPatD$ (this stands for the list of orbits of dimension 2). 

{\bf Example for $C_3$ and $\mathcal S=\mbox{`QQQQSSSIII'}$:} The result is as follows:
$$\mbox{`QQQQQSSII'}\mapsto\mbox{`AQQQLSSII'}\mapsto\mbox{`ALQALSSII'}\mapsto\begin{array}{c}\mbox{`ALSALSSII'}\\ \mbox{`ALIALSSII'}\end{array},
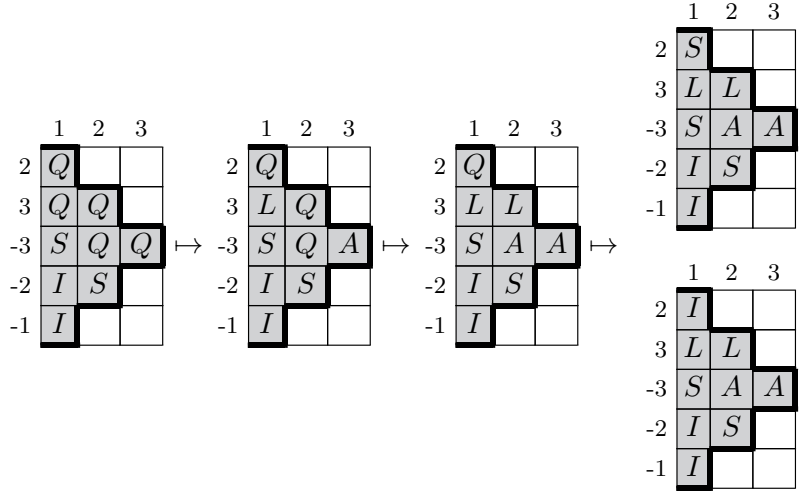
$$
$$
{\Autonumfalse\mymatrix{
\lNote{2}Q\Note{1}\gray\Top{2pt}\Rt{2pt}&\Note{2}\Note{2}\Bot{2pt}\pho&\Note{3}\pho\\
\lNote{3}Q\gray&Q\gray \Rt{2pt}&\Bot{2pt}\\
\lNote{-3}S\gray&Q\gray&Q\gray \Rt{2pt}\Bot{2pt}\\
\lNote{-2}I\gray&S\gray\Bot{2pt}\Rt{2pt}&\\
\lNote{-1}I\gray\Bot{2pt}\Rt{2pt}&& \\
}}\mapsto
{\Autonumfalse\mymatrix{
\lNote{2}Q\Note{1}\gray\Top{2pt}\Rt{2pt}&\Note{2}\Note{2}\Bot{2pt}\pho&\Note{3}\pho\\
\lNote{3}L\gray&Q\gray \Rt{2pt}&\Bot{2pt}\\
\lNote{-3}S\gray&Q\gray&A\gray \Rt{2pt}\Bot{2pt}\\
\lNote{-2}I\gray&S\gray\Bot{2pt}\Rt{2pt}&\\
\lNote{-1}I\gray\Bot{2pt}\Rt{2pt}&& \\
}}\mapsto
{\Autonumfalse\mymatrix{
\lNote{2}Q\Note{1}\gray\Top{2pt}\Rt{2pt}&\Note{2}\Note{2}\Bot{2pt}\pho&\Note{3}\pho\\
\lNote{3}L\gray&L\gray \Rt{2pt}&\Bot{2pt}\\
\lNote{-3}S\gray&A\gray&A\gray \Rt{2pt}\Bot{2pt}\\
\lNote{-2}I\gray&S\gray\Bot{2pt}\Rt{2pt}&\\
\lNote{-1}I\gray\Bot{2pt}\Rt{2pt}&& \\
}}\mapsto
\begin{array}{c}
{\Autonumfalse\mymatrix{
\lNote{2}S\Note{1}\gray\Top{2pt}\Rt{2pt}&\Note{2}\Note{2}\Bot{2pt}\pho&\Note{3}\pho\\
\lNote{3}L\gray&L\gray \Rt{2pt}&\Bot{2pt}\\
\lNote{-3}S\gray&A\gray&A\gray \Rt{2pt}\Bot{2pt}\\
\lNote{-2}I\gray&S\gray\Bot{2pt}\Rt{2pt}&\\
\lNote{-1}I\gray\Bot{2pt}\Rt{2pt}&& \\
}}\\
~\\
{\Autonumfalse\mymatrix{
\lNote{2}I\Note{1}\gray\Top{2pt}\Rt{2pt}&\Note{2}\Note{2}\Bot{2pt}\pho&\Note{3}\pho\\
\lNote{3}L\gray&L\gray \Rt{2pt}&\Bot{2pt}\\
\lNote{-3}S\gray&A\gray&A\gray \Rt{2pt}\Bot{2pt}\\
\lNote{-2}I\gray&S\gray\Bot{2pt}\Rt{2pt}&\\
\lNote{-1}I\gray\Bot{2pt}\Rt{2pt}&& \\
}}
\end{array}.$$
All letters `Q' were eliminated on all the paths and thus `QQQQSSSIII' is abelian. 

As the output we add string `QQQQSSSIII' to $LrgAQt$, we add strings `ALSALSSII', `ALIALSSII' to ClsPatD and add nothing to $LrgNAQt$.

The good news for dimension 2 is that for all Lie algebras in question we get no output to $LrgNAQt$. This is not the case for dimension 4 and 6 and thus the list of strings ClsPatD provides the desired list $S_1, S_2,\ldots$ thanks to Lemma~\ref{L:tpS} and Proposition~\ref{P:gdid31}b) applied recursively for the case of dimension~2.

\sect{Orbits of dimension 4}\label{S:dim4}\addcontentsline{toc}{subsection}{\ref{S:dim4}. Orbits of dimension 4}
In this section we would like to enlist all extensive orbits of dimension 4. The technique is very similar to the technique developed for dimension 2. 

Thanks to Proposition~\ref{PlWd} the only algebras for which there is a chance to find extensive orbits of dimension 2 are $A_2, A_3, A_4, A_5, B_2=C_2, B_3, B_4, B_5, C_3, C_4, C_5, D_4, D_5, D_6, G_2, F_4$ (most of them indeed allow such orbits with the exceptions of $A_2, B_2=C_2$ --- they are too small). 
The Dynkin diagrams of types $G_2, F_4$ can not occur as subdiagrams of the Dynkin diagrams of the classical series and hence the classification of orbits of dimension 4 for the classical series does not involve these exceptional cases.  

For each such a case we provide a list of subsets $S_1, S_2,\ldots$ of $\Phi^+$ such that each extensive coadjoint orbit of dimension 4 in $\nt^*$ intersects $V(S_1)\sqcup V(S_2)\sqcup\ldots$ in exactly one point. 
\theo{\label{Tdim4} The following strings enlists \textup(in the above sense\textup) all extensive orbits of dimension~\textup4 for Lie algebras of types $A_2, A_3, A_4, A_5,  B_2, B_3, B_4, B_5, C_3, C_4, C_5, D_4, D_5, D_6$.}{}
\begin{longtable}{|c|l|c|l|c|l|}
\hline \multicolumn{6}{|c|}{Table~\label{T:dim4P1}\ref{T:dim4P1}}\\
\hline Type& String&Type& String&Type& String\\
\hline $A_2$& ---&$B_2$& ---&$D_4$&LASSLAISSIII\\
\hline $A_3$& LSLAAS&$C_2$&---&---&LASILAISSIII\\
\hline ---& LILAAS&$B_3$&LALSSSASI&---&LAISLAISSIII\\
\hline $A_4$& LALASSSSII&---&LALSSIASI&---&LAIILAISSIII\\
\hline ---& LALASISSII&---&LSLSIAASI&---&ALSSALISISII\\
\hline ---& LSLAIASSII&---&LILSIAASI&---&ALSIALISISII\\
\hline ---& LILAIASSII&---&LALISSASI&---&ALISALISISII\\
\hline ---& ALALSSSISI&---&LALISIASI&---&ALIIALISISII\\
\hline ---&ALALSISISI&---&LSLIIAASI&---&ASSLAILISSII\\
\hline ---&ALSLSAIISI&---&LILIIAASI&---&AISLAILISSII\\
\hline ---&ALILSAIISI&---&LSLAASSII&---&ASILAILISSII\\
\hline ---&ALLASSSIII&---&LILAASSII&---&AIILAILISSII\\
\hline $A_5$&ALALSSSSSISIIII&---&LLAASSIII&&\\
\hline ---&ALALSSISSISIIII&---&LSLAAISII&&\\
\hline ---&ALALISSSSISIIII&---&LILAAISII&&\\
\hline ---&ALALISISSISIIII&$C_3$&SSLILIAAS&&\\
\hline ---&ALSLASSISISIIII&---&ISLILIAAS&&\\
\hline ---&ALSLASIISISIIII&---&SILILIAAS&&\\
\hline ---&ALILASSISISIIII&---&IILILIAAS&&\\
\hline ---&ALILASIISISIIII&---&LASLASSII&&\\
\hline ---&ALALSSSSSIIIIII&---&LAILASSII&&\\
\hline ---&ALALISSSSIIIIII&---&LSLAAISII&&\\
\hline &&---&LILAAISII&&\\
\hline
\end{longtable}

\begin{longtable}{|c|l|c|l|}
\hline \multicolumn{4}{|c|}{Table~\label{T:dim4P2}\ref{T:dim4P2}}\\
\hline Type& String&Type& String\\
\hline $B_4$&ALALSSSIISIIIIII&$C_4$&SLALASSSISIIIIII\\
\hline ---&ALALSISIISIIIIII&---&SLALAISSISIIIIII\\
\hline ---&ALSLSAIIISIIIIII&---&SLSLAAISISIIIIII\\
\hline ---&ALILSAIIISIIIIII&---&SLILAAISISIIIIII\\
\hline ---&LSLAASSSSISIIIII&---&ILALASSSISIIIIII\\
\hline ---&LSLAASSSIISIIIII&---&ILALAISSISIIIIII\\
\hline ---&LSLAAISSSISIIIII&---&ILSLAAISISIIIIII\\
\hline ---&LSLAAISSIISIIIII&---&ILILAAISISIIIIII\\
\hline ---&LILAASSSSISIIIII&---&ALALSSSIISIIIIII\\
\hline ---&LILAASSSIISIIIII&---&ALALSISIISIIIIII\\
\hline ---&LILAAISSSISIIIII&---&ALSLSAIIISIIIIII\\
\hline ---&LILAAISSIISIIIII&---&ALILSAIIISIIIIII\\
\hline ---&LALASSSISISIIIII&---&SLLAASSSIIIIIIII\\
\hline ---&LALASSSIIISIIIII&---&ILLAASSSIIIIIIII\\
\hline ---&LALASISISISIIIII&---&LALASSSISIIIIIII\\
\hline ---&LALASISIIISIIIII&---&LALASISISIIIIIII\\
\hline ---&LSLAISSIAISIIIII&---&LSLAIASISIIIIIII\\
\hline ---&LSLAIISIAISIIIII&---&LILAIASISIIIIIII\\
\hline ---&LILAISSIAISIIIII&---&ALLASSSIIIIIIIII\\
\hline ---&LILAIISIAISIIIII&---&SSLSILSIAIAISIII\\
\hline ---&LSLAASSSSIIIIIII&---&ISLSILSIAIAISIII\\
\hline ---&LSLAAISSSIIIIIII&---&SILSILSIAIAISIII\\
\hline ---&LILAASSSSIIIIIII&---&IILSILSIAIAISIII\\
\hline ---&LILAAISSSIIIIIII&---&SSLIILSIAIAISIII\\
\hline ---&LSLAASSSIIIIIIII&---&ISLIILSIAIAISIII\\
\hline ---&LILAASSSIIIIIIII&---&SILIILSIAIAISIII\\
\hline ---&LALASSSISIIIIIII&---&IILIILSIAIAISIII\\
\hline ---&LALASISISIIIIIII&&\\
\hline ---&LSLAIASISIIIIIII&&\\
\hline ---&LILAIASISIIIIIII&&\\
\hline ---&ALLASSSIIIIIIIII&&\\
\hline $B_5$&ALALSSSSSIISIIIIIIIIIIIII&&\\
\hline ---&ALALSSISSIISIIIIIIIIIIIII&&\\
\hline ---&ALALISSSSIISIIIIIIIIIIIII&&\\
\hline ---&ALALISISSIISIIIIIIIIIIIII&&\\
\hline ---&ALSLASSISIISIIIIIIIIIIIII&&\\
\hline ---&ALSLASIISIISIIIIIIIIIIIII&&\\
\hline ---&ALILASSISIISIIIIIIIIIIIII&&\\
\hline ---&ALILASIISIISIIIIIIIIIIIII&&\\
\hline ---&ALALSSSSSIIIIIIIIIIIIIIII&&\\
\hline ---&ALALISSSSIIIIIIIIIIIIIIII&&\\
\hline
\end{longtable}

\begin{longtable}{|c|l|c|l|}
\hline \multicolumn{4}{|c|}{Table~\label{T:dim4P3}\ref{T:dim4P3}}\\
\hline Type& String&Type& String\\
\hline $C_5$&SLALASSSSSISIIIIIIIIIIIII&$D_5$&ASLALSSSSIIISIIIIIII\\
\hline ---&SLALASSISSISIIIIIIIIIIIII&---&ASLALSSISIIISIIIIIII\\
\hline ---&SLALAISSSSISIIIIIIIIIIIII&---&ASLSLSSAIIIISIIIIIII\\
\hline ---&SLALAISISSISIIIIIIIIIIIII&---&ASLILSSAIIIISIIIIIII\\
\hline ---&ILALASSSSSISIIIIIIIIIIIII&---&AILALSSSSIIISIIIIIII\\
\hline ---&ILALASSISSISIIIIIIIIIIIII&---&AILALSSISIIISIIIIIII\\
\hline ---&ILALAISSSSISIIIIIIIIIIIII&---&AILSLSSAIIIISIIIIIII\\
\hline ---&ILALAISISSISIIIIIIIIIIIII&---&AILILSSAIIIISIIIIIII\\
\hline ---&SLSLAAISSSISIIIIIIIIIIIII&---&ASLLASSSSIIIIIIIIIII\\
\hline ---&SLSLAAIISSISIIIIIIIIIIIII&---&AILLASSSSIIIIIIIIIII\\
\hline ---&ILSLAAISSSISIIIIIIIIIIIII&---&ASSLSAILSISSIIIIIIII\\
\hline ---&ILSLAAIISSISIIIIIIIIIIIII&---&AISLSAILSISSIIIIIIII\\
\hline ---&SLILAAISSSISIIIIIIIIIIIII&---&ASILSAILSISSIIIIIIII\\
\hline ---&SLILAAIISSISIIIIIIIIIIIII&---&AIILSAILSISSIIIIIIII\\
\hline ---&ILILAAISSSISIIIIIIIIIIIII&---&ASSLIAILSISSIIIIIIII\\
\hline ---&ILILAAIISSISIIIIIIIIIIIII&---&AISLIAILSISSIIIIIIII\\
\hline ---&ALALSSSSSIISIIIIIIIIIIIII&---&ASILIAILSISSIIIIIIII\\
\hline ---&ALALSSISSIISIIIIIIIIIIIII&---&AIILIAILSISSIIIIIIII\\
\hline ---&ALALISSSSIISIIIIIIIIIIIII&$D_6$&ASLALSSSSSSIIISIIIIIIIIIIIIIII\\
\hline ---&ALALISISSIISIIIIIIIIIIIII&---&ASLALSSSISSIIISIIIIIIIIIIIIIII\\
\hline ---&ALSLASSISIISIIIIIIIIIIIII&---&ASLALISSSSSIIISIIIIIIIIIIIIIII\\
\hline ---&ALSLASIISIISIIIIIIIIIIIII&---&ASLALISSISSIIISIIIIIIIIIIIIIII\\
\hline ---&ALILASSISIISIIIIIIIIIIIII&---&ASLSLASSSISIIISIIIIIIIIIIIIIII\\
\hline ---&ALILASIISIISIIIIIIIIIIIII&---&ASLSLASSIISIIISIIIIIIIIIIIIIII\\
\hline ---&SLALASSSSSIIIIIIIIIIIIIII&---&ASLILASSSISIIISIIIIIIIIIIIIIII\\
\hline ---&SLALAISSSSIIIIIIIIIIIIIII&---&ASLILASSIISIIISIIIIIIIIIIIIIII\\
\hline ---&ILALASSSSSIIIIIIIIIIIIIII&---&AILALSSSSSSIIISIIIIIIIIIIIIIII\\
\hline ---&ILALAISSSSIIIIIIIIIIIIIII&---&AILALSSSISSIIISIIIIIIIIIIIIIII\\
\hline ---&ALALSSSSSIIIIIIIIIIIIIIII&---&AILALISSSSSIIISIIIIIIIIIIIIIII\\
\hline ---&ALALISSSSIIIIIIIIIIIIIIII&---&AILALISSISSIIISIIIIIIIIIIIIIII\\
\hline &&---&AILSLASSSISIIISIIIIIIIIIIIIIII\\
\hline &&---&AILSLASSIISIIISIIIIIIIIIIIIIII\\
\hline &&---&AILILASSSISIIISIIIIIIIIIIIIIII\\
\hline &&---&AILILASSIISIIISIIIIIIIIIIIIIII\\
\hline &&---&ASLALSSSSSSIIIIIIIIIIIIIIIIIII\\
\hline &&---&ASLALISSSSSIIIIIIIIIIIIIIIIIII\\
\hline &&---&AILALSSSSSSIIIIIIIIIIIIIIIIIII\\
\hline &&---&AILALISSSSSIIIIIIIIIIIIIIIIIII\\
\hline
\end{longtable}
\sst{The algorithm for dimension 4}\label{SS:adim4}
We will now show how one can get the desired strings of Theorem~\ref{Tdim4}. 
We extensively use the notation and the concepts of Section~\ref{S:dim2}. 
In particular, we work with strings $\mathcal S$ and sets $S, Z$ for all the needed root systems (they are enlisted in the statement of Theorem~\ref{Tdim4}). 

We first do the analogue of Step 1 of the proof of Theorem~\ref{Tdim2} for the respective root system with the following modifications: 1)
if $4>|\Sing(\delta_0)|$ in the above notation then $MaskD$ set to be the string of `Q' letters only; 2) we replace $\SuppD_2$ by $\SuppD_4$. 

Next we do the straightforward analogue of Steps 2 and 3.1 of the proof of Theorem~\ref{Tdim2} for the respective root system with the following modification: we add a string from $ClsPat$ to $ClsPatD$ if it contains exactly 2 letters `A'.

For $A_2, A_3, A_4, A_5, B_2, B_3, B_4, B_5, C_3$ we have that $LrgNAQt$ is empty and thus everything works as in the case of dimension 2. 

For $C_4, C_5, D_4, D_5, D_6$ we have that $LrgNAQt$ is not empty and we need to deal with it in the next steps.

{\bf Step 3.2 (recursive):} Note that AL-moves involve a choice of a pair of roots $(\beta, \delta)$ on each step and this choice does matter. It can lead to the following effect: a given string $\mathcal S$ is abelian with some such choices and non-abelian with many others. 
On Step 3.1 we use some choice of AL-moves (in a more or less random way) and most of the strings $\mathcal S\in LrgQt$ are abelian with respect to this~way. 

On Step 3.2 we pick a non-abelian string $\mathcal S\in LrgNAQt$ and check whether or not it is trully non-abelian: we are trying to apply AL-moves and I/S-moves in all possible ways and see whether or not it is abelian in at least one such a way. 
If it is the case we add string $\mathcal S$ to the list $LrgAQt$. 
We also save the list of strings produced out of $\mathcal S$ through ``an abelian path'' to $ClsPat$. 
We do it only for the first ``abelian path'' which we encounter. 
Then we copy strings from $ClsPat$ to $ClsPatD$ with exactly 2 letters `A'. 
We add $\mathcal S$ to the list $TrueNAQt$ in another case. 

It is clear that the approach of Step 3.2 is better then approach of Step 3.1. And one can try to apply it to all strings of Step 2. The issue here is that the approach of Step 3.2 is much more computationally expensive ``per string''.

The list $TrueNAQt$ is empty for $C_4$, $C_5$, $D_5$ and this finishes the consideration of these cases. 

{\bf Step 3.3:} We have the list $TrueNAQt$ of non-abelian strings and it can be so that some of them provide no orbits of dimension 4. We would like to find out at least some such strings using a certain procedure. We essentially repeat the operations of Step 3.2 with a few minor modifications. 

To each string we attach a nonnegative integer which provides a lower estimate for the dimensions of orbits available for $\nt_{S; Z}^*$. 

When we apply an I/S-move we pick the minimum for the estimates of the respective modified strings. It stands for the following concept: with an I/S-move we split the subset of orbits into two pieces and pick as a lower estimate for the whole subset the minimum of the estimates available for subsets. 

When we apply an AL-move we pick 2 plus the lower estimate of the respective modified string. It works fine because of Proposition~\ref{P:gdid31}. 

If we can apply an AL-move in several different ways we choose the best (=the largest) lower estimate among all the cases.

Finally, we need to explain which lower estimate we attach to strings to which we can not apply neither I/S-move nor AL-move. In this case we attach to this string 0 or 2: 0 if $\nt_{S; Z}^*$ contains 0-dimensional orbits and 2 otherwise (recall that all coadjoint orbits are even-dimensional). 

The 0-dimensional orbits are well-known: they coincide with $\nt_{S\setminus(S+S)}^*,$ where $S+S$ denotes the pointwise sum of sets of vectors. 

If $Z\cap (S+S)\ne{\varnothing}$ then $\nt_{S\setminus(S+S)}^*\cap \nt_{S; Z}^*=\{0\}$ and so all orbits of $\nt_{S; Z}^*$ are of dimension 2 or more. 
In the other case we have $\nt_{S\setminus(S+S)}^*\cap \nt_{S; Z}^*=\nt_{S\setminus(S+S); Z}^*$ and hence all orbits of $\nt_{S; Z}^*$ are of dimension 0 or more. 

Finally, we delete from $TrueNAQt$ all strings for which the lower estimate for the dimensions of the respective orbits is more than 4. 

After Step 3.3 we have that
\begin{center}$TrueNAQt$=$\begin{cases}\mbox{`QQQQQQQSSSII'}&\mbox{~for~}D_4\\\mbox{`QQQQQQQQQSSSSSIIIIIIIIIIIIIIII'}&\mbox{~for~}D_6\\{\varnothing}&\mbox{all~other~cases}\end{cases}$.\end{center}


{\bf Step 4:} We manually deal with all remaining cases of $TrueNAQt$.

{\bf `QQQQQQQSSSII'~for~$D_4$.} This case is considered in~\cite{GoodwinLeMagaard} as a family $\mathcal F_{8, 9, 10}$. 
All elements of this family stands for orbits of dimension 6. 
We denote this case $D_4(1)$. 
The picture behind this case is as follows
$$~{\Autonumfalse\mymatrix{
\lNote{2}Q \Note{1}\Rt{2pt}&\Note{2}\Bot{2pt}\pho& \Note{3}\pho\\
\lNote{3}Q&Q\Rt{2pt}&\Bot{2pt}\pho\\
\lNote{4}S&Q&Q\Rt{2pt}\\
\lNote{-4}S&Q&Q\Rt{2pt}\Bot{2pt}\\
\lNote{-3}I&S\Rt{2pt}\Bot{2pt}&\\
\lNote{-2}I\Rt{2pt}\Bot{2pt}&&\\
}}
\mapsto
{\Autonumfalse\mymatrix{
\lNote{2}A \Note{1}\Rt{2pt}&\Note{2}\Bot{2pt}\pho& \Note{3}\pho\\
\lNote{3}L&Q\Rt{2pt}&\Bot{2pt}\pho\\
\lNote{4}S&L&A\Rt{2pt}\\
\lNote{-4}S&L&A\Rt{2pt}\Bot{2pt}\\
\lNote{-3}I&S\Rt{2pt}\Bot{2pt}&\\
\lNote{-2}I\Rt{2pt}\Bot{2pt}&&\\
}}.
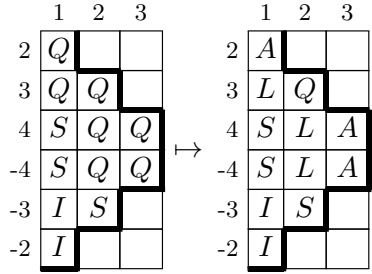
$$

{\bf `QQQQQQQQQSSSSSIIIIIIIIIIIIIIII'~for~$D_6$.} 
We denote this case $D_6(4)$. 
Applying AL-move we have: 
$$
{\Autonumfalse\mymatrix{
\lNote{2}Q \Note{1}\Rt{2pt}&
\Note{3}\Bot{2pt}\pho& \Note{3}\pho&\Note{4}\pho&\Note{5}\pho\\
\lNote{3}\yellow S&Q \Rt{2pt}&\Bot{2pt}\pho&&\\
\lNote{4}I&\yellow S&\red Q\Rt{2pt}&\Bot{2pt}&\\
\lNote{5}I&I&\red Q&\red Q\Rt{2pt}&\Bot{2pt}\\
\lNote{6}I&I&\red S&\red Q&\red Q\Rt{2pt}\\
\lNote{-6}I&I&\red S&\red Q&\red Q\Rt{2pt}\Bot{2pt}\\
\lNote{-5}I&I&I&\red S\Rt{2pt}\Bot{2pt}&\\
\lNote{-4}I&I&I\Rt{2pt}\Bot{2pt}&&\\
\lNote{-3}I&I \Note{2}\Rt{2pt}\Bot{2pt}&\pho&&\\
\lNote{-2}I\Rt{2pt}&\pho&\pho&&\\
}}
\hspace{10pt}\mapsto\hspace{10pt}
{\Autonumfalse\mymatrix{
\lNote{2}L_1 \Note{1}\Rt{2pt}&
\Note{2}\Bot{2pt}\pho& \Note{3}\pho&\Note{4}\pho&\Note{5}\pho\\
\lNote{3}\yellow S&A_1 \Rt{2pt}&\Bot{2pt}\pho&&\\
\lNote{4}I&\yellow S&\red Q\Rt{2pt}&\Bot{2pt}&\\
\lNote{5}I&I&\red Q&\red Q\Rt{2pt}&\Bot{2pt}\\
\lNote{6}I&I&\red S&\red Q&\red Q\Rt{2pt}\\
\lNote{-6}I&I&\red S&\red Q&\red Q\Rt{2pt}\Bot{2pt}\\
\lNote{-5}I&I&I&\red S\Rt{2pt}\Bot{2pt}&\\
\lNote{-4}I&I&I\Rt{2pt}\Bot{2pt}&&\\
\lNote{-3}I&I \Note{2}\Rt{2pt}\Bot{2pt}&\pho&&\\
\lNote{-2}I\Rt{2pt}&\pho&\pho&&\\
}}.
$$
We denote the new string $\mathcal S_{red}$.
We have that $\nt_{S_{red}}$ is a direct sum of the two-dimensional abelian Lie algebra standing for ``yellow'' letters `S' and the 10-dimensional Lie algebra $\nt^R$ standing for ``red'' letters. 
It is easy to see that $\nt^R$ is isomorphic to $\nt_{S}$ from the above case $D_4(1)$.
This implies that all orbits of case $D_6(4)$ are of dimension 8 ($>4$).
\sect{Orbits of dimension 6}\label{S:dim6}\addcontentsline{toc}{subsection}{\ref{S:dim6}. Orbits of dimension 6}
In this section we would like to enlist all extensive orbits of dimension 6. The technique is very similar to the technique developed for dimensions 2, 4. 

Thanks to Proposition~\ref{PlWd} the only algebras for which there is a chance to find extensive orbits of dimension 6 are $A_n, 2\le n\le 7,  B_n, 2\le n\le 7, C_n, 2\le n\le 7, D_n, 4\le n\le 8, G_2, F_4, E_6, E_7, E_8$ (most of them indeed allow such orbits with the exceptions of $A_2, B_2=C_2$ --- they are too small). 
The Dynkin diagrams of type $G_2, F_4, E_6, E_7, E_8$ can not occur as subdiagrams of the Dynkin diagrams of the classical series and thus the classification of orbits of dimension 6 for the classical series does not involve these exceptional cases.  
For each such a case we provide a list of subsets $S_1, S_2,\ldots$ of $\Phi^+$ such that each extensive coadjoint orbit of dimension 6 in $\nt^*$ intersects $V(S_1)\sqcup V(S_2)\sqcup\ldots$ in exactly one point. 
\mtheo{\label{Tdim6} The strings of \textup19 Tables of Section~\textup{\ref{S:clstr6}} enlists \textup(in the above sense\textup) all extensive orbits of dimension 6 for Lie algebras of types $A_n, 2\le n\le 7,  B_n, 2\le n\le 7, C_n, 2\le n\le 7, D_n$\textup, $4\le n\le 8$.}
The respective Tables are very long alltogether and we place them in a separate section, Section~\ref{S:clstr6}. 
\sst{The algorithm for dimension 6}
We will now show how one can get the desired strings of Theorem~\ref{Tdim6}. 
We extensively use the notation and the concepts of Sections~\ref{S:dim2} and~\ref{S:dim4}. 
In particular, we work with strings $\mathcal S$ and sets $S, Z$ for all the needed root systems (they are enlisted in the statement of Theorem~\ref{Tdim6}). 

We first do the analogue of Step 1 of the proof of Theorem~\ref{Tdim2} for the respective root system with the following modifications: 1)
if $6>|\Sing(\delta_0)|$ in the above notation then $MaskD$ set to be the string of `Q' letters only; 2) we replace $\SuppD_2$ by $\SuppD_6$.

Next we do the straightforward analogue of Steps 2 and 3.1 of the proof of Theorem~\ref{Tdim2} for the respective root system with the following modification: we add a string from $ClsPat$ to $ClsPatD$ if it contains exactly 3 letters `A'.

After that we do the straighforward analogue of Step 3.2 of the proof of Theorem~\ref{Tdim4} with the modification as above. We next do Step 3.3 of the proof of Theorem~\ref{Tdim2} removing only strings which stand for sets where all orbits are of dimension 8 or more. 

For $A_n, B_n$ and $C_n$ cases we have that $TrueNAQt$ is empty after Step 3.3 and thus everything works as in the main cases of dimension 4. 
We left to deal with the remaining $D_n$-cases. The number of them is as follows:
$$\begin{array}{|c|c|c|c|c|c|c|c|c|c|c|c|c|c|c|}\hline
Dim=6&D_4&D_5&D_6&D_7&D_8\\
\hline {\rm \# Special~cases} &1&2&4&3&4\\
\hline
\end{array}.$$
The only case in $D_4$ stands for the family $\mathcal F_{8, 9, 10}$ of 6-dimensional orbits described in~\cite[Table~4]{GoodwinLeMagaard}. 
We recall that we denote this case $D_4(1)$. 
A set-section for the respective action is given by $V({\{\alpha_8, \alpha_9, \alpha_{10}\}})\sqcup V({\{\alpha_3, \alpha_8, \alpha_9, \alpha_{10}\}}),$ see~\cite[Table~4]{GoodwinLeMagaard}. 

It turns out that for $D_5, D_6, D_7, D_8$ all the remaining cases provide no orbits of dimension 6 and we check it case by case. 

{\bf The two cases of $D_5$.} They stand for the strings 
\begin{center}$M_1:=$`QQQQQQQQQQQQQSSSIIII' and $M_2:=$
`QQQQQQQQQSQQQISSIIII'.\end{center} 
We denote these cases $D_5(1)$ and $D_5(2)$ respectively. 
We will deal with $D_5(1)$ first. 
We would like to show that this case provides no orbits of dimension 6. 
Applying two AL-moves we obtain
$$M_1\to
~{\Autonumfalse\mymatrix{
\lNote{2}\red Q \Note{1}\Rt{2pt}&
\Note{2}\Bot{2pt}\pho& \Note{3}\pho&\Note{4}\pho\\
\lNote{3}A_2&A_1 \Rt{2pt}&\Bot{2pt}\pho&\\
\lNote{4}\red Q&\red Q&\yellow Q\Rt{2pt}&\Bot{2pt}\\
\lNote{5}\red S&\red Q&\yellow Q&\red Q\Rt{2pt}\\
\lNote{-5}\red S&\red Q&L_2&\red Q\Rt{2pt}\Bot{2pt}\\
\lNote{-4}I&\red S&L_1\Rt{2pt}\Bot{2pt}&\\
\lNote{-3}I&I \Note{2}\Rt{2pt}\Bot{2pt}&\pho&\\
\lNote{-2}I\Rt{2pt}&\pho&\pho&\\
}}.
$$
It is easy to see that the subalgebra $\nt^R$ defined by the red cells is isomorphic to the subalgebra $\nt_S$ from the subcase $D_4(1)$. 
If $N.f\in\nt_{S(M_1)}^*$ is $Z$-saturated then the restriction of $N.f$ to $\nt^R$ will consist of several orbits from the family $\mathcal F_{8, 9, 10}$ of~\cite{GoodwinLeMagaard}. 
This implies that $\dim N.f\ge 6+4$. 
Therefore for $\nt_{S(M_1)}^*$ all $Z$-saturated orbits are of dimensions 6+4 ($>6$) or more. 
This completes case $D_5(1)$.

In all subsequent cases we denote by $\nt^R$ the subalgebra defined by the red cells on the respective pictures and $N^R$ the respective subgroup of $N$. 
Moreover, for most cases we use a similar strategy: we apply first $d$ (with $d=0, 1,$ or $2$) AL-moves and then project the respective orbits on the subalgebra $\nt^R$ defined by the red cells of the respective picture. For this subalgebra we already knew from a previous case that the $Z$-saturated orbits are of dimension $X$ or more and then we deduce from it that the dimensions of the needed orbits are $X+2d$ or more. 
In all the cases we will have that $X+2d>6$ and this completes the argument.

$$\begin{array}{l}\text{Case~}D_5(2)\\\nt^R\to D_4(1)\\X=6, d=1\end{array}\to
~{\Autonumfalse\mymatrix{
\lNote{2}A_1 \Note{1}\Rt{2pt}&
\Note{2}\Bot{2pt}\pho& \Note{3}\pho&\Note{4}\pho\\
\lNote{3}\red Q&\yellow Q \Rt{2pt}&\Bot{2pt}\pho&\\
\lNote{4}\red Q& \yellow Q&\red Q\Rt{2pt}&\Bot{2pt}\\
\lNote{5}\red S&\yellow Q&\red Q&\red Q\Rt{2pt}\\
\lNote{-5}\red S& L_1&\red Q&\red Q\Rt{2pt}\Bot{2pt}\\
\lNote{-4}I&I&\red S\Rt{2pt}\Bot{2pt}&\\
\lNote{-3}I&I \Note{2}\Rt{2pt}\Bot{2pt}&\pho&\\
\lNote{-2}I\Rt{2pt}&\pho&\pho&\\
}}.
$$

{\bf The four cases of $D_6$.} They stand for strings $M_1, M_2, M_3, M_4:$
$$\begin{tabular}{ll} `QQQQQQQQQQQQQQQSSSSIIIIIIIIIII',&`QQQQQQQQQQSSQQQIISSIIIIIIIIIII',\\
`QQQQQQQQQQSSSSSIIIIIIIIIIIIIII',&`QQQQQQQQQSSSSSIIIIIIIIIIIIIIII'\end{tabular}.$$
We denote these cases by $D_6(1), D_6(2), D_6(3), D_6(4)$  respectively. 
We will deal with $D_6(1)$ first. We would like to show that this case provides no orbits of dimension 6. Assume to the contrary that there exists a $Z$-saturated orbit $N_{S(M_1)}.f$. Then its restriction to the subalgebra $\nt^R$ is a union of $Z$-saturated orbits of $\nt^R$.
$$M_1\to~{\Autonumfalse\mymatrix{
\lNote{2}\red Q \Note{1}\Rt{2pt}&
\Note{2}\Bot{2pt}\pho& \Note{3}\pho&\Note{4}\pho&\Note{5}\pho\\
\lNote{3}\red Q&\red Q \Rt{2pt}&\Bot{2pt}\pho&&\\
\lNote{4}\red S&\red Q&\red Q\Rt{2pt}&\Bot{2pt}&\\
\lNote{5}I&\yellow Q&\yellow Q&\yellow Q\Rt{2pt}&\Bot{2pt}\\
\lNote{6}I&\red S&\red Q&\yellow Q&\red Q\Rt{2pt}\\
\lNote{-6}I&\red S&\red Q&\yellow Q&\red Q\Rt{2pt}\Bot{2pt}\\
\lNote{-5}I&I&\red S&\yellow Q\Rt{2pt}\Bot{2pt}&\\
\lNote{-4}I&I&I\Rt{2pt}\Bot{2pt}&&\\
\lNote{-3}I&I \Note{2}\Rt{2pt}\Bot{2pt}&\pho&&\\
\lNote{-2}I\Rt{2pt}&\pho&\pho&&\\
}}.
$$
We can apply 4 AL-moves to the $\nt^R$:
$$
~{\Autonumfalse\mymatrix{
\lNote{2}\red Q \Note{1}\Rt{2pt}&
\Note{2}\Bot{2pt}\pho& \Note{3}\pho&\Note{4}\pho&\Note{5}\pho\\
\lNote{3}\red Q&\red Q \Rt{2pt}&\Bot{2pt}\pho&&\\
\lNote{4}\red S&\red Q&\red Q\Rt{2pt}&\Bot{2pt}&\\
\lNote{5}I&\yellow Q&\yellow Q&\yellow Q\Rt{2pt}&\Bot{2pt}\\
\lNote{6}I&\red S&\red Q&\yellow Q&\red Q\Rt{2pt}\\
\lNote{-6}I&\red S&\red Q&\yellow Q&\red Q\Rt{2pt}\Bot{2pt}\\
\lNote{-5}I&I&\red S&\yellow Q\Rt{2pt}\Bot{2pt}&\\
\lNote{-4}I&I&I\Rt{2pt}\Bot{2pt}&&\\
\lNote{-3}I&I \Note{2}\Rt{2pt}\Bot{2pt}&\pho&&\\
\lNote{-2}I\Rt{2pt}&\pho&\pho&&\\
}}
\to
~{\Autonumfalse\mymatrix{
\lNote{2}\red L_2 \Note{1}\Rt{2pt}&
\Note{2}\Bot{2pt}\pho& \Note{3}\pho&\Note{4}\pho&\Note{5}\pho\\
\lNote{3}\red L_1&\red Q \Rt{2pt}&\Bot{2pt}\pho&&\\
\lNote{4}\red S&\red A_2&\red A_1\Rt{2pt}&\Bot{2pt}&\\
\lNote{5}I&\yellow I&\yellow I&\yellow I\Rt{2pt}&\Bot{2pt}\\
\lNote{6}I&\red S&\red A_3&\yellow I&\red L_3\Rt{2pt}\\
\lNote{-6}I&\red S&\red A_4&\yellow I&\red L_4\Rt{2pt}\Bot{2pt}\\
\lNote{-5}I&I&\red S&\yellow I\Rt{2pt}\Bot{2pt}&\\
\lNote{-4}I&I&I\Rt{2pt}\Bot{2pt}&&\\
\lNote{-3}I&I \Note{2}\Rt{2pt}\Bot{2pt}&\pho&&\\
\lNote{-2}I\Rt{2pt}&\pho&\pho&&\\
}}.
$$
Hence all $Z$-saturated orbits of $\nt^R$ are of dimension 8 or more. 
This completes case $D_6(1)$.

The arguments for $D_6(2), D_6(3)$ are more straightforward: 
$$\begin{array}{l}\text{Case~}D_6(2)\\ \nt^R\to D_5(2)\\X=8, d=0\end{array}\to~{\Autonumfalse\mymatrix{
\lNote{2}\yellow Q \Note{1}\Rt{2pt}&
\Note{2}\Bot{2pt}\pho& \Note{3}\pho&\Note{4}\pho&\Note{5}\pho\\
\lNote{3}\yellow S&\red Q \Rt{2pt}&\Bot{2pt}\pho&&\\
\lNote{4}I&\red Q&\red Q\Rt{2pt}&\Bot{2pt}&\\
\lNote{5}I&\red Q&\red Q&\red Q\Rt{2pt}&\Bot{2pt}\\
\lNote{6}I&\red S&\red Q&\red Q&\red Q\Rt{2pt}\\
\lNote{-6}I&\red S&\red Q&\red Q&\red Q\Rt{2pt}\Bot{2pt}\\
\lNote{-5}I&I&I&\red S\Rt{2pt}\Bot{2pt}&\\
\lNote{-4}I&I&I\Rt{2pt}\Bot{2pt}&&\\
\lNote{-3}I&I \Note{2}\Rt{2pt}\Bot{2pt}&\pho&&\\
\lNote{-2}I\Rt{2pt}&\pho&\pho&&\\
}},\hspace{10pt}\begin{array}{l}\text{Case~}D_6(3)\\\nt^R\to D_4(1)\\X=6, d=1\end{array}\to~{\Autonumfalse\mymatrix{
\lNote{2}L_1 \Note{1}\Rt{2pt}&
\Note{2}\Bot{2pt}\pho& \Note{3}\pho&\Note{4}\pho&\Note{5}\pho\\
\lNote{3}\yellow S&A_1 \Rt{2pt}&\Bot{2pt}\pho&&\\
\lNote{4}I&\yellow Q&\red Q\Rt{2pt}&\Bot{2pt}&\\
\lNote{5}I&\yellow S&\red Q&\red Q\Rt{2pt}&\Bot{2pt}\\
\lNote{6}I&I&\red S&\red Q&\red Q\Rt{2pt}\\
\lNote{-6}I&I&\red S&\red Q&\red Q\Rt{2pt}\Bot{2pt}\\
\lNote{-5}I&I&I&\red S\Rt{2pt}\Bot{2pt}&\\
\lNote{-4}I&I&I\Rt{2pt}\Bot{2pt}&&\\
\lNote{-3}I&I \Note{2}\Rt{2pt}\Bot{2pt}&\pho&&\\
\lNote{-2}I\Rt{2pt}&\pho&\pho&&\\
}}.
$$
The case $D_6(4)$ was considered in subsection~\ref{SS:adim4} where we have shown that the respective orbits are of dimension 8 ($>6$). 

{\bf The three cases of $D_7$.} They stand for the strings $M_1, M_2, M_3$:
$$\mbox{`QQQQQQQQQQQQQSQQQISISSIIIIIIIIIIIIIIIIIIII'},$$
$$\mbox{`QQQQQQQQQQQSSSQQQIIISSIIIIIIIIIIIIIIIIIIII'},$$
$$\mbox{`QQQQQQQQQQSQQSSSIISIIIIIIIIIIIIIIIIIIIIIII'}.$$
We denote the respective cases by $D_7(1), D_7(2), D_7(3)$. 
$$\begin{array}{l}\text{Case~}D_7(1)\\\nt^R\to D_5(2)\\X=8, d=0\end{array}\to~{\Autonumfalse\mymatrix{
\lNote{2}\yellow Q \Note{1}\Rt{2pt}&
\Note{2}\Bot{2pt}\pho& \Note{3}\pho&\Note{4}\pho&\Note{5}\pho&\Note{6}\pho\\
\lNote{3}\yellow Q&\yellow Q \Rt{2pt}&\Bot{2pt}\pho&&&\\
\lNote{4}\yellow S&\yellow Q& \red Q\Rt{2pt}&\Bot{2pt}&&\\
\lNote{5}I&I&\red Q&\red Q\Rt{2pt}&\Bot{2pt}&\\
\lNote{6}I&I&\red Q&\red Q&\red Q\Rt{2pt}&\Bot{2pt}\\
\lNote{7}I&I&\red S&\red Q&\red Q&\red Q\Rt{2pt}\\
\lNote{-7}I&I&\red S&\red Q&\red Q&\red Q\Rt{2pt}\Bot{2pt}\\
\lNote{-6}I&I&I&I&\red S\Rt{2pt}\Bot{2pt}&\\
\lNote{-5}I&I&I&I\Rt{2pt}\Bot{2pt}&&\\
\lNote{-4}I&I&I\Rt{2pt}\Bot{2pt}&&&\\
\lNote{-2}I&I \Note{2}\Rt{2pt}\Bot{2pt}&\pho&&&\\
\lNote{-2}I\Rt{2pt}&\pho&\pho&&&\\
}},\hspace{10pt}
\begin{array}{l}\text{Case~}D_7(2)\\\nt^R\to D_5(2)\\X=8, d=0\end{array}\to~{\Autonumfalse\mymatrix{
\lNote{2}\yellow Q \Note{1}\Rt{2pt}&
\Note{2}\Bot{2pt}\pho& \Note{3}\pho&\Note{4}\pho&\Note{5}\pho&\Note{6}\pho\\
\lNote{3}\yellow S&\yellow Q \Rt{2pt}&\Bot{2pt}\pho&&&\\
\lNote{4}I&\yellow S& \red Q\Rt{2pt}&\Bot{2pt}&&\\
\lNote{5}I&I&\red Q&\red Q\Rt{2pt}&\Bot{2pt}&\\
\lNote{6}I&I&\red Q&\red Q&\red Q\Rt{2pt}&\Bot{2pt}\\
\lNote{7}I&I&\red S&\red Q&\red Q&\red Q\Rt{2pt}\\
\lNote{-7}I&I&\red S&\red Q&\red Q&\red Q\Rt{2pt}\Bot{2pt}\\
\lNote{-6}I&I&I&I&\red S\Rt{2pt}\Bot{2pt}&\\
\lNote{-5}I&I&I&I\Rt{2pt}\Bot{2pt}&&\\
\lNote{-4}I&I&I\Rt{2pt}\Bot{2pt}&&&\\
\lNote{-2}I&I \Note{2}\Rt{2pt}\Bot{2pt}&\pho&&&\\
\lNote{-2}I\Rt{2pt}&\pho&\pho&&&\\
}}.
$$
$$\begin{array}{l}\text{Case~}D_7(3)\\\nt^R\to D_4(1)\\X=6, d=2\end{array}\to~{\Autonumfalse\mymatrix{
\lNote{2}L_2 \Note{1}\Rt{2pt}&
\Note{2}\Bot{2pt}\pho& \Note{3}\pho&\Note{4}\pho&\Note{5}\pho&\Note{6}\pho\\
\lNote{3}L_1&\yellow Q \Rt{2pt}&\Bot{2pt}\pho&&&\\
\lNote{4}\yellow S&A_2& A_1\Rt{2pt}&\Bot{2pt}&&\\
\lNote{5}I&I&\yellow S&\red Q\Rt{2pt}&\Bot{2pt}&\\
\lNote{6}I&I&I&\red Q&\red Q\Rt{2pt}&\Bot{2pt}\\
\lNote{7}I&I&I&\red S&\red Q&\red Q\Rt{2pt}\\
\lNote{-7}I&I&I&\red S&\red Q&\red Q\Rt{2pt}\Bot{2pt}\\
\lNote{-6}I&I&I&I&\red S\Rt{2pt}\Bot{2pt}&\\
\lNote{-5}I&I&I&I\Rt{2pt}\Bot{2pt}&&\\
\lNote{-4}I&I&I\Rt{2pt}\Bot{2pt}&&&\\
\lNote{-3}I&I \Note{2}\Rt{2pt}\Bot{2pt}&\pho&&&\\
\lNote{-2}I\Rt{2pt}&\pho&\pho&&&\\
}}.
$$
{\bf The four cases of $D_8$.} They stand for the strings $M_1, M_2, M_3, M_4:$
$$`QQQQQQQQQQQQQQSSQQQISIISSIIIIIIIIIIIIIIIIIIIIIIIIIIIIIII',$$
$$`QQQQQQQQQQQQSSSSQQQIIIISSIIIIIIIIIIIIIIIIIIIIIIIIIIIIIII',$$
$$`QQQQQQQQQQQSQQSSSSIISIIIIIIIIIIIIIIIIIIIIIIIIIIIIIIIIIII',$$
$$`QQQQQQQQQQQSSSSSSSIIIIIIIIIIIIIIIIIIIIIIIIIIIIIIIIIIIIII'.$$
We denote the respective cases by $D_8(1), D_8(2), D_8(3), D_8(4)$. 
$$\begin{array}{l}\text{Case~}D_8(1)\\\nt^R\to D_5(2)\\X=8, d=0\end{array}\to ~{\Autonumfalse\mymatrix{
\lNote{2}\yellow Q \Note{1}\Rt{2pt}&
\Note{2}\Bot{2pt}\pho& \Note{3}\pho&\Note{4}\pho&\Note{5}\pho&\Note{6}\pho&\Note{7}\pho\\
\lNote{3}\yellow S& \yellow Q \Rt{2pt}&\Bot{2pt}\pho&&&&\\
\lNote{4}I&\yellow Q&\yellow Q\Rt{2pt}&\Bot{2pt}&&&\\
\lNote{5}I&\yellow S&\yellow Q&\red Q\Rt{2pt}&\Bot{2pt}&&\\
\lNote{6}I&I&I&\red Q&\red Q\Rt{2pt}&\Bot{2pt}&\\
\lNote{7}I&I&I&\red Q&\red Q&\red Q\Rt{2pt}&\Bot{2pt}\\
\lNote{8}I&I&I&\red S&\red Q&\red Q&\red Q\Rt{2pt}\\
\lNote{-8}I&I&I&\red S&\red Q&\red Q&\red Q\Rt{2pt}\Bot{2pt}\\
\lNote{-7}I&I&I&I&I&\red S\Rt{2pt}\Bot{2pt}&\\
\lNote{-6}I&I&I&I&I\Rt{2pt}\Bot{2pt}&&\\
\lNote{-5}I&I&I&I\Rt{2pt}\Bot{2pt}&&&\\
\lNote{-4}I&I&I\Rt{2pt}\Bot{2pt}&&&&\\
\lNote{-3}I&I \Note{2}\Rt{2pt}\Bot{2pt}&\pho&&&&\\
\lNote{-2}I\Rt{2pt}&\pho&\pho&&&&\\
}},
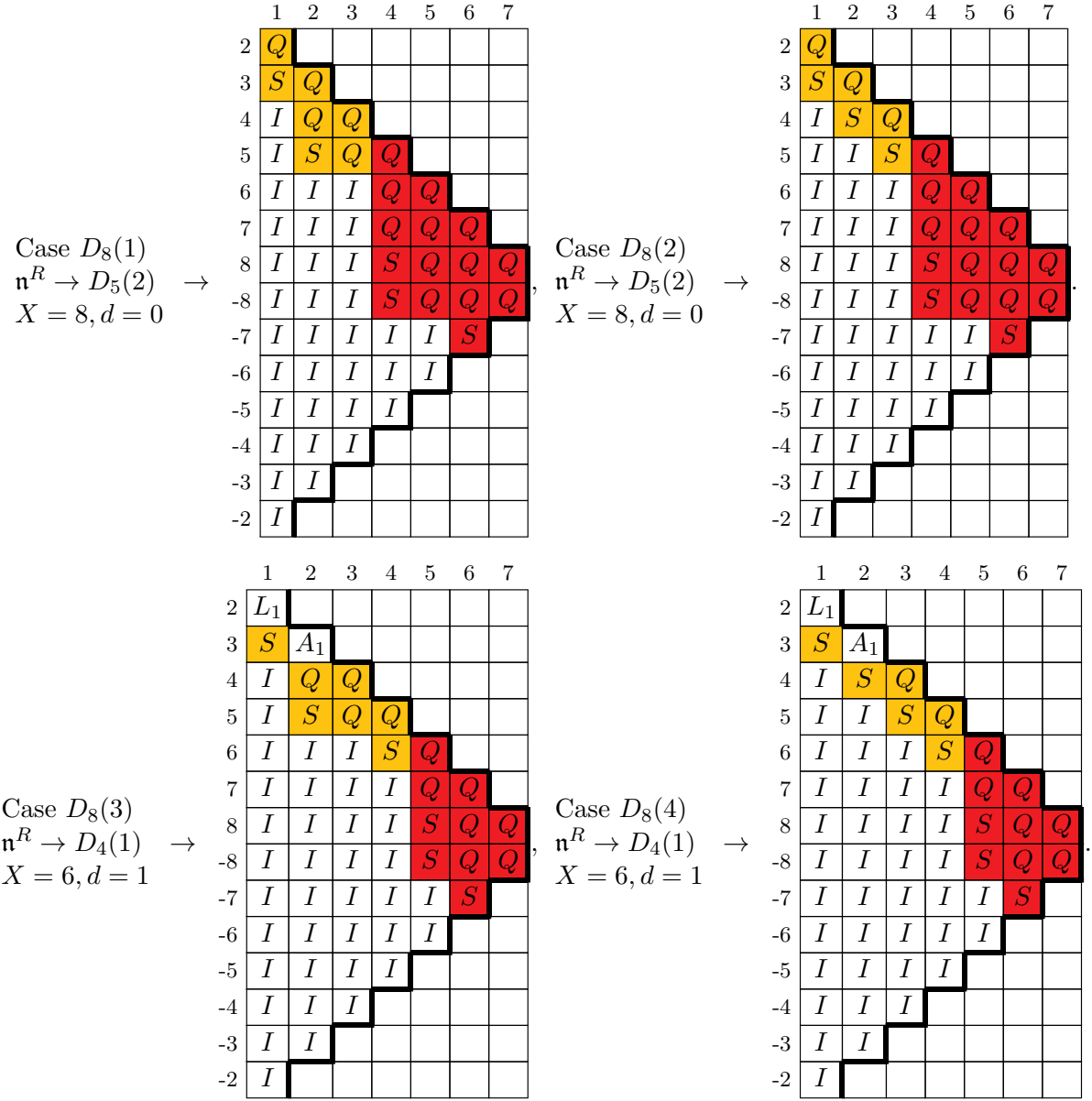
\begin{array}{l}\text{Case~}D_8(2)\\\nt^R\to D_5(2)\\X=8, d=0\end{array}\to ~{\Autonumfalse\mymatrix{
\lNote{2}\yellow Q \Note{1}\Rt{2pt}&
\Note{2}\Bot{2pt}\pho& \Note{3}\pho&\Note{4}\pho&\Note{5}\pho&\Note{6}\pho&\Note{7}\pho\\
\lNote{3}\yellow S& \yellow Q \Rt{2pt}&\Bot{2pt}\pho&&&&\\
\lNote{4}I&\yellow S&\yellow Q\Rt{2pt}&\Bot{2pt}&&&\\
\lNote{5}I&I&\yellow S&\red Q\Rt{2pt}&\Bot{2pt}&&\\
\lNote{6}I&I&I&\red Q&\red Q\Rt{2pt}&\Bot{2pt}&\\
\lNote{7}I&I&I&\red Q&\red Q&\red Q\Rt{2pt}&\Bot{2pt}\\
\lNote{8}I&I&I&\red S&\red Q&\red Q&\red Q\Rt{2pt}\\
\lNote{-8}I&I&I&\red S&\red Q&\red Q&\red Q\Rt{2pt}\Bot{2pt}\\
\lNote{-7}I&I&I&I&I&\red S\Rt{2pt}\Bot{2pt}&\\
\lNote{-6}I&I&I&I&I\Rt{2pt}\Bot{2pt}&&\\
\lNote{-5}I&I&I&I\Rt{2pt}\Bot{2pt}&&&\\
\lNote{-4}I&I&I\Rt{2pt}\Bot{2pt}&&&&\\
\lNote{-3}I&I \Note{2}\Rt{2pt}\Bot{2pt}&\pho&&&&\\
\lNote{-2}I\Rt{2pt}&\pho&\pho&&&&\\
}}.$$
$$\begin{array}{l}\text{Case~}D_8(3)\\\nt^R\to D_4(1)\\X=6, d=1\end{array}\to ~{\Autonumfalse\mymatrix{
\lNote{2}L_1 \Note{1}\Rt{2pt}&
\Note{2}\Bot{2pt}\pho& \Note{3}\pho&\Note{4}\pho&\Note{5}\pho&\Note{6}\pho&\Note{7}\pho\\
\lNote{3}\yellow S& A_1 \Rt{2pt}&\Bot{2pt}\pho&&&&\\
\lNote{4}I&\yellow Q&\yellow Q\Rt{2pt}&\Bot{2pt}&&&\\
\lNote{5}I&\yellow S&\yellow Q&\yellow Q\Rt{2pt}&\Bot{2pt}&&\\
\lNote{6}I&I&I&\yellow S&\red Q\Rt{2pt}&\Bot{2pt}&\\
\lNote{7}I&I&I&I&\red Q&\red Q\Rt{2pt}&\Bot{2pt}\\
\lNote{8}I&I&I&I&\red S&\red Q&\red Q\Rt{2pt}\\
\lNote{-8}I&I&I&I&\red S&\red Q&\red Q\Rt{2pt}\Bot{2pt}\\
\lNote{-7}I&I&I&I&I&\red S\Rt{2pt}\Bot{2pt}&\\
\lNote{-6}I&I&I&I&I\Rt{2pt}\Bot{2pt}&&\\
\lNote{-5}I&I&I&I\Rt{2pt}\Bot{2pt}&&&\\
\lNote{-4}I&I&I\Rt{2pt}\Bot{2pt}&&&&\\
\lNote{-3}I&I \Note{2}\Rt{2pt}\Bot{2pt}&\pho&&&&\\
\lNote{-2}I\Rt{2pt}&\pho&\pho&&&&\\
}},
\begin{array}{l}\text{Case~}D_8(4)\\\nt^R\to D_4(1)\\X=6, d=1\end{array}\to ~{\Autonumfalse\mymatrix{
\lNote{2}L_1 \Note{1}\Rt{2pt}&
\Note{2}\Bot{2pt}\pho& \Note{3}\pho&\Note{4}\pho&\Note{5}\pho&\Note{6}\pho&\Note{7}\pho\\
\lNote{3}\yellow S& A_1 \Rt{2pt}&\Bot{2pt}\pho&&&&\\
\lNote{4}I&\yellow S&\yellow Q\Rt{2pt}&\Bot{2pt}&&&\\
\lNote{5}I&I&\yellow S&\yellow Q\Rt{2pt}&\Bot{2pt}&&\\
\lNote{6}I&I&I&\yellow S&\red Q\Rt{2pt}&\Bot{2pt}&\\
\lNote{7}I&I&I&I&\red Q&\red Q\Rt{2pt}&\Bot{2pt}\\
\lNote{8}I&I&I&I&\red S&\red Q&\red Q\Rt{2pt}\\
\lNote{-8}I&I&I&I&\red S&\red Q&\red Q\Rt{2pt}\Bot{2pt}\\
\lNote{-7}I&I&I&I&I&\red S\Rt{2pt}\Bot{2pt}&\\
\lNote{-6}I&I&I&I&I\Rt{2pt}\Bot{2pt}&&\\
\lNote{-5}I&I&I&I\Rt{2pt}\Bot{2pt}&&&\\
\lNote{-4}I&I&I\Rt{2pt}\Bot{2pt}&&&&\\
\lNote{-3}I&I \Note{2}\Rt{2pt}\Bot{2pt}&\pho&&&&\\
\lNote{-2}I\Rt{2pt}&\pho&\pho&&&&\\
}}.$$
\sect{The classifying strings for orbits of dimension 6}\label{S:clstr6}\addcontentsline{toc}{subsection}{\ref{S:clstr6}. The classifying strings for orbits of dimension 6}
For some tables the respective strings are too short for the respective root system. In this case the convention is that we fill the end of the respective string with letters `I' (this shorten the answer a lot: for example for the case of $D_8$ this shorten the answer by 35 letters per input to 21 letters per input).
\begin{longtable}{llllll}
\multicolumn{6}{c}{Table~\label{T:dim6P1}\ref{T:dim6P1}: $A_4$ case}\\
\hline
ASSAAILLLS&
ASIAAILLLS&
AISAAILLLS\\
AIIAAILLLS&
LAASLALSSI&
LAAILALSSI\\
\hline\end{longtable}

\begin{longtable}{llllll}
\multicolumn{6}{c}{Table~\label{T:dim6P2}\ref{T:dim6P2}: $A_5$ case}\\
\hline
ASLALALSSSIISII&
ASLALALISSIISII&
AILALALSSSIISII&
AILALALISSIISII\\
ASSALAILSILISII&
AISALAILSILISII&
ASIALAILSILISII&
AIIALAILSILISII\\
LLASSLAAISSSIII&
LLASILAAISSSIII&
LLAISLAAISSSIII&
LLAIILAAISSSIII\\
ASALALLSSSISIII&
ASALALLISSISIII&
ASASALLLISISIII&
ASAIALLLISISIII\\
AIALALLSSSISIII&
AIALALLISSISIII&
AIASALLLISISIII&
AIAIALLLISISIII\\
ASAALLLSSSIIIII&
AIAALLLSSSIIIII&
LALSASSLAIISISI&
LALSASILAIISISI\\
LALIASSLAIISISI&
LALIASILAIISISI&
LASSASAILILIISI&
LASIASAILILIISI\\
LAISASAILILIISI&
LAIIASAILILIISI&
LAASASSLLIISIII&
LAAIASSLLIISIII\\\hline\end{longtable}

\begin{longtable}{llllll}
\multicolumn{6}{c}{Table~\label{T:dim6P3}\ref{T:dim6P3}: $A_6$ case}\\
\hline
ALALALSSSSSSISIII&
ALALALSISSSSISIII&
ALALALSSSISSISIII&
ALALALSISISSISIII\\
ASALALILSSSSISIII&
ASALALILSISSISIII&
AIALALILSSSSISIII&
AIALALILSISSISIII\\
ASASALLLISSSISIII&
ASASALLLIISSISIII&
AIASALLLISSSISIII&
AIASALLLIISSISIII\\
ASAIALLLISSSISIII&
ASAIALLLIISSISIII&
AIAIALLLISSSISIII&
AIAIALLLIISSISIII\\
ALALALSSSSSSIIIII&
ALALALSISSSSIIIII&
ASALALILSSSSIIIII&
AIALALILSSSSIIIII\\
LALALASSSSSISISII&
LALALASSSISISISII&
LALASASSSLIISISII&
LALAIASSSLIISISII\\
LALALASISSSISISII&
LALALASISISISISII&
LALASASISLIISISII&
LALAIASISLIISISII\\
LASALASLISSISISII&
LASALASLIISISISII&
LASASASLILIISISII&
LASAIASLILIISISII\\
LAIALASLISSISISII&
LAIALASLIISISISII&
LAIASASLILIISISII&
LAIAIASLILIISISII\\
LALAALSSSSSISIIII&
LALAALSISSSISIIII&
LASAALSLISSISIIII&
LAIAALSLISSISIIII\\
LAALALSSSSSIISIII&
LAALALSSSISIISIII&
LAASALSSILSIISIII&
LAAIALSSILSIISIII\\
LALALASSSSSIIISII&
LALALASSSISIIISII&
LALASASSSLIIIISII&
LALAIASSSLIIIISII\\
LALAALSSSSSIIIIII&
LALSASSSLASIISIIS&
LALSASSILASIISIIS&
LALIASSSLASIISIIS\\
LALIASSILASIISIIS&
LALSAISSLASIISIIS&
LALSAISILASIISIIS&
LALIAISSLASIISIIS\\
LALIAISILASIISIIS&
LASLALSAISSISIIIS&
LASLALSAIISISIIIS&
LASSALSAILSIIIIIS\\
LASIALSAILSIIIIIS&
LAILALSAISSISIIIS&
LAILALSAIISISIIIS&
LAISALSAILSIIIIIS\\
LAIIALSAILSIIIIIS&\\
\hline
\end{longtable}

\begin{longtable}{llllll}
\multicolumn{6}{c}{Table~\label{T:dim6P4}\ref{T:dim6P4}: $A_7$ case}\\
\hline
LALALASSSSSSSISIS&
LALALASSSSISSISIS&
LALALAISSSSSSISIS\\
LALALAISSSISSISIS&
LALASALSSSSISISIS&
LALASALSSSIISISIS\\
LALAIALSSSSISISIS&
LALAIALSSSIISISIS&
LALALASSISSSSISIS\\
LALALASSISISSISIS&
LALALAISISSSSISIS&
LALALAISISISSISIS\\
LALASALSISSISISIS&
LALASALSISIISISIS&
LALAIALSISSISISIS\\
LALAIALSISIISISIS&
LASALALSSISSSISIS&
LASALALSSISISISIS\\
LASALALSIISSSISIS&
LASALALSIISISISIS&
LASASALSLIISSISIS\\
LASASALSLIIISISIS&
LASAIALSLIISSISIS&
LASAIALSLIIISISIS\\
LAIALALSSISSSISIS&
LAIALALSSISISISIS&
LAIALALSIISSSISIS\\
LAIALALSIISISISIS&
LAIASALSLIISSISIS&
LAIASALSLIIISISIS\\
LAIAIALSLIISSISIS&
LAIAIALSLIIISISIS&
LALALASSSSSSSISII\\
LALALAISSSSSSISII&
LALALASSISSSSISII&
LALALAISISSSSISII\\
LASALALSSISSSISII&
LASALALSIISSSISII&
LAIALALSSISSSISII\\
LAIALALSIISSSISII&
LALALASSSSSSSIIIS&
LALALASSSSISSIIIS\\
LALALAISSSSSSIIIS&
LALALAISSSISSIIIS&
LALASALSSSSISIIIS\\
LALASALSSSIISIIIS&
LALAIALSSSSISIIIS&
LALAIALSSSIISIIIS\\
LALALASSSSSSSIIII&
LALALAISSSSSSIIII&\\
\hline\end{longtable}

\begin{longtable}{llllll}
\multicolumn{6}{c}{Table~\label{T:dim6P6}\ref{T:dim6P6}: $B_3$ case}\\
\hline
SASAALLLS&
SAIAALLLS&
IASAALLLS&
IAIAALLLS\\
\hline\end{longtable}

\begin{longtable}{llllll}
\multicolumn{6}{c}{Table~\label{T:dim6P7}\ref{T:dim6P7}: $B_4$ case}\\
\hline
LASLAALSSSIISI&
LASLAALSISIISI&
LAILAALSSSIISI&
LAILAALSISIISI\\
SASLAALILSIISI&
SAILAALILSIISI&
IASLAALILSIISI&
IAILAALILSIISI\\
ALSASLASISILIS&
ALSASLAIISILIS&
ALIASLASISILIS&
ALIASLAIISILIS\\
ASSASIALIILLIS&
ASIASIALIILLIS&
AISASIALIILLIS&
AIIASIALIILLIS\\
ALSAILASISILIS&
ALSAILAIISILIS&
ALIAILASISILIS&
ALIAILAIISILIS\\
ASSAIIALIILLIS&
ASIAIIALIILLIS&
AISAIIALIILLIS&
AIIAIIALIILLIS\\
ASSAAILSLLISII&
ASIAAILSLLISII&
AISAAILSLLISII&
AIIAAILSLLISII\\
LALSAALSSSIIII&
LALIAALSSSIIII&
AASALLLSISIIII&
AAIALLLSISIIII\\
ASSAAILILLISII&
ASIAAILILLISII&
AISAAILILLISII&
AIIAAILILLISII\\
LAASLALISSIIII&
LAAILALISSIIII&&\\\hline\end{longtable}

\begin{longtable}{llllll}
\multicolumn{6}{c}{Table~\label{T:dim6P8}\ref{T:dim6P8}: $B_5$ case}\\
\hline
LALSASSLAIIISIISII&
LALSASILAIIISIISII&
LALIASSLAIIISIISII\\
LALIASILAIIISIISII&
LASSASAILIILIIISII&
LASIASAILIILIIISII\\
LAISASAILIILIIISII&
LAIIASAILIILIIISII&
ASASALSLLSSISSIIII\\
ASASALSLLSIISSIIII&
ASASALILLSSISSIIII&
ASASALILLSIISSIIII\\
AIASALSLLSSISSIIII&
AIASALSLLSIISSIIII&
AIASALILLSSISSIIII\\
AIASALILLSIISSIIII&
ALASASSLLISISSIIII&
ALASASSLLIIISSIIII\\
ALASASILLISISSIIII&
ALASASILLIIISSIIII&
ASASAISLLILISSIIII\\
ASASAIILLILISSIIII&
AIASAISLLILISSIIII&
AIASAIILLILISSIIII\\
ASAIALSLLSSISSIIII&
ASAIALSLLSIISSIIII&
ASAIALILLSSISSIIII\\
ASAIALILLSIISSIIII&
AIAIALSLLSSISSIIII&
AIAIALSLLSIISSIIII\\
AIAIALILLSSISSIIII&
AIAIALILLSIISSIIII&
ALAIASSLLISISSIIII\\
ALAIASSLLIIISSIIII&
ALAIASILLISISSIIII&
ALAIASILLIIISSIIII\\
ASAIAISLLILISSIIII&
ASAIAIILLILISSIIII&
AIAIAISLLILISSIIII\\
AIAIAIILLILISSIIII&
ASALALLSSSSISIIIII&
ASALALLISSSISIIIII\\
ASASALLLISSISIIIII&
ASAIALLLISSISIIIII&
AIALALLSSSSISIIIII\\
AIALALLISSSISIIIII&
AIASALLLISSISIIIII&
AIAIALLLISSISIIIII\\
ALALALSSSSIISIIIII&
ALALALSISSIISIIIII&
ALASALSLISIISIIIII\\
ALAIALSLISIISIIIII&
LLASSLAAIISSSIIIII&
LLASILAAIISSSIIIII\\
LLAISLAAIISSSIIIII&
LLAIILAAIISSSIIIII&
ASALALLSSISISIIIII\\
ASALALLISISISIIIII&
ASASALLLIISISIIIII&
ASAIALLLIISISIIIII\\
AIALALLSSISISIIIII&
AIALALLISISISIIIII&
AIASALLLIISISIIIII\\
AIAIALLLIISISIIIII&
LAASASSLLIIISIIIII&
LAAIASSLLIIISIIIII\\
ALAALSSSSSLIISIIII&
ALAALSSSSILIISIIII&
ASAALSISSLLIISIIII\\
AIAALSISSLLIISIIII&
ALAALISSSSLIISIIII&
ALAALISSSILIISIIII\\
ASAALIISSLLIISIIII&
AIAALIISSLLIISIIII&
ASLALALSSSSIIISIII\\
ASLALALISSSIIISIII&
AILALALSSSSIIISIII&
AILALALISSSIIISIII\\
ASSALAILSSILIISIII&
AISALAILSSILIISIII&
ASIALAILSSILIISIII\\
AIIALAILSSILIISIII&
ASAALLLSSSSIIIIIII&
AIAALLLSSSSIIIIIII\\
AALALLSSSSISIIIIII&
AALALLSISSISIIIIII&
AASALLILSSISIIIIII\\
AAIALLILSSISIIIIII&
AALALLSSSSIIIIIIII&
ASLALALSSISIIISIII\\
ASLALALISISIIISIII&
AILALALSSISIIISIII&
AILALALISISIIISIII\\
ASSALAILSIILIISIII&
AISALAILSIILIISIII&
ASIALAILSIILIISIII\\
AIIALAILSIILIISIII&
ASAALLLSSISIIIIIII&
AIAALLLSSISIIIIIII\\
ALSASSLASSISIILIIS&
ALSASSLASIISIILIIS&
ALIASSLASSISIILIIS\\
ALIASSLASIISIILIIS&
ALSAISLASSISIILIIS&
ALSAISLASIISIILIIS\\
ALIAISLASSISIILIIS&
ALIAISLASIISIILIIS&
ASLALSISSAIIISLIIS\\
ASLALSIISAIIISLIIS&
ASSALSIASLIIIILIIS&
ASIALSIASLIIIILIIS\\
AILALSISSAIIISLIIS&
AILALSIISAIIISLIIS&
AISALSIASLIIIILIIS\\
AIIALSIASLIIIILIIS&
ALSASILASSISIILIIS&
ALSASILASIISIILIIS\\
ALIASILASSISIILIIS&
ALIASILASIISIILIIS&
ALSAIILASSISIILIIS\\
ALSAIILASIISIILIIS&
ALIAIILASSISIILIIS&
ALIAIILASIISIILIIS\\
ASLALIISSAIIISLIIS&
ASLALIIISAIIISLIIS&
ASSALIIASLIIIILIIS\\
ASIALIIASLIIIILIIS&
AILALIISSAIIISLIIS&
AILALIIISAIIISLIIS\\
AISALIIASLIIIILIIS&
AIIALIIASLIIIILIIS&\\
\hline
\end{longtable}

\begin{longtable}{llllll}
\multicolumn{6}{c}{Table~\label{T:dim6P10}\ref{T:dim6P10}: $B_6$ case}\\
\hline
ASALALLSSSSSSISISII&
ASALALLSSSSSIISISII&
ASALALLISSSSSISISII\\
ASALALLISSSSIISISII&
AIALALLSSSSSSISISII&
AIALALLSSSSSIISISII\\
AIALALLISSSSSISISII&
AIALALLISSSSIISISII&
ALALALSSSSSISISISII\\
ALALALSSSSSIIISISII&
ALALALSISSSISISISII&
ALALALSISSSIIISISII\\
ASALALISSSSILISISII&
ASALALIISSSILISISII&
AIALALISSSSILISISII\\
AIALALIISSSILISISII&
ASALALLSSISSSISISII&
ASALALLSSISSIISISII\\
ASALALLISISSSISISII&
ASALALLISISSIISISII&
AIALALLSSISSSISISII\\
AIALALLSSISSIISISII&
AIALALLISISSSISISII&
AIALALLISISSIISISII\\
ALALALSSSISISISISII&
ALALALSSSISIIISISII&
ALALALSISISISISISII\\
ALALALSISISIIISISII&
ASALALISSISILISISII&
ASALALIISISILISISII\\
AIALALISSISILISISII&
AIALALIISISILISISII&
ASASALLSILSSSISISII\\
ASASALLSILSSIISISII&
ASASALLIILSSSISISII&
ASASALLIILSSIISISII\\
AIASALLSILSSSISISII&
AIASALLSILSSIISISII&
AIASALLIILSSSISISII\\
AIASALLIILSSIISISII&
ALASALSSILSISISISII&
ALASALSSILSIIISISII\\
ALASALSIILSISISISII&
ALASALSIILSIIISISII&
ASASALISILSILISISII\\
ASASALIIILSILISISII&
AIASALISILSILISISII&
AIASALIIILSILISISII\\
ASAIALLSILSSSISISII&
ASAIALLSILSSIISISII&
ASAIALLIILSSSISISII\\
ASAIALLIILSSIISISII&
AIAIALLSILSSSISISII&
AIAIALLSILSSIISISII\\
AIAIALLIILSSSISISII&
AIAIALLIILSSIISISII&
ALAIALSSILSISISISII\\
ALAIALSSILSIIISISII&
ALAIALSIILSISISISII&
ALAIALSIILSIIISISII\\
ASAIALISILSILISISII&
ASAIALIIILSILISISII&
AIAIALISILSILISISII\\
AIAIALIIILSILISISII&
ASALALLSSSSSSISIIII&
ASALALLISSSSSISIIII\\
AIALALLSSSSSSISIIII&
AIALALLISSSSSISIIII&
ASALALLSSISSSISIIII\\
ASALALLISISSSISIIII&
AIALALLSSISSSISIIII&
AIALALLISISSSISIIII\\
ASASALLLISSSSISIIII&
ASASALLLIISSSISIIII&
AIASALLLISSSSISIIII\\
AIASALLLIISSSISIIII&
ASAIALLLISSSSISIIII&
ASAIALLLIISSSISIIII\\
AIAIALLLISSSSISIIII&
AIAIALLLIISSSISIIII&
ASALALLSSSSSIISIIII\\
AIALALLSSSSSIISIIII&
ASALALLSSISSIISIIII&
AIALALLSSISSIISIIII\\
ALASALLSISSSIISIIII&
ALASALLSIISSIISIIII&
ALAIALLSISSSIISIIII\\
ALAIALLSIISSIISIIII&
ALALALSSSSSISISIIII&
ALALALSISSSISISIIII\\
ALALALSSSISISISIIII&
ALALALSISISISISIIII&
ASALALILSSSISISIIII\\
ASALALILSISISISIIII&
AIALALILSSSISISIIII&
AIALALILSISISISIIII\\
ASASALLLISSISISIIII&
ASASALLLIISISISIIII&
AIASALLLISSISISIIII\\
AIASALLLIISISISIIII&
ASAIALLLISSISISIIII&
ASAIALLLIISISISIIII\\
AIAIALLLISSISISIIII&
AIAIALLLIISISISIIII&
LAALALSSSSSIIISIIII\\
LAALALSSSISIIISIIII&
LAASALSSILSIIISIIII&
LAAIALSSILSIIISIIII\\
LALALASSSSSIISISIII&
LALALASSSISIISISIII&
LALASASSSLIIISISIII\\
LALAIASSSLIIISISIII&
LALALASISSSIISISIII&
LALALASISISIISISIII\\
LALASASISLIIISISIII&
LALAIASISLIIISISIII&
LASALASLISSIISISIII\\
LASALASLIISIISISIII&
LASASASLILIIISISIII&
LASAIASLILIIISISIII\\
LAIALASLISSIISISIII&
LAIALASLIISIISISIII&
LAIASASLILIIISISIII\\
LAIAIASLILIIISISIII&
LALALASSSSSIIIISIII&
LALALASSSISIIIISIII\\
LALASASSSLIIIIISIII&
LALAIASSSLIIIIISIII&
ASALALLSSSSSSIIISII\\
ASALALLSSSSSIIIISII&
ASALALLISSSSSIIISII&
ASALALLISSSSIIIISII\\
AIALALLSSSSSSIIISII&
AIALALLSSSSSIIIISII&
AIALALLISSSSSIIISII\\
AIALALLISSSSIIIISII&
ALALALSSSSSISIIISII&
ALALALSSSSSIIIIISII\\
ALALALSISSSISIIISII&
ALALALSISSSIIIIISII&
ASALALISSSSILIIISII\\
ASALALIISSSILIIISII&
AIALALISSSSILIIISII&
AIALALIISSSILIIISII\\
ASALALLSSSSSSIIIIII&
ASALALLISSSSSIIIIII&
AIALALLSSSSSSIIIIII\\
AIALALLISSSSSIIIIII&
ASALALLSSSSSIIIIIII&
AIALALLSSSSSIIIIIII\\
ALALALSSSSSISIIIIII&
ALALALSISSSISIIIIII&
ASALALILSSSISIIIIII\\
AIALALILSSSISIIIIII&
LALAALSSSSSIISIIIII&
LALAALSISSSIISIIIII\\
LASAALSLISSIISIIIII&
LAIAALSLISSIISIIIII&
LALAALSSSSSIIIIIIII\\
LALSASSSLASIIISIIIS&
LALSASSILASIIISIIIS&
LALIASSSLASIIISIIIS\\
LALIASSILASIIISIIIS&
LALSAISSLASIIISIIIS&
LALSAISILASIIISIIIS\\
LALIAISSLASIIISIIIS&
LALIAISILASIIISIIIS&
LASLALSAISSIISIIIIS\\
LASLALSAIISIISIIIIS&
LASSALSAILSIIIIIIIS&
LASIALSAILSIIIIIIIS\\
LAILALSAISSIISIIIIS&
LAILALSAIISIISIIIIS&
LAISALSAILSIIIIIIIS\\
LAIIALSAILSIIIIIIIS&\\
\hline
\end{longtable}

\begin{longtable}{llllll}
\multicolumn{6}{c}{Table~\label{T:dim6P11}\ref{T:dim6P11}: $B_7$ case}\\
\hline
LALALASSSSSSSIISIS&
LALALASSSSISSIISIS&
LALALAISSSSSSIISIS\\
LALALAISSSISSIISIS&
LALASALSSSSISIISIS&
LALASALSSSIISIISIS\\
LALAIALSSSSISIISIS&
LALAIALSSSIISIISIS&
LALALASSISSSSIISIS\\
LALALASSISISSIISIS&
LALALAISISSSSIISIS&
LALALAISISISSIISIS\\
LALASALSISSISIISIS&
LALASALSISIISIISIS&
LALAIALSISSISIISIS\\
LALAIALSISIISIISIS&
LASALALSSISSSIISIS&
LASALALSSISISIISIS\\
LASALALSIISSSIISIS&
LASALALSIISISIISIS&
LASASALSLIISSIISIS\\
LASASALSLIIISIISIS&
LASAIALSLIISSIISIS&
LASAIALSLIIISIISIS\\
LAIALALSSISSSIISIS&
LAIALALSSISISIISIS&
LAIALALSIISSSIISIS\\
LAIALALSIISISIISIS&
LAIASALSLIISSIISIS&
LAIASALSLIIISIISIS\\
LAIAIALSLIISSIISIS&
LAIAIALSLIIISIISIS&
LALALASSSSSSSIIIIS\\
LALALASSSSISSIIIIS&
LALALAISSSSSSIIIIS&
LALALAISSSISSIIIIS\\
LALASALSSSSISIIIIS&
LALASALSSSIISIIIIS&
LALAIALSSSSISIIIIS\\
LALAIALSSSIISIIIIS&
LALALASSSSSSSIISII&
LALALAISSSSSSIISII\\
LALALASSISSSSIISII&
LALALAISISSSSIISII&
LASALALSSISSSIISII\\
LASALALSIISSSIISII&
LAIALALSSISSSIISII&
LAIALALSIISSSIISII\\
LALALASSSSSSSIIIII&
LALALAISSSSSSIIIII&\\
\hline
\end{longtable}

\begin{longtable}{llllll}
\multicolumn{6}{c}{Table~\label{T:dim6P12}\ref{T:dim6P12}: $C_3$ case}\\
\hline
SAALASLLS&
IAALASLLS&
ALASAILLS&
SAAALLLSI&
IAAALLLSI\\\hline\end{longtable}

\begin{longtable}{llllll}
\multicolumn{6}{c}{Table~\label{T:dim6P13}\ref{T:dim6P13}: $C_4$ case}\\
\hline
SSSAIIAIIAILILLS&
ISSAIIAIIAILILLS&
SISAIIAIIAILILLS\\
IISAIIAIIAILILLS&
SSIAIIAIIAILILLS&
ISIAIIAIIAILILLS\\
SIIAIIAIIAILILLS&
IIIAIIAIIAILILLS&
SAALLALSSSSISIII\\
SAALLALSSSIISIII&
SAALLALISSSISIII&
SAALLALISSIISIII\\
SAALSALSISLISIII&
SAALSALIISLISIII&
SAALIALSISLISIII\\
SAALIALIISLISIII&
IAALLALSSSSISIII&
IAALLALSSSIISIII\\
IAALLALISSSISIII&
IAALLALISSIISIII&
IAALSALSISLISIII\\
IAALSALIISLISIII&
IAALIALSISLISIII&
IAALIALIISLISIII\\
LAASLALSSSSIIIII&
LAASLALISSSIIIII&
LAAILALSSSSIIIII\\
LAAILALISSSIIIII&
SAALSALLISSIIIII&
SAALIALLISSIIIII\\
IAALSALLISSIIIII&
IAALIALLISSIIIII&
SAALALSSLISIIIII\\
IAALALSSLISIIIII&
SAALALSILISIIIII&
IAALALSILISIIIII\\
ALLAAISSSLISIIII&
ALLAAISSILISIIII&
ALSAAIISLLISIIII\\
ALIAAIISLLISIIII&
LAASLALSSSIIIIII&
LAAILALSSSIIIIII\\
AAALLLSSSIIIIIII&
ASSAAILILLISIIII&
ASIAAILILLISIIII\\
AISAAILILLISIIII&
AIIAAILILLISIIII&
LAASLALISSIIIIII\\
LAAILALISSIIIIII&
SAASLASSLILISIII&
IAASLASSLILISIII\\
ALASSASILILISIII&
SAAILASSLILISIII&
IAAILASSLILISIII\\
ALAISASILILISIII&&
\\\hline\end{longtable}

\begin{longtable}{llllll}
\multicolumn{6}{c}{Table~\label{T:dim6P14}\ref{T:dim6P14}: $C_5$ case}\\
\hline
SALSALSLASIISIISIIIIII&
SALSALILASIISIISIIIIII\\
SALIALSLASIISIISIIIIII&
SALIALILASIISIISIIIIII\\
SASSALAILSILIIISIIIIII&
SASIALAILSILIIISIIIIII\\
SAISALAILSILIIISIIIIII&
SAIIALAILSILIIISIIIIII\\
IALSALSLASIISIISIIIIII&
IALSALILASIISIISIIIIII\\
IALIALSLASIISIISIIIIII&
IALIALILASIISIISIIIIII\\
IASSALAILSILIIISIIIIII&
IASIALAILSILIIISIIIIII\\
IAISALAILSILIIISIIIIII&
IAIIALAILSILIIISIIIIII\\
LALSASSLAIIISIISIIIIII&
LALSASILAIIISIISIIIIII\\
LALIASSLAIIISIISIIIIII&
LALIASILAIIISIISIIIIII\\
LASSASAILIILIIISIIIIII&
LASIASAILIILIIISIIIIII\\
LAISASAILIILIIISIIIIII&
LAIIASAILIILIIISIIIIII\\
SLALSIAALISSSIIISIIIII&
SLALIIAALISSSIIISIIIII\\
SLALSIAALIISSIIISIIIII&
SLALIIAALIISSIIISIIIII\\
ILALSIAALISSSIIISIIIII&
ILALIIAALISSSIIISIIIII\\
ILALSIAALIISSIIISIIIII&
ILALIIAALIISSIIISIIIII\\
SAASALSLLSIISIIIIIIIII&
SAAIALSLLSIISIIIIIIIII\\
IAASALSLLSIISIIIIIIIII&
IAAIALSLLSIISIIIIIIIII\\
LLASSLAAIISSSIIIIIIIII&
LLASILAAIISSSIIIIIIIII\\
LLAISLAAIISSSIIIIIIIII&
LLAIILAAIISSSIIIIIIIII\\
ASALALLSSISISIIIIIIIII&
ASALALLISISISIIIIIIIII\\
ASASALLLIISISIIIIIIIII&
ASAIALLLIISISIIIIIIIII\\
AIALALLSSISISIIIIIIIII&
AIALALLISISISIIIIIIIII\\
AIASALLLIISISIIIIIIIII&
AIAIALLLIISISIIIIIIIII\\
LAASASSLLIIISIIIIIIIII&
LAAIASSLLIIISIIIIIIIII\\
SSAALIASSILIILIISIIIII&
ISAALIASSILIILIISIIIII\\
SIAALIASSILIILIISIIIII&
IIAALIASSILIILIISIIIII\\
ALSALALSSSSIIISIIIIIII&
ALSALALISSSIIISIIIIIII\\
ALIALALSSSSIIISIIIIIII&
ALIALALISSSIIISIIIIIII\\
ALSALAISSSILIISIIIIIII&
ALSALAIISSILIISIIIIIII\\
ALIALAISSSILIISIIIIIII&
ALIALAIISSILIISIIIIIII\\
ALSALALSSSSIIIIIIIIIII&
ALIALALSSSSIIIIIIIIIII\\
ASLALALSSISIIISIIIIIII&
ASLALALISISIIISIIIIIII\\
AILALALSSISIIISIIIIIII&
AILALALISISIIISIIIIIII\\
ASSALAILSIILIISIIIIIII&
AISALAILSIILIISIIIIIII\\
ASIALAILSIILIISIIIIIII&
AIIALAILSIILIISIIIIIII\\
ASAALLLSSISIIIIIIIIIII&
AIAALLLSSISIIIIIIIIIII\\
SSASAIALSILISLIISIIIII&
ISASAIALSILISLIISIIIII\\
SIASAIALSILISLIISIIIII&
IIASAIALSILISLIISIIIII\\
SSASAIALIILISLIISIIIII&
ISASAIALIILISLIISIIIII\\
SIASAIALIILISLIISIIIII&
IIASAIALIILISLIISIIIII\\
SSAIAIALSILISLIISIIIII&
ISAIAIALSILISLIISIIIII\\
SIAIAIALSILISLIISIIIII&
IIAIAIALSILISLIISIIIII\\
SSAIAIALIILISLIISIIIII&
ISAIAIALIILISLIISIIIII\\
SIAIAIALIILISLIISIIIII&
IIAIAIALIILISLIISIIIII\\
SSSASIIASIIAIILIILILIS&
ISSASIIASIIAIILIILILIS\\
SISASIIASIIAIILIILILIS&
IISASIIASIIAIILIILILIS\\
SSIASIIASIIAIILIILILIS&
ISIASIIASIIAIILIILILIS\\
SIIASIIASIIAIILIILILIS&
IIIASIIASIIAIILIILILIS\\
SSSAIIIASIIAIILIILILIS&
ISSAIIIASIIAIILIILILIS\\
SISAIIIASIIAIILIILILIS&
IISAIIIASIIAIILIILILIS\\
SSIAIIIASIIAIILIILILIS&
ISIAIIIASIIAIILIILILIS\\
SIIAIIIASIIAIILIILILIS&
IIIAIIIASIIAIILIILILIS\\
\hline
\end{longtable}

\begin{longtable}{llllll}
\multicolumn{6}{c}{Table~\label{T:dim6P15}\ref{T:dim6P15}: $C_6$ case}\\
\hline
SAALALLSSSSSIISIIIIII&
IAALALLSSSSSIISIIIIII\\
SAALALLSSISSIISIIIIII&
IAALALLSSISSIISIIIIII\\
SAASALLSILSSIISIIIIII&
IAASALLSILSSIISIIIIII\\
SAAIALLSILSSIISIIIIII&
IAAIALLSILSSIISIIIIII\\
ALALALSSSSSISISIIIIII&
ALALALSISSSISISIIIIII\\
ALALALSSSISISISIIIIII&
ALALALSISISISISIIIIII\\
ASALALILSSSISISIIIIII&
ASALALILSISISISIIIIII\\
AIALALILSSSISISIIIIII&
AIALALILSISISISIIIIII\\
ASASALLLISSISISIIIIII&
ASASALLLIISISISIIIIII\\
AIASALLLISSISISIIIIII&
AIASALLLIISISISIIIIII\\
ASAIALLLISSISISIIIIII&
ASAIALLLIISISISIIIIII\\
AIAIALLLISSISISIIIIII&
AIAIALLLIISISISIIIIII\\
LAALALSSSSSIIISIIIIII&
LAALALSSSISIIISIIIIII\\
LAASALSSILSIIISIIIIII&
LAAIALSSILSIIISIIIIII\\
SALALALSSSSSISISIIIII&
SALALALSSISSISISIIIII\\
SALASALSSLISISISIIIII&
SALAIALSSLISISISIIIII\\
SALALALISSSSISISIIIII&
SALALALISISSISISIIIII\\
SALASALISLISISISIIIII&
SALAIALISLISISISIIIII\\
SASALALLISSSISISIIIII&
SASALALLIISSISISIIIII\\
SASASALLILISISISIIIII&
SASAIALLILISISISIIIII\\
SAIALALLISSSISISIIIII&
SAIALALLIISSISISIIIII\\
SAIASALLILISISISIIIII&
SAIAIALLILISISISIIIII\\
IALALALSSSSSISISIIIII&
IALALALSSISSISISIIIII\\
IALASALSSLISISISIIIII&
IALAIALSSLISISISIIIII\\
IALALALISSSSISISIIIII&
IALALALISISSISISIIIII\\
IALASALISLISISISIIIII&
IALAIALISLISISISIIIII\\
IASALALLISSSISISIIIII&
IASALALLIISSISISIIIII\\
IASASALLILISISISIIIII&
IASAIALLILISISISIIIII\\
IAIALALLISSSISISIIIII&
IAIALALLIISSISISIIIII\\
IAIASALLILISISISIIIII&
IAIAIALLILISISISIIIII\\
SALALALSSSSSIIISIIIII&
SALALALSSISSIIISIIIII\\
SALASALSSLISIIISIIIII&
SALAIALSSLISIIISIIIII\\
IALALALSSSSSIIISIIIII&
IALALALSSISSIIISIIIII\\
IALASALSSLISIIISIIIII&
IALAIALSSLISIIISIIIII\\
LALALASSSSSIISISIIIII&
LALALASSSISIISISIIIII\\
LALASASSSLIIISISIIIII&
LALAIASSSLIIISISIIIII\\
LALALASISSSIISISIIIII&
LALALASISISIISISIIIII\\
LALASASISLIIISISIIIII&
LALAIASISLIIISISIIIII\\
LASALASLISSIISISIIIII&
LASALASLIISIISISIIIII\\
LASASASLILIIISISIIIII&
LASAIASLILIIISISIIIII\\
LAIALASLISSIISISIIIII&
LAIALASLIISIISISIIIII\\
LAIASASLILIIISISIIIII&
LAIAIASLILIIISISIIIII\\
LALALASSSSSIIIISIIIII&
LALALASSSISIIIISIIIII\\
LALASASSSLIIIIISIIIII&
LALAIASSSLIIIIISIIIII\\
SALAALLSSSSSISIIIIIII&
SALAALLISSSSISIIIIIII\\
SASAALLLISSSISIIIIIII&
SAIAALLLISSSISIIIIIII\\
IALAALLSSSSSISIIIIIII&
IALAALLISSSSISIIIIIII\\
IASAALLLISSSISIIIIIII&
IAIAALLLISSSISIIIIIII\\
SAALALLSSSSSIIIIIIIII&
IAALALLSSSSSIIIIIIIII\\
ALALALSSSSSISIIIIIIII&
ALALALSISSSISIIIIIIII\\
ASALALILSSSISIIIIIIII&
AIALALILSSSISIIIIIIII\\
LALAALSSSSSIISIIIIIII&
LALAALSISSSIISIIIIIII\\
LASAALSLISSIISIIIIIII&
LAIAALSLISSIISIIIIIII\\
LALAALSSSSSIIIIIIIIII&
SASLALSALSSSISIIIISII\\
SASLALSALISSISIIIISII&
SAILALSALSSSISIIIISII\\
SAILALSALISSISIIIISII&
SASLALIALSSSISIIIISII\\
SASLALIALISSISIIIISII&
SAILALIALSSSISIIIISII\\
SAILALIALISSISIIIISII&
IASLALSALSSSISIIIISII\\
IASLALSALISSISIIIISII&
IAILALSALSSSISIIIISII\\
IAILALSALISSISIIIISII&
IASLALIALSSSISIIIISII\\
IASLALIALISSISIIIISII&
IAILALIALSSSISIIIISII\\
IAILALIALISSISIIIISII&
SALSALLSIASSIISIIISII\\
SALSALLIIASSIISIIISII&
SASSALLAILSSIIIIIISII\\
SAISALLAILSSIIIIIISII&
IALSALLSIASSIISIIISII\\
IALSALLIIASSIISIIISII&
IASSALLAILSSIIIIIISII\\
IAISALLAILSSIIIIIISII&
SALIALLSIASSIISIIISII\\
SALIALLIIASSIISIIISII&
SASIALLAILSSIIIIIISII\\
SAIIALLAILSSIIIIIISII&
IALIALLSIASSIISIIISII\\
IALIALLIIASSIISIIISII&
IASIALLAILSSIIIIIISII\\
IAIIALLAILSSIIIIIISII&
LALSASSSLASIIISIIISII\\
LALSASSILASIIISIIISII&
LALIASSSLASIIISIIISII\\
LALIASSILASIIISIIISII&
LALSAISSLASIIISIIISII\\
LALSAISILASIIISIIISII&
LALIAISSLASIIISIIISII\\
LALIAISILASIIISIIISII&
LASLALSAISSIISIIIISII\\
LASLALSAIISIISIIIISII&
LASSALSAILSIIIIIIISII\\
LASIALSAILSIIIIIIISII&
LAILALSAISSIISIIIISII\\
LAILALSAIISIISIIIISII&
LAISALSAILSIIIIIIISII\\
LAIIALSAILSIIIIIIISII&
SSASALIASSSILISILIIIS\\
ISASALIASSSILISILIIIS&
SIASALIASSSILISILIIIS\\
IIASALIASSSILISILIIIS&
SSASALIASISILISILIIIS\\
ISASALIASISILISILIIIS&
SIASALIASISILISILIIIS\\
IIASALIASISILISILIIIS&
SSAIALIASSSILISILIIIS\\
ISAIALIASSSILISILIIIS&
SIAIALIASSSILISILIIIS\\
IIAIALIASSSILISILIIIS&
SSAIALIASISILISILIIIS\\
ISAIALIASISILISILIIIS&
SIAIALIASISILISILIIIS\\
IIAIALIASISILISILIIIS&
SSASALIAISSILISILIIIS\\
ISASALIAISSILISILIIIS&
SIASALIAISSILISILIIIS\\
IIASALIAISSILISILIIIS&
SSASALIAIISILISILIIIS\\
ISASALIAIISILISILIIIS&
SIASALIAIISILISILIIIS\\
IIASALIAIISILISILIIIS&
SSAIALIAISSILISILIIIS\\
ISAIALIAISSILISILIIIS&
SIAIALIAISSILISILIIIS\\
IIAIALIAISSILISILIIIS&
SSAIALIAIISILISILIIIS\\
ISAIALIAIISILISILIIIS&
SIAIALIAIISILISILIIIS\\
IIAIALIAIISILISILIIIS&
SSASALIASSSILIIILIIIS\\
ISASALIASSSILIIILIIIS&
SIASALIASSSILIIILIIIS\\
IIASALIASSSILIIILIIIS&
SSAIALIASSSILIIILIIIS\\
ISAIALIASSSILIIILIIIS&
SIAIALIASSSILIIILIIIS\\
IIAIALIASSSILIIILIIIS&\\
\hline
\end{longtable}

\begin{longtable}{llllll}
\multicolumn{6}{c}{Table~\label{T:dim6P16}\ref{T:dim6P16}: $C_7$ case}\\
\hline
SALALALSSSSSSSISIS&
SALALALSSISSSSISIS&
SALALALISSSSSSISIS\\
SALALALISISSSSISIS&
IALALALSSSSSSSISIS&
IALALALSSISSSSISIS\\
IALALALISSSSSSISIS&
IALALALISISSSSISIS&
SALALALSSSSISSISIS\\
SALALALSSISISSISIS&
SALALALISSSISSISIS&
SALALALISISISSISIS\\
IALALALSSSSISSISIS&
IALALALSSISISSISIS&
IALALALISSSISSISIS\\
IALALALISISISSISIS&
SASALALLISSSSSISIS&
SASALALLIISSSSISIS\\
IASALALLISSSSSISIS&
IASALALLIISSSSISIS&
SASALALLISSISSISIS\\
SASALALLIISISSISIS&
IASALALLISSISSISIS&
IASALALLIISISSISIS\\
SAIALALLISSSSSISIS&
SAIALALLIISSSSISIS&
IAIALALLISSSSSISIS\\
IAIALALLIISSSSISIS&
SAIALALLISSISSISIS&
SAIALALLIISISSISIS\\
IAIALALLISSISSISIS&
IAIALALLIISISSISIS&
SALASALLSSISSSISIS\\
SALASALLSSIISSISIS&
SALASALLISISSSISIS&
SALASALLISIISSISIS\\
SASASALLLIISSSISIS&
SASASALLLIIISSISIS&
SAIASALLLIISSSISIS\\
SAIASALLLIIISSISIS&
IALASALLSSISSSISIS&
IALASALLSSIISSISIS\\
IALASALLISISSSISIS&
IALASALLISIISSISIS&
IASASALLLIISSSISIS\\
IASASALLLIIISSISIS&
IAIASALLLIISSSISIS&
IAIASALLLIIISSISIS\\
SALAIALLSSISSSISIS&
SALAIALLSSIISSISIS&
SALAIALLISISSSISIS\\
SALAIALLISIISSISIS&
SASAIALLLIISSSISIS&
SASAIALLLIIISSISIS\\
SAIAIALLLIISSSISIS&
SAIAIALLLIIISSISIS&
IALAIALLSSISSSISIS\\
IALAIALLSSIISSISIS&
IALAIALLISISSSISIS&
IALAIALLISIISSISIS\\
IASAIALLLIISSSISIS&
IASAIALLLIIISSISIS&
IAIAIALLLIISSSISIS\\
IAIAIALLLIIISSISIS&
SALALALSSSSSSSIIIS&
SALALALISSSSSSIIIS\\
IALALALSSSSSSSIIIS&
IALALALISSSSSSIIIS&
SALALALSSSSISSIIIS\\
SALALALISSSISSIIIS&
IALALALSSSSISSIIIS&
IALALALISSSISSIIIS\\
SALASALLSSISSSIIIS&
SALASALLSSIISSIIIS&
IALASALLSSISSSIIIS\\
IALASALLSSIISSIIIS&
SALAIALLSSISSSIIIS&
SALAIALLSSIISSIIIS\\
IALAIALLSSISSSIIIS&
IALAIALLSSIISSIIIS&
LALALASSSSSSSIISIS\\
LALALASSSSISSIISIS&
LALALAISSSSSSIISIS&
LALALAISSSISSIISIS\\
LALASALSSSSISIISIS&
LALASALSSSIISIISIS&
LALAIALSSSSISIISIS\\
LALAIALSSSIISIISIS&
LALALASSISSSSIISIS&
LALALASSISISSIISIS\\
LALALAISISSSSIISIS&
LALALAISISISSIISIS&
LALASALSISSISIISIS\\
LALASALSISIISIISIS&
LALAIALSISSISIISIS&
LALAIALSISIISIISIS\\
LASALALSSISSSIISIS&
LASALALSSISISIISIS&
LASALALSIISSSIISIS\\
LASALALSIISISIISIS&
LASASALSLIISSIISIS&
LASASALSLIIISIISIS\\
LASAIALSLIISSIISIS&
LASAIALSLIIISIISIS&
LAIALALSSISSSIISIS\\
LAIALALSSISISIISIS&
LAIALALSIISSSIISIS&
LAIALALSIISISIISIS\\
LAIASALSLIISSIISIS&
LAIASALSLIIISIISIS&
LAIAIALSLIISSIISIS\\
LAIAIALSLIIISIISIS&
LALALASSSSSSSIIIIS&
LALALASSSSISSIIIIS\\
LALALAISSSSSSIIIIS&
LALALAISSSISSIIIIS&
LALASALSSSSISIIIIS\\
LALASALSSSIISIIIIS&
LALAIALSSSSISIIIIS&
LALAIALSSSIISIIIIS\\
SALALALSSSSSSSISII&
SALALALSSISSSSISII&
SALALALISSSSSSISII\\
SALALALISISSSSISII&
IALALALSSSSSSSISII&
IALALALSSISSSSISII\\
IALALALISSSSSSISII&
IALALALISISSSSISII&
SASALALLISSSSSISII\\
SASALALLIISSSSISII&
IASALALLISSSSSISII&
IASALALLIISSSSISII\\
SAIALALLISSSSSISII&
SAIALALLIISSSSISII&
IAIALALLISSSSSISII\\
IAIALALLIISSSSISII&
SALALALSSSSSSSIIII&
SALALALISSSSSSIIII\\
IALALALSSSSSSSIIII&
IALALALISSSSSSIIII&
LALALASSSSSSSIISII\\
LALALAISSSSSSIISII&
LALALASSISSSSIISII&
LALALAISISSSSIISII\\
LASALALSSISSSIISII&
LASALALSIISSSIISII&
LAIALALSSISSSIISII\\
LAIALALSIISSSIISII&
LALALASSSSSSSIIIII&
LALALAISSSSSSIIIII\\
\hline
\end{longtable}

\begin{longtable}{llllll}
\multicolumn{6}{c}{Table~\label{T:dim6P17}\ref{T:dim6P17}: $D_4$ case}\\
\hline
AALASSSSLLS&
AALASSSILLS&
AALAISSSLLS&
AALAISSILLS\\
AALALSISSLS&
AALALSIISLS&
AALASSILILS&
AALAISILILS\\
AALALISSSLS&
AALALISSILS&
AALASISILLS&
AALAIISILLS\\
AASASIILLLS&
AASAIIILLLS&
AAIASIILLLS&
AAIAIIILLLS\\
AALALLSSSII&
AALALLSSISI&
AAALLLSSIII&
AALALSLISSI\\
ALAALSLISII&
LAAASLLIISI&
AASALLLSSSI&
\\\hline\end{longtable}

\begin{longtable}{llllll}
\multicolumn{6}{c}{Table~\label{T:dim6P18}\ref{T:dim6P18}: $D_5$ case}\\
\hline
AASALLLSSSSSISII&
AASALLLSSSISISII&
AASALLLISSSSISII\\
AASALLLISSISISII&
AASALLLSSSSIISII&
AASALLLSSSIIISII\\
AASALLLISSSIISII&
AASALLLISSIIISII&
AAIALLLSSSSSISII\\
AAIALLLSSSISISII&
AAIALLLISSSSISII&
AAIALLLISSISISII\\
AAIALLLSSSSIISII&
AAIALLLSSSIIISII&
AAIALLLISSSIISII\\
AAIALLLISSIIISII&
AALALSSSSISLISII&
AALALSSSSIILISII\\
AALALSSISISLISII&
AALALSSISIILISII&
AALALSILSISSISII\\
AALALSILSIISISII&
AALALSISSILIISII&
AALALSIISILIISII\\
AALALISLSISSISII&
AALALISLSISIISII&
AALALISSSIILISII\\
AALALISISIILISII&
AASALIISSILLISII&
AASALIIISILLISII\\
AAIALIISSILLISII&
AAIALIIISILLISII&
ALSLAALISSSIIISI\\
ALSLAALISSIIIISI&
ALSSAALIISLIIISI&
ALSIAALIISLIIISI\\
ALILAALISSSIIISI&
ALILAALISSIIIISI&
ALISAALIISLIIISI\\
ALIIAALIISLIIISI&
AASALLLSSSSSIIII&
AASALLLISSSSIIII\\
AAIALLLSSSSSIIII&
AAIALLLISSSSIIII&
AASALLLSSSSIIIII\\
AASALLLISSSIIIII&
AAIALLLSSSSIIIII&
AAIALLLISSSIIIII\\
LASLALAISSISIIIS&
LASLALAISSIIIIIS&
LASSALAIISILIIIS\\
LASIALAIISILIIIS&
LAILALAISSISIIIS&
LAILALAISSIIIIIS\\
LAISALAIISILIIIS&
LAIIALAIISILIIIS&
AASALLLSSSISIIII\\
AASALLLISSISIIII&
AAIALLLSSSISIIII&
AAIALLLISSISIIII\\
AASALLLSSSIIIIII&
AAIALLLSSSIIIIII&
LSSSALIIAILIAISS\\
LISSALIIAILIAISS&
LSISALIIAILIAISS&
LIISALIIAILIAISS\\
LSSIALIIAILIAISS&
LISIALIIAILIAISS&
LSIIALIIAILIAISS\\
LIIIALIIAILIAISS&
LLAASLSASISSSIII&
LLAASLSAIISSSIII\\
LLAAILSASISSSIII&
LLAAILSAIISSSIII&
LSAASLIALISSSIII\\
LSAAILIALISSSIII&
LIAASLIALISSSIII&
LIAAILIALISSSIII\\
LLAASLSASISISIII&
LLAAILSASISISIII&
LLAASLSAIISISIII\\
LLAAILSAIISISIII&
ALAALSSLSISIIIII&
ALAALISLSISIIIII\\
LLAASSLASIISSIII&
LLAAISLASIISSIII&
LLAASSLAIIISSIII\\
LLAAISLAIIISSIII&
LAAALSSLSIISIIII&
LAAALSILSIISIIII\\
AALASLSLSISSIIII&
AALAILSLSISSIIII&\\
\hline
\end{longtable}

\begin{longtable}{llllll}
\multicolumn{6}{c}{Table~\label{T:dim6P19}\ref{T:dim6P19}: $D_6$ case}\\
\hline
ALSSALALISSSSIIIISII&
ALSSALALIISSSIIIISII\\
ALSIALALISSSSIIIISII&
ALSIALALIISSSIIIISII\\
ALISALALISSSSIIIISII&
ALISALALIISSSIIIISII\\
ALIIALALISSSSIIIISII&
ALIIALALIISSSIIIISII\\
ALSSALALISSSIIIIISII&
ALSSALALIISSIIIIISII\\
ALISALALISSSIIIIISII&
ALISALALIISSIIIIISII\\
ALSIALALISSSIIIIISII&
ALSIALALIISSIIIIISII\\
ALIIALALISSSIIIIISII&
ALIIALALIISSIIIIISII\\
ALSSALALISSSSIIIIIII&
ALSIALALISSSSIIIIIII\\
ALISALALISSSSIIIIIII&
ALIIALALISSSSIIIIIII\\
LASSALLAISSSISIIIISI&
LASSALLAIISSISIIIISI\\
LASIALLAISSSISIIIISI&
LASIALLAIISSISIIIISI\\
LAISALLAISSSISIIIISI&
LAISALLAIISSISIIIISI\\
LAIIALLAISSSISIIIISI&
LAIIALLAIISSISIIIISI\\
LASSALLAISSSIIIIIISI&
LASSALLAIISSIIIIIISI\\
LAISALLAISSSIIIIIISI&
LAISALLAIISSIIIIIISI\\
LASIALLAISSSIIIIIISI&
LASIALLAIISSIIIIIISI\\
LAIIALLAISSSIIIIIISI&
LAIIALLAIISSIIIIIISI\\
LASSALLAISSSISIIIIII&
LASIALLAISSSISIIIIII\\
LAISALLAISSSISIIIIII&
LAIIALLAISSSISIIIIII\\
LSLASSLIAAIISSSSIIII&
LSLASILIAAIISSSSIIII\\
LSLAISLIAAIISSSSIIII&
LSLAIILIAAIISSSSIIII\\
LILASSLIAAIISSSSIIII&
LILASILIAAIISSSSIIII\\
LILAISLIAAIISSSSIIII&
LILAIILIAAIISSSSIIII\\
LSSALALIASSISSISIIII&
LISALALIASSISSISIIII\\
LSSALALIAISISSISIIII&
LISALALIAISISSISIIII\\
LSSASALIALIISSISIIII&
LISASALIALIISSISIIII\\
LSSAIALIALIISSISIIII&
LISAIALIALIISSISIIII\\
LSIALALIASSISSISIIII&
LIIALALIASSISSISIIII\\
LSIALALIAISISSISIIII&
LIIALALIAISISSISIIII\\
LSIASALIALIISSISIIII&
LIIASALIALIISSISIIII\\
LSIAIALIALIISSISIIII&
LIIAIALIALIISSISIIII\\
LSSAALLIASSISSIIIIII&
LISAALLIASSISSIIIIII\\
LSIAALLIASSISSIIIIII&
LIIAALLIASSISSIIIIII\\
LSALSASSSLAIIIISIIIS&
LSALSASSILAIIIISIIIS\\
LSALIASSSLAIIIISIIIS&
LSALIASSILAIIIISIIIS\\
LSASSASSAILIIILIIIIS&
LSASIASSAILIIILIIIIS\\
LSAISASSAILIIILIIIIS&
LSAIIASSAILIIILIIIIS\\
LIALSASSSLAIIIISIIIS&
LIALSASSILAIIIISIIIS\\
LIALIASSSLAIIIISIIIS&
LIALIASSILAIIIISIIIS\\
LIASSASSAILIIILIIIIS&
LIASIASSAILIIILIIIIS\\
LIAISASSAILIIILIIIIS&
LIAIIASSAILIIILIIIIS\\
LSAASASSSLLIIIISIIII&
LSAAIASSSLLIIIISIIII\\
LIAASASSSLLIIIISIIII&
LIAAIASSSLLIIIISIIII\\
LSSSASLIIASILIAIISSI&
LISSASLIIASILIAIISSI\\
LSISASLIIASILIAIISSI&
LIISASLIIASILIAIISSI\\
LSSIASLIIASILIAIISSI&
LISIASLIIASILIAIISSI\\
LSIIASLIIASILIAIISSI&
LIIIASLIIASILIAIISSI\\
LSSSAILIIASILIAIISSI&
LISSAILIIASILIAIISSI\\
LSISAILIIASILIAIISSI&
LIISAILIIASILIAIISSI\\
LSSIAILIIASILIAIISSI&
LISIAILIIASILIAIISSI\\
LSIIAILIIASILIAIISSI&
LIIIAILIIASILIAIISSI\\
\hline
\end{longtable}

\begin{longtable}{llllll}
\multicolumn{6}{c}{Table~\label{T:dim6P20}\ref{T:dim6P20}: $D_7$ case}\\
\hline
LSALALASSSSSSIIISISIIII&
LSALALASSSSISIIISISIIII\\
LSALASASSSSLIIIISISIIII&
LSALAIASSSSLIIIISISIIII\\
LSALALASSISSSIIISISIIII&
LSALALASSISISIIISISIIII\\
LSALASASSISLIIIISISIIII&
LSALAIASSISLIIIISISIIII\\
LSASALASSLISSIIISISIIII&
LSASALASSLIISIIISISIIII\\
LSASASASSLILIIIISISIIII&
LSASAIASSLILIIIISISIIII\\
LSAIALASSLISSIIISISIIII&
LSAIALASSLIISIIISISIIII\\
LSAIASASSLILIIIISISIIII&
LSAIAIASSLILIIIISISIIII\\
LIALALASSSSSSIIISISIIII&
LIALALASSSSISIIISISIIII\\
LIALASASSSSLIIIISISIIII&
LIALAIASSSSLIIIISISIIII\\
LIALALASSISSSIIISISIIII&
LIALALASSISISIIISISIIII\\
LIALASASSISLIIIISISIIII&
LIALAIASSISLIIIISISIIII\\
LIASALASSLISSIIISISIIII&
LIASALASSLIISIIISISIIII\\
LIASASASSLILIIIISISIIII&
LIASAIASSLILIIIISISIIII\\
LIAIALASSLISSIIISISIIII&
LIAIALASSLIISIIISISIIII\\
LIAIASASSLILIIIISISIIII&
LIAIAIASSLILIIIISISIIII\\
LSALAALSSSSSSIIISIIIIII&
LSALAALSSISSSIIISIIIIII\\
LSASAALSSLISSIIISIIIIII&
LSAIAALSSLISSIIISIIIIII\\
LIALAALSSSSSSIIISIIIIII&
LIALAALSSISSSIIISIIIIII\\
LIASAALSSLISSIIISIIIIII&
LIAIAALSSLISSIIISIIIIII\\
LSAALALSSSSSSIIIISIIIII&
LSAALALSSSSISIIIISIIIII\\
LSAASALSSSILSIIIISIIIII&
LSAAIALSSSILSIIIISIIIII\\
LIAALALSSSSSSIIIISIIIII&
LIAALALSSSSISIIIISIIIII\\
LIAASALSSSILSIIIISIIIII&
LIAAIALSSSILSIIIISIIIII\\
LSALALASSSSSSIIIIISIIII&
LSALALASSSSISIIIIISIIII\\
LSALASASSSSLIIIIIISIIII&
LSALAIASSSSLIIIIIISIIII\\
LIALALASSSSSSIIIIISIIII&
LIALALASSSSISIIIIISIIII\\
LIALASASSSSLIIIIIISIIII&
LIALAIASSSSLIIIIIISIIII\\
LSALAALSSSSSSIIIIIIIIII&
LIALAALSSSSSSIIIIIIIIII\\
LSSASALLIASSSISSISIIIII&
LISASALLIASSSISSISIIIII\\
LSSASALLIASISISSISIIIII&
LISASALLIASISISSISIIIII\\
LSIASALLIASSSISSISIIIII&
LIIASALLIASSSISSISIIIII\\
LSIASALLIASISISSISIIIII&
LIIASALLIASISISSISIIIII\\
LSSAIALLIASSSISSISIIIII&
LISAIALLIASSSISSISIIIII\\
LSSAIALLIASISISSISIIIII&
LISAIALLIASISISSISIIIII\\
LSIAIALLIASSSISSISIIIII&
LIIAIALLIASSSISSISIIIII\\
LSIAIALLIASISISSISIIIII&
LIIAIALLIASISISSISIIIII\\
LSSASALLIAISSISSISIIIII&
LISASALLIAISSISSISIIIII\\
LSSASALLIAIISISSISIIIII&
LISASALLIAIISISSISIIIII\\
LSIASALLIAISSISSISIIIII&
LIIASALLIAISSISSISIIIII\\
LSIASALLIAIISISSISIIIII&
LIIASALLIAIISISSISIIIII\\
LSSAIALLIAISSISSISIIIII&
LISAIALLIAISSISSISIIIII\\
LSSAIALLIAIISISSISIIIII&
LISAIALLIAIISISSISIIIII\\
LSIAIALLIAISSISSISIIIII&
LIIAIALLIAISSISSISIIIII\\
LSIAIALLIAIISISSISIIIII&
LIIAIALLIAIISISSISIIIII\\
LSSASALLIASSSISSIIIIIII&
LISASALLIASSSISSIIIIIII\\
LSIASALLIASSSISSIIIIIII&
LIIASALLIASSSISSIIIIIII\\
LSSAIALLIASSSISSIIIIIII&
LISAIALLIASSSISSIIIIIII\\
LSIAIALLIASSSISSIIIIIII&
LIIAIALLIASSSISSIIIIIII\\
LSALSASSSSLASIIIISIIIIS&
LSALSASSSILASIIIISIIIIS\\
LSALIASSSSLASIIIISIIIIS&
LSALIASSSILASIIIISIIIIS\\
LSALSAISSSLASIIIISIIIIS&
LSALSAISSILASIIIISIIIIS\\
LSALIAISSSLASIIIISIIIIS&
LSALIAISSILASIIIISIIIIS\\
LSASLALSSAISSIIISIIIIIS&
LSASLALSSAIISIIISIIIIIS\\
LSASSALSSAILSIIIIIIIIIS&
LSASIALSSAILSIIIIIIIIIS\\
LSAILALSSAISSIIISIIIIIS&
LSAILALSSAIISIIISIIIIIS\\
LSAISALSSAILSIIIIIIIIIS&
LSAIIALSSAILSIIIIIIIIIS\\
LIALSASSSSLASIIIISIIIIS&
LIALSASSSILASIIIISIIIIS\\
LIALIASSSSLASIIIISIIIIS&
LIALIASSSILASIIIISIIIIS\\
LIALSAISSSLASIIIISIIIIS&
LIALSAISSILASIIIISIIIIS\\
LIALIAISSSLASIIIISIIIIS&
LIALIAISSILASIIIISIIIIS\\
LIASLALSSAISSIIISIIIIIS&
LIASLALSSAIISIIISIIIIIS\\
LIASSALSSAILSIIIIIIIIIS&
LIASIALSSAILSIIIIIIIIIS\\
LIAILALSSAISSIIISIIIIIS&
LIAILALSSAIISIIISIIIIIS\\
LIAISALSSAILSIIIIIIIIIS&
LIAIIALSSAILSIIIIIIIIIS\\
\hline
\end{longtable}

\begin{longtable}{llllll}
\multicolumn{6}{c}{Table~\label{T:dim6P21}\ref{T:dim6P21}: $D_8$ case}\\
\hline
LSALALASSSSSSSSIIISIS&
LSALALASSSSSISSIIISIS\\
LSALALAISSSSSSSIIISIS&
LSALALAISSSSISSIIISIS\\
LSALASALSSSSSISIIISIS&
LSALASALSSSSIISIIISIS\\
LSALAIALSSSSSISIIISIS&
LSALAIALSSSSIISIIISIS\\
LSALALASSSISSSSIIISIS&
LSALALASSSISISSIIISIS\\
LSALALAISSISSSSIIISIS&
LSALALAISSISISSIIISIS\\
LSALASALSSISSISIIISIS&
LSALASALSSISIISIIISIS\\
LSALAIALSSISSISIIISIS&
LSALAIALSSISIISIIISIS\\
LSASALALSSSISSSIIISIS&
LSASALALSSSISISIIISIS\\
LSASALALSSIISSSIIISIS&
LSASALALSSIISISIIISIS\\
LSASASALSSLIISSIIISIS&
LSASASALSSLIIISIIISIS\\
LSASAIALSSLIISSIIISIS&
LSASAIALSSLIIISIIISIS\\
LSAIALALSSSISSSIIISIS&
LSAIALALSSSISISIIISIS\\
LSAIALALSSIISSSIIISIS&
LSAIALALSSIISISIIISIS\\
LSAIASALSSLIISSIIISIS&
LSAIASALSSLIIISIIISIS\\
LSAIAIALSSLIISSIIISIS&
LSAIAIALSSLIIISIIISIS\\
LIALALASSSSSSSSIIISIS&
LIALALASSSSSISSIIISIS\\
LIALALAISSSSSSSIIISIS&
LIALALAISSSSISSIIISIS\\
LIALASALSSSSSISIIISIS&
LIALASALSSSSIISIIISIS\\
LIALAIALSSSSSISIIISIS&
LIALAIALSSSSIISIIISIS\\
LIALALASSSISSSSIIISIS&
LIALALASSSISISSIIISIS\\
LIALALAISSISSSSIIISIS&
LIALALAISSISISSIIISIS\\
LIALASALSSISSISIIISIS&
LIALASALSSISIISIIISIS\\
LIALAIALSSISSISIIISIS&
LIALAIALSSISIISIIISIS\\
LIASALALSSSISSSIIISIS&
LIASALALSSSISISIIISIS\\
LIASALALSSIISSSIIISIS&
LIASALALSSIISISIIISIS\\
LIASASALSSLIISSIIISIS&
LIASASALSSLIIISIIISIS\\
LIASAIALSSLIISSIIISIS&
LIASAIALSSLIIISIIISIS\\
LIAIALALSSSISSSIIISIS&
LIAIALALSSSISISIIISIS\\
LIAIALALSSIISSSIIISIS&
LIAIALALSSIISISIIISIS\\
LIAIASALSSLIISSIIISIS&
LIAIASALSSLIIISIIISIS\\
LIAIAIALSSLIISSIIISIS&
LIAIAIALSSLIIISIIISIS\\
LSALALASSSSSSSSIIISII&
LSALALAISSSSSSSIIISII\\
LSALALASSSISSSSIIISII&
LSALALAISSISSSSIIISII\\
LSASALALSSSISSSIIISII&
LSASALALSSIISSSIIISII\\
LSAIALALSSSISSSIIISII&
LSAIALALSSIISSSIIISII\\
LIALALASSSSSSSSIIISII&
LIALALAISSSSSSSIIISII\\
LIALALASSSISSSSIIISII&
LIALALAISSISSSSIIISII\\
LIASALALSSSISSSIIISII&
LIASALALSSIISSSIIISII\\
LIAIALALSSSISSSIIISII&
LIAIALALSSIISSSIIISII\\
LSALALASSSSSSSSIIIIIS&
LSALALASSSSSISSIIIIIS\\
LSALALAISSSSSSSIIIIIS&
LSALALAISSSSISSIIIIIS\\
LSALASALSSSSSISIIIIIS&
LSALASALSSSSIISIIIIIS\\
LSALAIALSSSSSISIIIIIS&
LSALAIALSSSSIISIIIIIS\\
LIALALASSSSSSSSIIIIIS&
LIALALASSSSSISSIIIIIS\\
LIALALAISSSSSSSIIIIIS&
LIALALAISSSSISSIIIIIS\\
LIALASALSSSSSISIIIIIS&
LIALASALSSSSIISIIIIIS\\
LIALAIALSSSSSISIIIIIS&
LIALAIALSSSSIISIIIIIS\\
LSALALASSSSSSSSIIIIII&
LSALALAISSSSSSSIIIIII\\
LIALALASSSSSSSSIIIIII&
LIALALAISSSSSSSIIIIII\\
\hline
\end{longtable}

\sect{The number of irreducible characters for finite unipotent groups\label{sect:number_chars}}\addcontentsline{toc}{subsection}{\ref{sect:number_chars}. The number of irreducible characters for finite unipotent groups}

Throughout this section, $\Fp_q$ is the finite field with $q$ elements of sufficiently large characteristic~$p$ (in fact, we will assume that $p\geq n$). We also assume that $\Fp$ is the algebraic closure of $\Fp_q$. One can define $\gt(q)$, $\bt(q)$, $\nt(q)$, $\nt^*(q)$ and $N(q)$ as $\gt$, $\bt$, $\nt$, $\nt^*$ and $N$ from Section~\ref{Sfdcn}
respectively using $\Fp_q$ instead of $\Fp$. Note that $N(q)$ is a subgroup of $N$, while $\nt^*(q)$ is an $\Fp_q$-subspaces of $\nt^*$. Note also that $N(q)$ is a Sylow $p$-subgroup in the corresponding classical Chevalley group $G(q)$ over $\Fp_q$. Given a coadjoint $N(q)$-orbit $\Omega(q)\subset\nt^*(q)$, we denote by $\Omega\subset\nt^*$ the coadjoint $N$-orbit of an arbitrary linear form from~$\Omega(q)$. It is well known that if $\dim\Omega=2e$ then $|\Omega(q)|=q^{2e}$. 

Fix a nontrivial group homomorphism $\theta\colon\Fp_q\to\Cp^{\times}$, where $\Cp^{\times}$ is the multiplicative group of $\Cp$. Let $\ln\colon N(q)\to\nt(q)$ be the inverse map to the exponential map. Given a coadjoint $N(q)$-orbit $\Omega(q)\subset\nt^*(q)$, put
$$\chi_{\Omega(q)}(g)=\dfrac{1}{\sqrt{|\Omega(q)|}}\sum_{\lambda\in\Omega(q)}\theta(\lambda(\ln g)),~g\in N(q).$$
The orbit method claims that
the map $\Omega(q)\mapsto\chi_{\Omega(q)}$ establishes a bijection between the set of coadjoint $N(q)$-orbits on $\nt^*(q)$ and the set of irreducible complex characters of the group $N(q)$. It is clear that if $\dim\Omega=2e$ then the \emph{degree} $\deg\chi_{\Omega(q)}$ of $\chi_{\Omega(q)}$ (i.e., the complex dimension of the corresponding irreducible representation) equals $q^e$.

A longstanding G. Higman's conjecture \cite{Higman60} states that the number of conjugacy classes for~$U_n(q)$, the group of all upper-triangular matrices over $\Fp_q$ with 1's on the diagonal, is a polynomial in~$q$ (note that $N(q)=U_n(q)$ if $\Phi=A_{n-1}$). It is easy to check that the number of conjugacy classes for $N(q)$ coincides with the number of coadjoint $N(q)$-orbits, or, equivalently, with the number of irreducible characters of $N(q)$. Of course, it is interesting to consider an analogue of Higman's conjecture for an arbitrary $N(q)$, not only for $U_n(q)$; we will denote the number of irreducible characters of the group $N(q)$ by $O(q)$. Forty years after Higman, G. Lehrer conjectured \cite{Lehrer74} that, for $N(q)=U_n(q)$, even the number $O_e(q)$ of characters of degree $q^e$ is polynomial in $q$. In 2007, M. Isaacs made a stronger conjecture that, for $N(q)=U_n(q)$, $O_e(q)$ is in fact polynomial in $q-1$ with nonnegative integer coefficients.

In past twenty five years, a significant progress in studying these conjectures has been made for $N(q)=U_n(q)$ (we will write $O_{n,e}(q)$ instead of $O_e(q)$ in this case). Original Higman's conjecture was checked for $n\leq13$ in 2003 by A. Vera-Lopez and and J.M. Arregi in \cite{VeraLopezArregi03} by computer computations. In 2007, Isaacs proposed conjectural polynomials for $O_{n,e}(q)$ with $n\leq9$; he also computed explicitly $O_{n,\mu(n)}(q)$ and $O_{n,\mu(n)-1}(q)$, where $q^{\mu(n)}$ is the maximal possible degree of an irreducible character for $U_n(q)$. Note that these formulas were proved in 1997 by M. Marjoram in his unpublished thesis \cite{Marjoram97'},~\cite{Marjoram97}; the formula for $O_{n,\mu(n)-1}(q)$ also follows from A. Panov's classification of $2(\mu(n)-1)$-dimensional orbits presented in \cite{IgnatevPanov09}. In 2010, A. Evseev computed $O_{n,e}(q)$ for $n\leq13$ confirming Isaacs's formulas (and so Isaacs's conjecture). In 1999, Marjoram calculated $O_{n,1}(q)$ and $O_{n,2}(q)$ \cite{Marjoram99} (clearly, $O_{n,0}(q)=q^{n-1}$), while $O_{n,3}(q)$ was calculated by M. Loukaki in 2011 \cite{Loukaki11}. The results of Marjoram were independently rederived by T. Le in 2010 \cite{Le10}. E. Marberg verified Isaacs's conjecture for $e\leq8$ for arbitrary $n$ in 2011 \cite{Marberg11}.

For Sylow subgroups $N(q)$ in other Chevalley groups $G(q)$, the situation is as follows. In 1999, Marjoram computed $O_{\mu(\Phi)}(q)$ for orthogonal group (i.e., for $\Phi=B_n$ or $D_n$), where $q^{\mu(\Phi)}$ is the maximal possible degree of an irreducible character of $N(q)$, see \cite{Marjoram99}. For the symplectic case (i.e., for $\Phi=C_n$), the formula for $O_{\mu(\Phi)}(q)$ follows from the results of C.A.M. Andr\`e and A.-M. Neto on so-called supercharacters published in 2006 \cite{AndreNeto06}. In 2016, S.M. Goodwin, P. Mosch and G. R\"ohrle calculated $O_e(q)$ and proved Isaacs's conjecture for all possible $e$ and all finite Chevalley groups $G(q)$ of rank $\leq 8$, except $E_8$. We also would like to mention the papers \cite{GoodwinLeMagaardPaolini16} and \cite{LeMagaardPaolini20} of S.M. Goodwin, T. Le, K. Magaard and A. Paolini, where the case $F_4$ is considered without any restrictions on the characteristic, as well as the paper \cite{Szechtman06}, where F. Szechtman constructed certain explicit families of characters for symplectic and orthogonal cases.

Our classification of orbits of dimension $\leq6$ for classical root systems allows us to prove Isaacs's conjecture for all classical Chevalley groups. (Of course, for $\Phi=A_{n-1}$ we just rederived the results of Marjoram and Loukaki.) We consider the following theorem as the second main result of the paper.

\theop{\label{theo:Isaacs} Let $\Phi$ be a classical root system of rank $n$. Then\textup, for an arbitrary $n\geq1$\textup, $O_e(q)$ is a polynomial in $v=q-1$ with nonnegative integer coefficients for $0\leq e\leq3$. 
}{
Let us recall the classification of the coadjoint orbits of dimension 2 for $A_{n-1}$ based on Theorems~\ref{T:extdec},~\ref{Tdim2},~\ref{Tdim4},~\ref{Tdim6}. 

For every given orbit $N.f$ we have defined a Dynkin subdiagram $\Dyn(f)$ of $Dyn$. 
Thanks to Theorem~\ref{T:extdec} every orbit is a direct sum of extensive orbits for the respective subdiagrams of $\Dyn(f)$. The total dimension of $N.f$ shall be $2e\le6$ and thus we must have at most $e$ non-zero dimensional summands plus several 0-dimensional summands.

Thanks to Theorems~\ref{Tdim2},~\ref{Tdim4}, ~\ref{Tdim6} we have that $\Dyn(f)$ shall be a union of at most $e$ connected Dynkin diagrams of rank 8 or less with several disjoint points. 
We fix a certain such collection of diagrams.  
We can choose a place for such a collection in some nonegative amount of ways (it can be computed quite easily in every given case). 
For each such a placement Theorems~\ref{Tdim2},~\ref{Tdim4},~\ref{Tdim6} provides a list of representatives for the respective orbit. 
Every series of representatives is encoded by a collection of strings $\mathcal S_1, \mathcal S_2, \ldots$ (one string for one nontrivial connected component of $\Dyn(f)$) and it is easy to check that the amount of the representatives attached to $\mathcal S_1, \mathcal S_2, \ldots$ equals $v^{d_{\mathcal S_1}+d_{\mathcal S_2}+\ldots}$ where $d_{\mathcal S_i}$ is the number of letters 'S' in the string $\mathcal S_i$. 
We also can choose a 0-dimensional orbit for the remaining points in $v+1$ ways for each such a point.   
Thus we have $O_e(q)$ is equal to the sum of polymonomials of the form $v^*(v+1)^*$ with nonegative coefficients as needed.

}
\exam{i) We would like to evaluate explicitly $O_1(q)$ for $A_{n-1}$. We will check that $$O_1(q)=(v+1)^{n-4}(n(v^3+2v^2+v)-(3v^3+5v^2+2v))$$
for $n\ge3$. 

Indeed, for every given orbit $N.f$ we have defined a Dynkin subdiagram $\Dyn(f)$ of $A_{n-1}$. 
Thanks to Theorem~\ref{T:extdec} every orbit is a direct sum of extensive orbits for the respective subdiagrams of $\Dyn(f)$. The total dimension of $N.f$ shall be $2$ and thus we must have exactly one non-zero dimensional summand plus several 0-dimensional summands. 

Thanks to Theorems~\ref{Tdim2} we have that $\Dyn(f)$ shall be a union of $A_2$ or $A_3$ with several disjoint points. We can choose a place for $A_2$ in $n-2$ ways (plus $n-3$ disjoint points) and a place for $A_3$ in $n-3$ ways (plus $n-4$ disjoint points).  
Thus we have
\begin{equation}
\label{E:exfq}
\begin{split}
O_1(q)&=(n-2)(v+1)^{n-3}v+(n-3)(v+1)^{n-4}(v^3+v^2)\\
&=(v+1)^{n-4}(n(v^3+2v^2+v)-(3v^3+5v^2+2v))
\end{split}\end{equation}
for $n\ge3$. 
The left part of~(\ref{E:exfq}) is a polynomial in $v$ with nonnegative coefficients for every given $n\ge3$. 
If $n=2$ then $A_{n-1}$ consists of one point and the respective Lie algebra $\nt$ is abelian. Thus all its coadjoint orbits are zero-dimensional and therefore $O_1(q)=0$ for this case. Note that the answer provided by the formula~(\ref{E:exfq}) for $n=2$ equals $\frac{-v^2}{v+1}\ne0$.

ii) Arguing similarly, one can easily check that, for $B_n, n\ge 7,$
\begin{equation*}
\begin{split}
O_3(q)&=q^{n+3}+(2n-11)q^{n+2}+(n^2-12n+39)q^{n+1}+\dfrac{n^3-24n^2+185n-300}{6}q^n\\
&-\dfrac{n^3-15n^2+88n-176}{2}q^{n-1}+\dfrac{n^3-15n^2+76n-130}{2}q^{n-2}\\
&-\dfrac{n^3-21n^2+110n-174}{6}q^{n-3}-\dfrac{n^2-7n+12}{2}q^{n-4}\\
&=(v+1)^{n-4}\left(v^7+(2n-4)v^6+(n^2-6)v^5+\dfrac{n^3+6n^2+5n-60}{6}v^4\right.\\&+\left.\vphantom{\dfrac{n^3-37n+24}{6}v^4}\dfrac{n^3+9n^2-4n-42}{6}v^3+(n^2+n-5)v^2+(n-1)v\right).\\
\end{split}
\end{equation*}
}
\sect{Appendix: Lie--Dynkin nil-algebras}\label{SLDynL}\addcontentsline{toc}{subsection}{\ref{SLDynL}. Appendix: Lie--Dynkin nil-algebras}
We would like to extend Theorems~\ref{Tdim2},~\ref{Tdim4},~\ref{Tdim6}  and the description of orbits of dimension $2, 4, 6$ to the setup of Lie--Dynkin nil--algebras. They are natural infinite-dimensional analogues of nilradicals of Borel subalgebras of simple finite-dimensional Lie algebras.

Lie--Dynkin nil-algebras can be described by the formulas of Section~\ref{Sfdcn} after a replacement of $\{1,\ldots, n\}$ by an arbitrary countable set $S$ and the binary relation $<$ on $\{1,\ldots, n\}$ by an arbitrary linear order relation $\prec$ on this set, see~\cite{IgnatyevPetukhov21}. 
Thus for a given pair $(S, \prec)$ we define four different simple Lie algebras: $A_S(\prec), B_S(\prec), C_S(\prec), D_S(\prec)$ together with the respective set of positive roots $\Phi_A^+(S, \succ), \Phi_B^+(S, \succ), \Phi_C^+(S, \succ), \Phi_D^+(S, \succ)$. We will often omit $S$ in this notation writing $A(\prec)$, $B(\prec)$, $C(\prec)$, $D(\prec)$ and $\Phi_A(\succ), \Phi_B(\succ), \Phi_C(\succ), \Phi_D(\succ)$.

We denote by $\frak g(S)$ one of the Lie algebras defined above. By $\nt(S)$ we denote the respective Lie--Dynkin nil-algebra (which will be a subalgebra of $\frak g(S)$). 
Further, for a subset $S'$ of $S$ we denote by $\frak g(S')$ and $\nt(S')$ the respective subalgebras of $\frak g(S)$. 

In the setting of a finite-dimensional Lie algebra it is well known that the simple roots $\Pi$ of $\Phi^+$ are exactly the indecomposable vectors of $\Phi^+$, and that they generate $\Phi^+$ as a semigroup. 
This should provide a motivation for the following notation: 
\begin{center}$\Pi_{A}(\succ):=\{\text{the~indecomposable~elements~of~}\Phi^+_A(\succ)\}$ (the same for $B/C/D$ cases).\end{center}
Note that $\Pi_*(\succ)$ can be even empty in this setting (for example this will be the case if $(S, \succ)\cong (\mathbb Q, >)$) and thus it does not generate $\Phi^+_*(\succ)$ at least in some cases.

It is clear that an exhaustion of $S$ by finite sets $S'$ gives rise to the exhaustion of $\frak n(S)$ by finite-dimensional nilpotent Lie algebras $\frak n(S')$, and thus $\frak n(S)$ is locally nilpotent, see~\cite{IgnatyevPetukhov21}. 

This allows us to define the adjoint group $N(S)$ of $\nt(S)$ and thus to consider the adjoint (and coadjoint) orbits with respect to the action of $N(S)$ on $\nt(S)$ and $\nt(S)^*$. 
For $f\in\nt(S)^*$ we denote by $\mathcal O(f)$ the respective coadjoint orbit. 

Consider $f\in\nt(S)^*$. We define $\Supp{f}, \NSupp{f}, \Dyn(f), \Phi(f)$ by the same formulas as in the finite-dimensional case. 
Note that $\Dyn(f)$ can be empty and we accept this answer. 
From now one we assume that $S'$ denotes a finite set.
\lemmp{We have
$$\Supp{f}=\bigcup_{S'\subset S} \NSupp{f|_{\frak g(S')}}$$
$$\NSupp{f}=\bigcup_{S'\subset S} \NSupp{f|_{\frak g(S')}},$$
$$\Phi(f)=\bigcup_{S'\subset S} (\Phi(f|_{\frak g(S')})\cap\underline{\Phi}_*(\succ)).$$
}{
All statements follows easily from the respective definitions.
}
We denote by $\Dyn_{\bar 1}(f)$ the subgraph of $\Dyn(f)$ from which we eliminated all the connected components consisting of a single vertex with no edges.  
\propp{\label{P:crfd}We have $\dim \mathcal O(f)<\infty$ if and only if $|\Dyn_{\bar 1}(f)|<\infty$ and $\Supp{f}\subset\Phi(f)$.}{
This follows from Proposition~\ref{PlWd}.
}
\lemmp{\label{L:fincr}If $\dim \mathcal O(f)<\infty$ then $\mathcal O(f)\cong N(S').f|_{\nt(S')}$ as algebraic varieties for any finite subset $S'$ such that $\Dyn_{\bar 1}(f)\subset\Phi(S')$ where we consider $\Dyn_{\bar 1}(f)$ as a set of the respective simple roots in this formula.}{
This is implied by a combination of~\cite[Lemma~3.11]{IgnatyevPetukhov21}, Theorem~\ref{T:extdec} and the fact that $\Dyn(f|_{\gt(S')})$ is $\Dyn(f)$ plus several disjoint points. 
}
Proposition~\ref{P:crfd} and Lemma~\ref{L:fincr} considered together reduce the description of the finite-dimensional coadjoint orbits of nil-Dynkin Lie algebras to the finite-dimensional case.


\sect{Appendix: Enumeration of roots by CHEVIE package}\label{S:chv}\addcontentsline{toc}{subsection}{\ref{S:chv}. Appendix: Enumeration of roots by CHEVIE package}
Throughout the paper we refer to the enumeration of roots by CHEVIE package for the classical Dynkin diagrams. 
For the convenience of the reader we provide tables in notation of Section~\ref{Sfdcn} giving this enumeration of roots for the classical root systems needed in the classification of extensive orbits of dimensions 2, 4, 6.      
$$A_2:
~{\Autonumfalse\mymatrix{
\lNote{2}\Note{1}\Rt{2pt}2& \Note{2}\Bot{2pt}\pho\\
\lNote{3}3&1 \Rt{2pt}\\
}},\hspace{10pt}
A_3:
~{\Autonumfalse\mymatrix{
\lNote{2}\Note{1}\Rt{2pt}3& \Note{2}\Bot{2pt}\pho&\Note{3}\pho\\
\lNote{3}5&2 \Rt{2pt}&\Bot{2pt}\\
\lNote{4}6&4&1 \Rt{2pt}\\
}},\hspace{10pt}
A_4:
~{\Autonumfalse\mymatrix{
\lNote{2}\Note{1}\Rt{2pt}4& \Note{2}\Bot{2pt}\pho&\Note{3}\pho&\Note{4}\pho\\
\lNote{3}7&3 \Rt{2pt}&\Bot{2pt}&\\
\lNote{4}9&6&2 \Rt{2pt}&\Bot{2pt}\\
\lNote{5}10&8&5 &1\Rt{2pt}\\
}},\hspace{10pt}
A_5:
~{\Autonumfalse\mymatrix{
\lNote{2}\Note{1}\Rt{2pt}5& \Note{2}\Bot{2pt}\pho&\Note{3}\pho&\Note{4}\pho&\Note{5}\pho\\
\lNote{3}9&4 \Rt{2pt}&\Bot{2pt}&&\\
\lNote{4}12&8&3 \Rt{2pt}&\Bot{2pt}&\\
\lNote{5}14&11&7 &2\Rt{2pt}&\Bot{2pt}\\
\lNote{6}15&13&10 &6&1\Rt{2pt}\\
}},\hspace{10pt}$$
$$A6:
~{\Autonumfalse\mymatrix{
\lNote{2}\Note{1}\Rt{2pt}6& \Note{2}\Bot{2pt}\pho&\Note{3}\pho&\Note{4}\pho&\Note{5}\pho&\Note{6}\pho\\
\lNote{3}11&5 \Rt{2pt}&\Bot{2pt}&&&\\
\lNote{4}15&10&4 \Rt{2pt}&\Bot{2pt}&&\\
\lNote{5}18&14&9 &3\Rt{2pt}&\Bot{2pt}&\\
\lNote{6}20&17&13 &8&2\Rt{2pt}&\Bot{2pt}\\
\lNote{7}21&19&16 &12&7&1\Rt{2pt}\\
}},
\hspace{10pt}
A_7:
~{\Autonumfalse\mymatrix{
\lNote{2}\Note{1}\Rt{2pt}7& \Note{2}\Bot{2pt}\pho&\Note{3}\pho&\Note{4}\pho&\Note{5}\pho&\Note{6}\pho&\Note{7}\pho\\
\lNote{3}13&6 \Rt{2pt}&\Bot{2pt}&&&&\\
\lNote{4}18&2&5 \Rt{2pt}&\Bot{2pt}&&&\\
\lNote{5}22&17&11 &4\Rt{2pt}&\Bot{2pt}&&\\
\lNote{6}25&21&16 &10&3\Rt{2pt}&\Bot{2pt}&\\
\lNote{7}27&24&20 &15&9&2\Rt{2pt}&\Bot{2pt}\\
\lNote{8}28&26&23 &19&14&8&1\Rt{2pt}\\
}},
$$
$$\begin{array}{c}B_2\\=\\C_2\end{array}:
~{\Autonumfalse\mymatrix{
\Note{1}\Rt{2pt}2& \Note{2}\Bot{2pt}\pho\\
3&1 \Rt{2pt}\Bot{2pt}\\
4\Rt{2pt}& \\
}},\hspace{10pt}
\begin{array}{c}B_3\\\sim\\ C_3\end{array}:
~{\Autonumfalse\mymatrix{
\Note{1}\Rt{2pt}3& \Note{2}\Bot{2pt}\pho& \Note{3}\pho\\
5&2 \Rt{2pt}&\Bot{2pt}\pho\\
7&4&\Rt{2pt}\Bot{2pt}1\\
8& \Rt{2pt}\Bot{2pt}6& \pho\\
\Rt{2pt}9& \pho& \pho\\
}},\hspace{10pt}
\begin{array}{c}B_4\\\sim\\ C_4\end{array}:
~{\Autonumfalse\mymatrix{
4 \Note{1}\Rt{2pt}&
\Note{2}\Bot{2pt}\pho& \Note{3}\pho&\Note{4}\pho\\
7&3 \Rt{2pt}&\Bot{2pt}\pho&\\
10&6&2\Rt{2pt}&\Bot{2pt}\\
12&9&5&1\Rt{2pt}\Bot{2pt}\\
14&11&8\Rt{2pt}\Bot{2pt}&\\
15&13 \Rt{2pt}\Bot{2pt}& \pho&\\
16\Rt{2pt}& \pho& \pho&\\
}},\hspace{10pt}
\begin{array}{c}B_5\\\sim\\ C_5\end{array}:
~{\Autonumfalse\mymatrix{
5 \Note{1}\Rt{2pt}&
\Note{2}\Bot{2pt}\pho& \Note{3}\pho&\Note{4}\pho&\Note{5}\pho\\
9&4 \Rt{2pt}&\Bot{2pt}\pho&&\\
13&8&3\Rt{2pt}&\Bot{2pt}&\\
16&12&7&2\Rt{2pt}&\Bot{2pt}\\
19&15&11&6&1\Rt{2pt}\Bot{2pt}\\
21&18&14&10\Rt{2pt}\Bot{2pt}&\\
23&20&17\Rt{2pt}\Bot{2pt}&&\\
24&22 \Note{2}\Rt{2pt}\Bot{2pt}&\pho&&\\
25\Rt{2pt}&\pho&\pho&&\\
}},$$
$$\begin{array}{c}B_6\\\sim\\ C_6\end{array}:
~{\Autonumfalse\mymatrix{
6 \Note{1}\Rt{2pt}&
\Note{2}\Bot{2pt}\pho& \Note{3}\pho&\Note{4}\pho&\Note{5}\pho&\Note{6}\pho\\
11&5 \Rt{2pt}&\Bot{2pt}\pho&&&\\
16&10&4\Rt{2pt}&\Bot{2pt}&&\\
20&15&9&3\Rt{2pt}&\Bot{2pt}&\\
24&19&14&8&2\Rt{2pt}&\Bot{2pt}\\
27&23&18&13&7&1\Rt{2pt}\Bot{2pt}\\
30&26&22&17&12\Rt{2pt}\Bot{2pt}&\\
32&29&25&21\Rt{2pt}\Bot{2pt}&&\\
34&31&28\Rt{2pt}\Bot{2pt}&&&\\
35&33 \Note{2}\Rt{2pt}\Bot{2pt}&\pho&&&\\
36\Rt{2pt}&\pho&\pho&&&\\
}},\hspace{10pt}
\begin{array}{c}B_7\\\sim\\ C_7\end{array}:
~{\Autonumfalse\mymatrix{
7 \Note{1}\Rt{2pt}&
\Note{2}\Bot{2pt}\pho& \Note{3}\pho&\Note{4}\pho&\Note{5}\pho&\Note{6}\pho&\Note{7}\pho\\
13&6 \Rt{2pt}&\Bot{2pt}\pho&&&&\\
19&12&5\Rt{2pt}&\Bot{2pt}&&&\\
24&18&11&4\Rt{2pt}&\Bot{2pt}&&\\
29&23&17&10&3\Rt{2pt}&\Bot{2pt}&\\
33&28&22&16&9&2\Rt{2pt}&\Bot{2pt}\\
37&32&27&21&15&8&1\Rt{2pt}\Bot{2pt}\\
40&36&31&26&20&14\Rt{2pt}\Bot{2pt}&\\
43&39&35&30&25\Rt{2pt}\Bot{2pt}&&\\
45&42&38&34\Rt{2pt}\Bot{2pt}&&&\\
47&44&41\Rt{2pt}\Bot{2pt}&&&&\\
48&46 \Note{2}\Rt{2pt}\Bot{2pt}&\pho&&&&\\
49\Rt{2pt}&\pho&\pho&&&&\\
}}.
$$

$$D_4:
~{\Autonumfalse\mymatrix{
4 \Note{1}\Rt{2pt}&
\Note{2}\Bot{2pt}\pho& \Note{3}\pho\\
7&3 \Rt{2pt}&\Bot{2pt}\pho\\
10&6&2\Rt{2pt}\\
9&5&1\Rt{2pt}\Bot{2pt}\\
11&8 \Rt{2pt}\Bot{2pt}&\\
12\Rt{2pt}& \pho& \pho\\
}},\hspace{10pt}
D_5:
~{\Autonumfalse\mymatrix{
5 \Note{1}\Rt{2pt}&
\Note{2}\Bot{2pt}\pho& \Note{3}\pho&\Note{4}\pho\\
9&4 \Rt{2pt}&\Bot{2pt}\pho&\\
13&8&3\Rt{2pt}&\Bot{2pt}\\
16&12&7&2\Rt{2pt}\\
15&11&6&1\Rt{2pt}\Bot{2pt}\\
18&14&10\Rt{2pt}\Bot{2pt}&\\
19&17 \Note{2}\Rt{2pt}\Bot{2pt}&\pho&\\
20\Rt{2pt}&\pho&\pho&\\
}},\hspace{10pt}
D_6:
~{\Autonumfalse\mymatrix{
6 \Note{1}\Rt{2pt}&
\Note{2}\Bot{2pt}\pho& \Note{3}\pho&\Note{4}\pho&\Note{5}\pho\\
11&5 \Rt{2pt}&\Bot{2pt}\pho&&\\
16&10&4\Rt{2pt}&\Bot{2pt}&\\
20&15&9&3\Rt{2pt}&\Bot{2pt}\\
24&19&14&8&2\Rt{2pt}\\
23&18&13&7&1\Rt{2pt}\Bot{2pt}\\
26&22&17&12\Rt{2pt}\Bot{2pt}&\\
28&25&21\Rt{2pt}\Bot{2pt}&&\\
29&27 \Note{2}\Rt{2pt}\Bot{2pt}&\pho&&\\
30\Rt{2pt}&\pho&\pho&&\\
}},\hspace{10pt}
$$
$$
D_7:
~{\Autonumfalse\mymatrix{
7 \Note{1}\Rt{2pt}&
\Note{2}\Bot{2pt}\pho& \Note{3}\pho&\Note{4}\pho&\Note{5}\pho&\Note{6}\pho\\
13&6 \Rt{2pt}&\Bot{2pt}\pho&&&\\
19&12&5\Rt{2pt}&\Bot{2pt}&&\\
24&18&11&4\Rt{2pt}&\Bot{2pt}&\\
29&23&17&10&3\Rt{2pt}&\Bot{2pt}\\
33&28&22&16&9&2\Rt{2pt}\\
32&27&21&15&8&1\Rt{2pt}\Bot{2pt}\\
36&31&26&20&14\Rt{2pt}\Bot{2pt}&\\
38&35&30&25\Rt{2pt}\Bot{2pt}&&\\
40&37&34\Rt{2pt}\Bot{2pt}&&&\\
41&39 \Note{2}\Rt{2pt}\Bot{2pt}&\pho&&&\\
42\Rt{2pt}&\pho&\pho&&&\\
}},\hspace{10pt}
D_8:
~{\Autonumfalse\mymatrix{
8 \Note{1}\Rt{2pt}&
\Note{2}\Bot{2pt}\pho& \Note{3}\pho&\Note{4}\pho&\Note{5}\pho&\Note{6}\pho&\Note{7}\pho\\
15&7 \Rt{2pt}&\Bot{2pt}\pho&&&&\\
22&14&6\Rt{2pt}&\Bot{2pt}&&&\\
28&21&13&5\Rt{2pt}&\Bot{2pt}&&\\
34&27&20&12&4\Rt{2pt}&\Bot{2pt}&\\
39&33&26&19&11&3\Rt{2pt}&\Bot{2pt}\\
44&38&32&25&18&10&2\Rt{2pt}\\
43&37&31&24&17&9&1\Rt{2pt}\Bot{2pt}\\
47&42&36&30&23&16\Rt{2pt}\Bot{2pt}&\\
50&46&41&35&29\Rt{2pt}\Bot{2pt}&&\\
52&49&45&40\Rt{2pt}\Bot{2pt}&&&\\
54&51&48\Rt{2pt}\Bot{2pt}&&&&\\
55&43 \Note{2}\Rt{2pt}\Bot{2pt}&\pho&&&&\\
56\Rt{2pt}&\pho&\pho&&&&\\
}}.
$$

\medskip\textsc{Mikhail Ignatev: National Research University Higher School of Economics,\break\indent Pokrovsky Boulevard 11, 109028, Moscow, Russia}

\emph{E-mail address}: \texttt{mihail.ignatev@gmail.com}

\medskip\textsc{Alexey Petukhov: Institute for Information Transmission Problems, Bolshoy\break\indent Karetniy 19--1, 127994, Moscow, Russia}

\emph{E-mail address}: \texttt{alex--2@yandex.ru}

\newpage
\end{document}